\documentclass[a4paper]{amsart} 
\usepackage[ascii]{inputenc} 

\usepackage{amssymb}
\usepackage{mathrsfs}
\usepackage[hidelinks]{hyperref}

\usepackage{xcolor}

\usepackage[shortlabels]{enumitem}
\setlist[enumerate,1]{label={(\Alph*)}}
\setlist[enumerate,2]{label={(\alph*)}}
\setlist[enumerate,3]{label={$\bullet_{\arabic*}$}}

\newenvironment{PROOF}[2][\proofname.]
   {\begin{proof}[#1]\renewcommand{\qedsymbol}{\textsquare$_{\rm #2}$}}
   {\end{proof}}

\newtheorem{theorem}{Theorem}[section] 
\newtheorem{claim}[theorem]{Claim}
\newtheorem{lemma}[theorem]{Lemma} 
\newtheorem{construction}[theorem]{Construction} 
 
\newtheorem{conclusion}[theorem]{Conclusion}
\newtheorem{observation}[theorem]{Observation}

\theoremstyle{definition}
\newtheorem{definition}[theorem]{Definition}

\newtheorem{fact}[theorem]{Fact}

\newtheorem{discussion}[theorem]{Discussion}
\newtheorem{convention}[theorem]{Convention}

\newtheorem{hypothesis}[theorem]{Hypothesis}

\theoremstyle{remark}
\newtheorem{remark}[theorem]{Remark}

\newtheorem{notation}[theorem]{Notation}

\DeclareMathOperator{\divides}{\big|}

\newcommand{\pinkappa}{\kappa}

\newcommand{\x}{\mathrm{x}}

\newcommand{\sseq}{\mathrm{sseq}}

\newcommand{\ba}{\mathsf{ba}}
\newcommand{\BA}{\mathsf{BA}}
\newcommand{\MA}{\mathsf{MA}}
\newcommand{\LO}{\mathrm{LO}}

\newcommand{\lev}{\mathrm{lev}}

\newcommand{\dcI}{{\dot{\mathcal I}}}
\newcommand{\isoI}{{\dot{\mathbb I}}}
\newcommand{\numbIE}{{\dot{I}\dot{E}}}

\newcommand{\ptr}{\mathrm{ptr}}
\newcommand{\rt}{\mathrm{rt}}
\newcommand{\trr}{\mathrm{trr}}
\newcommand{\nr}{\mathrm{nr}}
\newcommand{\vr}{\mathrm{vr}}

\newcommand{\Eq}{\mathrm{Eq}}
\newcommand{\Res}{\mathrm{Res}}
\newcommand{\Suc}{\mathrm{Suc}}
\newcommand{\Sur}{\mathrm{Sur}}

\newcommand{\cd}{\mathrm{cd}}

\newcommand{\leqdot}{\mathrel{\mathord{\leq}\!\!\raise
0.8 pt\hbox{$\cdot$}}}

\newcommand{\closure}{\mathrm{closure}}

\newcommand{\tthh}{\mathrm{th}}

\newcommand{\cf}{\mathrm{cf}}

\newcommand{\dom}{\mathrm{dom}}

\newcommand{\id}{\mathrm{id}}
\newcommand{\inc}{\mathrm{inc}}

\newcommand{\otp}{\mathrm{otp}}

\newcommand{\rang}{\mathrm{rang}}

\newcommand{\supminus}{{\text{\hspace{.05cm}--}}}

\newcommand{\qf}{\mathrm{qf}}

\newcommand{\Card}{\mathrm{Card}}

\newcommand{\lex}{\mathrm{lex}}

\newcommand{\set}{\mathrm{set}}

\newcommand{\tp}{\mathrm{tp}}

\newcommand{\tr}{\mathrm{tr}}

\newcommand{\base}{\mathrm{base}}

\newcommand{\bfB}{\mathbf{B}}

\newcommand{\bfc}{\mathbf{c}}

\newcommand{\bfE}{\mathbf{E}}
\newcommand{\bfe}{\mathbf{e}}

\newcommand{\bff}{\mathbf{f}}

\newcommand{\bfJ}{\mathbf{J}}
\newcommand{\bfj}{\mathbf{j}}

\newcommand{\bfn}{\mathbf{n}}

\newcommand{\bfR}{\mathbf{R}}

\newcommand{\bfS}{\mathbf{S}}

\newcommand{\bbL}{\mathbb{L}}
\newcommand{\bbM}{\mathbb{M}}

\newcommand{\bbP}{\mathbb{P}}

\newcommand{\bbZ}{\mathbb{Z}}

\newcommand{\mn}{\medskip\noindent}
\newcommand{\sn}{\smallskip\noindent}
\newcommand{\bn}{\bigskip\noindent}

\newcommand{\cM}{\mathscr{M}}

\newcommand{\cP}{\mathscr{P}}

\newcommand{\clB}{\mathcal{B}}

\newcommand{\clE}{\mathcal{E}}

\newcommand{\clH}{\mathcal{H}}

\newcommand{\clP}{\mathcal{P}}

\newcommand{\clS}{\mathcal{S}}
\newcommand{\clT}{\mathcal{T}}
\newcommand{\clU}{\mathcal{U}}

\newcommand{\gA}{\mathfrak{A}}
\newcommand{\gB}{\mathfrak{B}}

\newcommand{\eps}{\varepsilon}
\newcommand{\cl}{c\kern-.11ex \ell}
\newcommand{\lh}{{\ell\kern-.27ex g}}
\newcommand{\rest}{\restriction}

\newcommand{\caret}{{\char 94}}
\newcommand{\LL}{\langle}
\newcommand{\RR}{\rangle}

\newcommand{\subref}[1]{$_{\mathrm{\texttt{=}}\mathsf{L{#1}}}$}
\newcommand{\lepref}[1]{({<}\,{#1})}

\newcommand{\olsi}[1]{\,\overline{\!{#1}}} 

\newcommand*{\defeq}{\mathrel{\vcenter{\baselineskip0.5ex \lineskiplimit0pt\hbox{\scriptsize.}\hbox{\scriptsize.}}}=}

\usepackage{tcolorbox}

\newcount\skewfactor
\def\mathunderaccent#1#2 {\let\theaccent#1\skewfactor#2
\mathpalette\putaccentunder}
\def\putaccentunder#1#2{\oalign{$#1#2$\crcr\hidewidth
\vbox to.2ex{\hbox{$#1\skew\skewfactor\theaccent{}$}\vss}\hidewidth}}
\def\name{\mathunderaccent\tilde-3 }
\def\Name{\mathunderaccent\widetilde-3 }

\newbox\noforkbox \newdimen\forklinewidth
\setbox0\hbox{$\textstyle\bigcup$}
\setbox1\hbox to \wd0{\hfil\vrule width \forklinewidth depth \dp0
                        height \ht0 \hfil}
\setbox\noforkbox\hbox{\box1\box0\relax}
\def\unionstick{\mathop{\copy\noforkbox}\limits}
\def\nonfork#1#2_#3{#1\unionstick_{\textstyle #3}#2}
\def\nonforkin#1#2_#3^#4{#1\unionstick_{\textstyle #3}^{\textstyle
    #4}#2}
\setbox0\hbox{$\textstyle\bigcup$}
\setbox1\hbox to \wd0{\hfil{\sl /\/}\hfil}
\setbox2\hbox to \wd0{\hfil\vrule height \ht0 depth \dp0 width
                                \forklinewidth\hfil}
\newbox\doesforkbox
\setbox\doesforkbox\hbox{\box1\box0\relax}
\def\nunionstick{\mathop{\copy\doesforkbox}\limits}

\def\fork#1#2_#3{#1\nunionstick_{\textstyle #3}#2}
\def\forkin#1#2_#3^#4{#1\nunionstick_{\textstyle #3}^{\textstyle
    #4}#2}

\newcommand{\stickT}{%
\setbox255=\hbox{\raise1ex\hbox{$\hspace{0.2pt}\,\bullet\,$}}
\mathord{\rlap{\hbox to\wd255{\hss\hbox{$|$}\hss}}
\box255}
}
\newcommand{\stickS}{%
\setbox255=\hbox{\raise0.6ex\hbox{$\scriptstyle\bullet$}}
\mathord{\rlap{\hbox to\wd255{\hss\hbox{$\scriptstyle|$}\hss}}
\box255}
}

\title[Complicated index models and Boolean algebras]
{Building complicated index models and Boolean algebras}

\author {Saharon Shelah}
\address{Einstein Institute of Mathematics\\
Edmond J. Safra Campus, Givat Ram\\
The Hebrew University of Jerusalem\\
Jerusalem, 91904, Israel\\
 and \\
 Department of Mathematics\\
 Hill Center - Busch Campus \\ 
 Rutgers, The State University of New Jersey \\
 110 Frelinghuysen Road \\
 Piscataway, NJ 08854-8019 USA}
\email{shelah@math.huji.ac.il}
\urladdr{http://shelah.logic.at}

\thanks{For versions up to 2019, the author thanks Alice Leonhardt for the beautiful typing. In the latest version, the author thanks Craig Falls for generously funding typing services, and thanks Matt Grimes for the careful and beautiful typing. For their partial support of this research, the author would like to thank the Israel Science Foundation (=ISF) with grant 2320/23 (2023-2027) 
and grant 1838(19)  
(2019-1023).
and for various grants from the BSF (United States Israel Binational Foundation), the Israel Academy of Sciences and the NSF via Rutgers University.\\
This paper is number 511 in the author's publication list; it has existed 
(and been occasionally revised) for many years. It was mostly ready in the early nineties, and was made public to some extent. This was written as Chapter VII of 
the book \cite{Sh:e}, which hopefully will materialize some day, but in the 
meantime it is \cite{Sh:511}. The intention was to have \cite{Sh:E58} 
(revising \cite{Sh:229}) for Ch.I, \cite{Sh:421} for Ch.II, \cite{Sh:E59} for 
Ch.III, \cite{Sh:309} for Ch.IV, \cite{Sh:363} for Ch.V, \cite{Sh:331} for Ch.VI,
\cite{Sh:511} for Ch.VII, \cite{Sh:E60} (a revision of \cite{Sh:128}) for Ch.VIII,
\cite{Sh:E62} for the appendix, and probably \cite{Sh:757}, \cite{Sh:384}, \cite{Sh:482}, and \cite{Sh:800}. References like [Sh:511, 0.1\subref{x2}]
means that \textsf{x2} is the internal label of Lemma \ref{x2} in the TeXfile of [Sh:511].\\
The reader should note that the version in my website is usually more
up-to-date than the one in the mathematical archive.}

\makeatletter
\@namedef{subjclassname@2020}{\textup{2020} Mathematics Subject Classification}
\makeatother


\date {June 9, 2026} 

\begin{document}
\makeatletter\def\shfiuwefootnote{\gdef\@thefnmark{}\@footnotetext}\makeatother\shfiuwefootnote{Version 2026-06-09\_3. See \url{https://shelah.logic.at/papers/511/} for possible updates.}
\begin{abstract}
We build models using an indiscernible model sub-structures of
${}^{\kappa \ge}\lambda$ and related more complicated structures.  We
use this to build various Boolean algebras.
\end{abstract}

\maketitle
\numberwithin{equation}{section}
\setcounter{section}{-1}
\newpage

\section {Introduction}

We continue \cite{Sh:E59}, \cite{Sh:309}, and \cite{Sh:311} (improving 
\cite[III]{Sh:300}) on the one hand, and \cite{Sh:136} on the other. 
A starting idea was that the ``many pairwise non-isomorphic models'' 
proofs in Chapters VII and VIII of \cite{Sh:a}, \cite{Sh:c} (and earlier
\cite{Sh:12}, \cite{Sh:31}, \cite{Sh:51}) can be generalized to many contexts --- 
in particular, to building Boolean algebras (as in \cite{Sh:51}, \cite{Sh:136}).

In \cite{Sh:E59}, \cite{Sh:311} we build the so-called ``strongly unembeddable sequence of index models'' $\LL I_\alpha : \alpha < \lambda\RR$, and from there build `many models' or `models with few automorphisms' (or endomorphisms: e.g. for abelian groups and --- our central point here --- Boolean algebras) as was done earlier in \cite{Sh:136}.

The index models were mainly linear orders and trees with $\omega+1$ levels. In this paper, we deal with generalizations. (See also \cite{Sh:897}.)
 
We begin with an example that motivates our need to pass beyond the
framework of trees with $\omega+1$ levels. Suppose that we are asked to
construct a rigid Boolean algebra of cardinality $\lambda$. We can take a
sequence $\LL I_\alpha : \alpha < \lambda\RR$ exemplifying that
$K^\omega_\tr$ has the so-called full $(\lambda,\lambda,\aleph_0,\aleph_0)$-bigness 
property (see \cite[2.5\subref{2.3}]{Sh:E59}). (It says that each $I_\alpha$ is 
so-called ``strongly unembeddable" into 
$\sum \big\{I_\beta : \beta \in \lambda \setminus \{\alpha\} \big\}$.
These exist: e.g. $\lambda$ is regular and $I_\alpha$ codes $S_\alpha$, 
a stationary subset of\\ $\{\delta < \lambda : \cf(\delta) = \aleph_0\}$, 
with $\LL S_\alpha : \alpha < \lambda\RR$ pairwise disjoint.)

Now build a Boolean algebra
$\BA(I_\alpha)$ for each $\alpha$. We then construct a rigid
Boolean algebra $\bfB_\lambda$ by choosing an increasing continuous sequence
$\LL \bfB_\alpha:\alpha\leq\lambda\RR$, where $\bfB_0$ is trivial
and $\bfB_{\alpha+1}$ is obtained from $\bfB_\alpha$ by ``planting'' a copy
of  $\BA(I_\alpha)$ below $a_\alpha \in \bfB_\alpha$, and our bookkeeping 
will ensure that
$\bfB_\lambda \setminus \{0\} = \{a_\alpha : \alpha < \lambda\}$. 
This seems to be a reasonable strategy, and it works (see a little more below). 
Now suppose, moreover, that we are asked to construct a complete Boolean algebra
$\bfB$ of cardinality $\lambda$ with no non-trivial one-to-one endomorphism. 
We should assume that $\lambda^{\aleph_0}=\lambda$ (as the cardinality  
of any complete Boolean algebra satisfies this equality) and it is
natural to demand in
addition that $\bfB$ satisfies the ccc. It is not hard to modify the
construction above so that $\bfB_\lambda$ has the ccc, so let
$\bfB$ be its completion. 

Assume toward a contradiction that
$f : \bfB \to \bfB$ is a non-trivial, one-to-one endomorphism. We
can find $a\in \bfB\setminus\{0\}$ with $a\cap f(a)=0$ and $\alpha<\lambda$
such that $a=a_\alpha$. Then
$I_\alpha$ is embedded in $\bfB\rest a_\alpha$ in some sense, say by
$\eta \mapsto a^\alpha_\eta$. Hence $\eta\mapsto f(a^\alpha_\eta)$ is a similar
embedding into $\bfB\rest f(a)$ that is constructed from $\LL
I_\beta : \beta\neq\alpha\RR$ alone. It seems reasonable that the {demand} 
``$I_\alpha$ strongly unembeddable into $\sum \{I_\beta:\beta\neq \alpha\}$"
in the sense of Definition \cite[2.5\subref{2.3}]{Sh:E59} {can be used} 
to deduce a contradiction; this works in the case above (i.e. without the
completion demand). However in the present case $f(a^\alpha_\eta)$ 
is not in general a member of $\bfB_\lambda$, but rather is a countable
union $\bigcup\limits_{n<\omega} b^\alpha_{\eta,n}$ of members of 
$\bfB_\lambda$.  We would like to find an appropriate unembeddability
condition of $I_\alpha$ into $\sum\limits_{\beta\neq\alpha} I_\beta$ to
handle this complication. At some price, our original notion can be modified
to handle this complication when $\eta$ has finite length, but not when
$\eta$ has length $\omega$. Instead, in this latter case, we replace it by an
``approximation'' $b^\alpha_{\eta, n(\alpha,\eta)}>0$: this was part of the 
motivation of having the definition ``strongly finitary on $P^I_\omega$'' in
\cite[2.5\subref{2.3}]{Sh:E59}. Previously, we could use demands like
``$a^\alpha_{\eta\rest\ell} \geq a^\alpha_\eta$'' but now we have to use 
demands like $a^\alpha_\nu \cap a^\alpha_\eta=0$, $\lh(\eta) = \omega$, 
$\lh(\nu) < \omega$, but such demands tend to contradict the ccc.

Our solution is to replace subtrees of ${}^{\omega{\geq}}\lambda$ 
by index sets $I$ of the form
\[
I = I' \cup \big\{(\eta\rest n) \caret \LL\alpha_\ell\RR:
n<\omega,\ \eta\in I',\ \eta(n) = (\alpha_0,\alpha_1), \text{ and } \ell
\in \{0,1\}\big\},
\]
where $I'\subseteq {}^\omega\big\{ (\alpha_0,\alpha_1) : \alpha_0 < \alpha_1 < \lambda \big\}$, and choose $\BA(I)$ to be generated by
$\{a^I_\eta : \eta\in I\}$ freely except that
\[
\eta\in I'\ \wedge\ \eta(n)=(\alpha_0,\alpha_1)\ \Rightarrow\ 
a^I_{\eta\rest n \caret \LL\alpha_0\RR}-a^I_{\eta \rest
n \caret \LL\alpha_1\RR} \ge a^I_\eta.
\]

\sn
(Actually, to ensure the ccc it is better to use a more complicated
variant.) But now the bigness properties have to be proved in this
context. For other aims, we use subtrees of ${}^{\omega\geq}2$ of cardinality
$\kappa \in [\aleph_1, 2^{\aleph_0}),$ originally to deal with number
of non-isomorphic models.

In this work we deal with more complicated index sets as motivated above.

In \S1 we introduce classes like $K^\omega_{\tr(n)}$, which are close
to being trees with $\omega+1$ levels, together with bigness
properties (related to $\psi_{\tr(n)}$) for them. We then generalize this to $K_{\tr(n)}^\partial$.

In \S1(B), we prove some existence theorems of the form
``for many $\lambda$ there is a sequence $\LL I_\alpha : \alpha < \lambda\RR$, 
where each $I_\alpha\in K^\omega_{\tr(n)}$ has cardinality
$\lambda$ and is strongly $\psi_{\tr(n)}$-unembeddable into
$\sum\limits_{\beta\neq\alpha} I_\beta$.'' We also define ``super'' versions
of these bigness properties related to the ones in
\cite[1.1\subref{7.1},1.5\subref{7.3}]{Sh:331}. Note that Lemma \ref{1.12} is used by Asgharzadeh, Golshani, and the author in \cite{Sh:1232}.
weaken
In \S2 we construct Boolean algebras 
with few appropriate morphisms for several versions. 

In \S3 we construct a ccc Boolean algebra of cardinality $2^{\aleph_0}$ of
pre-given length (see Definition \ref{5.1}) such that any infinite homomorphic
image has cardinality $2^{\aleph_0}$. We use a Boolean algebra constructed 
from a single $I \in K^\omega_{\tr(h)}$ as in \S 2. As it happens, the 
complicated $I \in K^\omega_{\tr(h)}$ are not needed, just non-trivial ones. 
Our point is that $K^\omega_{\tr(h)}$ is not good just for the constructions 
in \S2, it is a quite versatile way to build structures with pre-assumed 
properties (not to speak of varying the index model).

The main result is (\ref{5.4}):
\mn
\begin{enumerate}
    \item[$(*)$]  For $\mu\in [\aleph_0,2^{\aleph_0})$, there is a ccc  
    Boolean algebra $\bfB$ with length $\mu$ (see Definition \ref{5.0A} 
    below) such that every infinite homomorphic image of $\bfB$ is of 
    cardinality $2^{\aleph_0}$. 
\end{enumerate}
\sn
If $\mu$ is a limit {cardinal} and $\cf(\mu) > \aleph_0$ we can demand the length is not obtained 
(see Definition \ref{5.0A}): if $\cf(\mu)=\aleph_0$ this is impossible.

Also, we can replace $\aleph_0$ here by any strong limit cardinal $\kappa$ of
cofinality $\aleph_0$ (see \ref{5.8}).

In \S 4 we deal with trees of the form $S\cup {}^{\omega>}2,$ where 
$S \subseteq {}^\omega 2$ is of cardinality $\lambda$.  

Note that \S 1, \S 2 are revised versions of parts of
\cite{Sh:136} and parallel to \cite{Sh:331}, and \S4 is a 
revised version of parts of \cite{Sh:262}.
The results in \S2 answer problems of Monk (presented in {Oberwolfach 1980}).

In \S3, we solve a problem of Boolean algebras of Monk on which the author
earlier gave a consistency result, using ideas from \S2. 

\S4 supersedes \cite[VIII 1.8]{Sh:a} and repeats \cite[1.2,1.3]{Sh:262}. 
Baldwin \cite{Bl89} has continued \cite[1.2-1.3]{Sh:262}. We
can apply this to models of $\varphi\in \bbL_{\aleph_1,\aleph_0}$,
probably using \cite{Sh:522}.

Recall that in \cite[Ch.VIII,1.8+1.7(2)]{Sh:a}, we proved that 
for pairs of first order
complete theories $(T,T_1)$ satisfying the hypothesis of Theorem \ref{3.1}
below
\[
\isoI(\lambda,T_1,T)\geq\min\{2^\lambda,\beth_2\}.
\]

\mn
We shall improve the result, replacing $\isoI(\lambda,T_1,T)$ 
by $\numbIE(\lambda,T_1,T)$. We improve the proof from 
\cite[VIII 1.8]{Sh:a}; in particular, we use
the trees $U_\eta$ defined in Fact \ref{3.3B}. They are subtrees of
${}^{\omega>}2$ as close to disjoint as we can manage.
 
We can use trees similar to $({}^{\omega\geq} 2,\lhd)$ with
finite or countable levels and heavier structure (i.e.~like pure conditions
in forcing notions as in \cite[\S2]{Sh:326}). As in 1.4(3),
we use here a weak form of representation: the amount of
similarity depends on the terms and formulas.

We can use such trees as in {\S}2 to build ``complicated,'' rigid-like
structures. In \cite[1.2,1.1(3)]{Sh:98} (more in \cite[1.4,
1.1]{Sh:105}) this was done for abelian groups: one step is getting $G
\supseteq \bbZ$ such that $G$ is $\aleph_1$-free of cardinality $\aleph_1$, 
$\bbZ$ not a direct summand of $G$). This was continued  in G\"obel 
and Shelah \cite{Sh:519} and G\"obel-Shelah-Str\"ongmann \cite{Sh:785}.


\begin{definition}\label{x2}
1) We say a structure $M$ is atomically $({<}\, \mu)$-stable \underline{when}: 
if\\ $A \subseteq M$ and $|A| < \mu$ then the set
$\{\tp_\qf(\bar a,A,M) : \bar a \in {}^{\omega >}\! M\}$  of possible types 
has cardinality $< \mu$.

\mn
2) We may write $\mu$-stable instead of `$({<}\,\mu^+)$-stable.'

\mn
3) For $\Theta$ a set of cardinals, when we write `$\Theta$-stable' we mean $\mu$-stable for every $\mu \in \Theta$.
\end{definition}

\bn
\centerline{*\qquad*\qquad*}

\bigskip
We thank Haim Horowitz, Mark Po\'or, and Miguel Cardona for many helpful comments.

\mn
\begin{notation}\label{z2}
1) $k,\ell,m,n$ will denote natural numbers.

\mn
2) Cardinals (infinite, if not stated otherwise) will be denoted by $\lambda,\mu,\kappa,\chi,\theta,\partial,\Upsilon$. $\partial$ will always denote a regular cardinal.

\mn
3) Ordinals will be denoted by $\alpha,\beta,\gamma,\delta,\eps,\xi,i,j$. $\delta$ will be a limit ordinal if not explicitly stated otherwise.
\end{notation}

\newpage

\section {Trees with structure}
\sectionmark{Trees}
\label{par1}

We deal here with ``relatives'' of $K^\omega_{\tr}$ which are more
complicated, strengthening our ability to carry out our constructions. 
We mean two things when we say `more complicated' --- namely:
\begin{enumerate}
    \item If $\eta \in \clT$ is not of the maximal level, there is a structure on the set of immediate successors. E.g.~pairs of members are used above it, so we have more induced non-embeddability {(involving the relevant $\psi$ and $\bfe$)}.
\sn
    \item $\clT$ has $\partial$-many levels, where $\partial$ may be $\aleph_0$ or be uncountable.
\end{enumerate}
This helps when using them for constructions, but we have to work more for the existence proofs (see \S1B).

In this section (and the next) we define and see what we can do for
$K^\partial_\ptr$, $\psi_\ptr$, $K^\partial_{\tr (n)}$, $\psi_{\tr(n)}$,
$K^\partial_{\tr(*)}$, $\psi_{\tr(*)}$ (which were introduced in
\cite{Sh:136}), getting the parallel of \cite[2.15\subref{7.11}]{Sh:331}; the main case is 
$\partial \defeq \aleph_0$. The reason for their 
introduction was for constructing certain Boolean algebras; 
we shall deal with these constructions later.

More specifically, \cite[2.2\subref{f5}]{Sh:E59} defines versions of ``$I$ {is} strongly $\varphi(\bar x,\bar y)$-unembeddable into $J$" and ``$K$ {has} [full and/or strong] $(\chi,\lambda,\mu,\kappa)$-bigness," so we can apply it to $(K,\varphi) = (K_\ptr^\partial,\psi_\ptr)$, or $(K_{\tr(h)}^\partial,\psi_{\tr(h)})$ or $(K_{\tr(h)}^\partial,\psi_{\tr(h)}')$, as defined in Definitions \ref{1.1},\ref{1.2} below. But below, essentially we choose more general $\varphi$-s represented by $\bfe$.

The relevant results are obtained by the existence of the super version, as in \cite{Sh:331} (see Definitions \ref{1.4},\ref{1.5}). 

Note that \S1 continues \cite[1.9\subref{7.5}]{Sh:331}.

\subsection{The frame}

\begin{definition}\label{1.1}
1) $K^\partial_\ptr$ is the class of $I$ such that:
\mn
\begin{enumerate}
    \item  The set of elements of $I$ is, for some linear order $J$, a subset of
\begin{align*}
    \set_{\ptr}[J] \defeq \big\{\eta : &\ \eta \text{ is a sequence of length $\leq\partial$, such that if}\\
    &\ i + 1 <\lh(\eta) \text{ then $\eta(i)$ has the form $\LL s,t\RR$ with }
    s <_J t,\\
    &\text{ and if } i + 1 = \lh(\eta) \text{ then } \eta(i)\in J \big\}.
\end{align*}
    Also, if $\eta\in I$, $i+1 < \lh(\eta)$, and $\eta(i)=\LL s,t\RR$ \underline{then} $(\eta\rest i) \caret \LL s\RR\in I$ and $(\eta\rest i) \caret \LL t\RR\in I$.
    Furthermore, the empty sequence belongs to $I$, and if $\delta < \lh(\eta)$ is a limit ordinal then $\eta \rest \delta \in I$.
\sn
    \item  The relations of $I$ are:
    \begin{enumerate}
        \item[$(\alpha)$]  $\eta\unlhd \nu$, meaning `$\eta$ is an initial segment of $\nu$' (i.e. $\eta = \nu \rest \lh(\eta)$); 
\sn
        \item[$(\beta)$]  $P_i \defeq \{\eta : \lh(\eta) = i\}$
\sn
        \item[$(\gamma)$]  ${<_1} \defeq$ 
        $$
        \big\{(\eta,\nu) : \text{for some } i < \partial,\ \lh(\eta) = \lh(\nu) = i+1,\ \eta(i) <_J \nu(i),\ \eta \rest i = \nu \rest i\big\}
        $$
        \item[$(\delta)$]  $\Eq_i \defeq \{\LL\eta,\nu\RR : \lh(\eta),\lh(\nu) \geq i\, \wedge\, \eta \rest i = \nu \rest i\}$ 
\sn
        \item[$(\eps)$]  $\Suc_L \defeq$ 
        \begin{align*}
            \big\{ \LL\eta,\nu\RR : &\text{ for some $i<\partial$ and } s<_J t \text{ we have } \eta \rest i = \nu \rest i, \\ 
            &\ i+1 = \lh(\eta) < \lh(\nu),\ \nu(i)=\LL s,t\RR \text{ and } \eta(i) = s \big\}
        \end{align*}
        \item[$(\zeta)$]  $\Suc_R \defeq$ 
        \begin{align*}
            \big\{ \LL\eta,\nu\RR : &\text{ for some $i<\partial$ and } s<_J t \text{ we have } \eta \rest i = \nu  \rest i,\\ 
            &\ i+1 = \lh(\eta) < \lh(\nu),\ \nu(i) = \LL s,t\RR \text{ and } \eta(i) = t \big\}
        \end{align*}
        \item[$(\eta)$]  An individual constant  $\LL\ \RR$.
\sn
        \item[$(\theta)$] Functions $\Res^L_\alpha,\Res^R_\alpha$ such that $\Res^L_\alpha(\eta) = (\eta \rest \alpha) \caret \LL s\RR$ and $\Res^R_\alpha(\eta) = (\eta\rest \alpha) \caret \LL t\RR$ when 
        $\eta(\alpha) = \LL s,t\RR$ and $\alpha+1<\lh(\eta)$, and 
        $\Res^L_\alpha (\eta)=\Res^R_\alpha (\eta)=\eta$ otherwise.
    \end{enumerate}
\end{enumerate}
If we omit the superscript $\partial$, we mean $K^{\aleph_0}_\ptr$.

\mn
2) Let
\begin{align*}
    \psi_\ptr(x_0,x_1;y_0,y_1) = \bigvee\limits_{i+1<\partial} & \big[P_{i+1}(x_1) \wedge P_{i+1}(y_1) \wedge P_\kappa(x_0) \wedge x_0 = y_0 \\ 
    &\ \wedge \Suc_L(x_1,x_0) \wedge \Suc_R(y_1,y_0) \wedge x_1 <_1 y_1 \big].
\end{align*}
This depends on $\partial$, but we usually suppress this parameter.

\sn
3) $I \in K^\partial_\ptr$ is standard \underline{iff} in (1)(A), $J$ is a set
of ordinals with the natural order, or at least a well-ordering (usually we
shall use those).
\end{definition}

\begin{definition}\label{1.2}
1) For $h : \partial \to \omega \setminus\{0\}$, the class
$K^\partial_{{\tr(h)}}$ is defined like $K^\partial_{\ptr}$, but replacing pairs
by increasing $h(i)$-tuples at level $i$. That is, $I \in K_{\tr(h)}^\partial$ \underline{iff}
\mn
\begin{enumerate}
    \item  the set of elements of $I$ is, for some linear order $J$, a subset of
\begin{align*}
    \set_{\tr(h)}(J) \defeq \big\{\eta : &\ \eta\text{ is a sequence of length} \leq \partial,\\
    &\text{ for }i+1<\lh(\eta),\ \eta(i)\text{ has the form }\LL s_0,\ldots, s_{h(i)-1}\RR\\
    &\text{ such that }s_0<_Js_1<_J\ldots<_J s_{h(i)-1}, \text{ and}\\
    &\text{ for }i+1=\lh(\eta),\ \eta(i)\in J\big\}.
\end{align*}
    Also, if $\eta\in I$, $i+1 < \lh(\eta)$, $m < h(i)$ and 
    $\eta(i) = \LL s_0,\ldots,s_{h(i)-1}\RR$ then 
    $(\eta\rest i) \caret \LL s_m\RR\in I$. Furthermore, the empty sequence belongs to $I$, and if $\delta < \lh(\eta)$ is a limit ordinal then $\eta\rest\delta\in I$.
\sn
    \item  The relations of $I$ are:
\sn
    \begin{enumerate}
        \item[$(\alpha)$]  $\eta \unlhd \nu$, which holds iff 
$\eta=\nu\rest\lh(\eta)$.
\sn
        \item[$(\beta)$]  $P_i \defeq \{\eta : \lh(\eta) = i\}$. (We let $P_{<\partial} \defeq \bigcup\limits_{i<\partial} P_i$.)
\sn
        \item[$(\gamma)$]  ${<_1} \defeq \{\LL\eta,\nu\RR : \lh(\eta) = \lh(\nu) = i+1,\ \eta(i) <_J \nu(i),\ \eta \rest i = \nu \rest i\}$
\sn
        \item[$(\delta)$]  $\Eq_i \defeq \{\LL\eta,\nu\RR : \lh(\eta),\lh(\nu) \geq i,\ \eta\rest i = \nu \rest i\}$
\sn
        \item[$(\eps)$]  For $m < h(i)$ and $ i < \partial$:
\begin{align*}
        \Suc_{i,m} = \{\LL\eta,\nu\RR : &\ \eta\rest i = \nu\rest i,\ \lh(\eta) = i+1 < \lh(\nu),\\
        &\ \nu(i)=\LL s_0,\ldots ,s_{h(i)-1}\RR, \eta(i)=s_m\}
\end{align*}

        \item[$(\zeta)$]  An individual constant $\LL\ \RR$.
\sn
        \item[$(\theta)$]   Functions $\Res^m_i$ such that
        $\Res^m_i(\eta) = (\eta\rest i) \caret \LL s_m\RR$ when 
        \[
            \eta(i) = \LL s_0,\ldots,s_{h(i)-1}\RR,\ i + 1 < \lh(\eta), \text{ and } m < h(i).
        \]
        Above, 
        \begin{enumerate}
            \item If $n \ge h(\lh(\eta(i)))$, $\eta(i) = \LL s_0,\ldots,s_{h(i)-1}\RR$, and $i + 1 < \lh(\eta)$ then we stipulate $\Res_i^m(\eta) = \eta$.
\sn
            \item If $\lh(\eta) \leq i$ then we stipulate $\Res_i^m(\eta) = \eta$.
        \end{enumerate}   
    \end{enumerate}
\end{enumerate}

\sn
1A) For $\Lambda \subseteq I \in K_{\tr(h)}^\partial$, let 
\begin{align*}
    \closure_I(\Lambda) \defeq \Big\{ \eta \in I : &\ \lh(\eta) \text{ is a limit ordinal and }  \\
    &\ \big[i < \lh(\eta) \, \wedge\, \ell < h(i) \Rightarrow \Res_i^k(\eta) \in \Lambda\} \big] \Big\}
\end{align*}

\sn
1B) For $I \in K_{\tr(h)}^\partial$, let
\begin{itemize}
    \item $\sseq_{i,m}(I) \defeq \big\{\eta \caret \LL\varrho\RR : \eta \in P_{i+1}^I,\ \varrho \in \inc_m(J) \big\}$
\sn
    \item $\sseq_i(I) \defeq \bigcup\limits_{m < h(i)} \sseq_{i,m}(I)$
\sn
    \item $\sseq(I) \defeq \bigcup\limits_{i < \partial} \sseq_i(I)$
\end{itemize}

\sn
2)  If $\bigwedge\limits_{i<\partial} [h(i)=n]$ we may write
$K^\partial_{\tr(n)}$, so for $n=2$ we get $K^\partial_\ptr$ up to some
renaming. Similarly for $\psi_\ptr$.

\sn
3) If $\bigwedge\limits_{i<\kappa} [h(i) = i \mod \omega]$ we may write
$K^\partial_{\tr(*)}$. We say ``$I \in K^\partial_{\tr(h)}$ is standard'' in
case the underlying set $J$ is well-ordered (usually a set of
ordinals). 

\sn
4) When we write $\eta(\alpha)(\ell)$, we mean $\eta(\alpha)$ if $\lh(\eta) = \alpha+1$ and $\Res^\ell_\alpha(\eta)$ if $\alpha+1 < \lh(\eta)$.
\end{definition}

\mn
\begin{definition}\label{1.2A}
Continuing \ref{1.2}, we define some formulas relevant to the class $K_{\tr(h)}^\partial$.

\sn
1) $\psi_{\tr(h)}(\bar x;\bar y)$, where $\bar x=(x_0,x_1)$, $\bar y = (y_0,y_1)$, 
    is\footnote{Below, the intention is  $y_0\rest i = x_\ell \rest i$ and $y_0(i) = \LL x_0(i),\ldots,x_{h(i)-1}(i)\RR$.}
\begin{align*}
(x_0=y_0) \wedge \ P_\partial(y_0) \wedge \bigvee\limits_{i<\partial} 
&\big[P_{i+1}(x_1) \wedge P_{i+1}(y_1) \wedge (x_1 <_1 y_1) \ \wedge \\
&\ \Suc_{i,0} (x_1,x_0) \wedge \Suc_{i,h(i)-1} (y_1,y_0) \big]\\
\end{align*}

\noindent
2) We define $\psi'_{\tr(h),i}(\bar x;\bar y)$ as follows, where $\lh(\bar x) = \lh(\bar y) = h(i) + 1$:
$$
x_0=y_0 \wedge P_\partial (y_0) \wedge 
\bigwedge\limits_{\ell=0}^{h(i)-2} \! [x_{\ell+1} = y_\ell]\ \wedge\\
\bigwedge\limits_{\ell=1}^{h(i)} \!\big[P_{i+1} (x_\ell) \wedge P_{i+1}
(y_\ell)  \big] \bigwedge\limits_{\ell<h(i)}\Suc_{i,\ell}(x_{\ell+1},x_0).
$$

\sn
3) Assume  $\alpha = \sup(\rang(h))$ and let 
$\bar x = \LL x_\ell : \ell < 1 + \alpha\RR$ and 
$\bar y = \LL y_\ell : \ell < 1 + \alpha\RR$ 
(noting $\alpha \leq \omega$). Then we define 
$\psi'_{\tr(h)}(\bar x;\bar y)$ as 
$$
\bigvee\limits_{i < \partial} \psi'_{\tr(h),i} \big( \bar x \rest (h(i) + 1);\bar y \rest (h(i) + 1) \big).
$$

\sn
4) For $\bar\bfe$ as in Definition \ref{1.4} below, we define $\psi_{\tr(h),\bar\bfe}(\bar x; \bar y)$ as 
$\bigvee\limits_{i < \partial} \psi_{\tr(h),\bar\bfe,i}$,
where 
\begin{itemize}
    \item $\lh(\bar x) = \lh(\bar y) = \sup (\rang(h))$
\sn
    \item $\psi_{\tr(h),\bar\bfe,i} \defeq \bigwedge \big\{ \psi_{\tr(h),i,u_1,u_2}\big( \bar x \rest (h(i) + 1);\bar y \rest (h(i) + 1) \big) : (u_1,u_2) \in \bfe_i \big\}$, where
\sn
    \item $\psi_{\tr(h),i,u_1,u_2} \defeq$
\begin{align*}
    \Big(&x_0=y_0 \wedge P_\partial (y_0) \wedge \big\{ x_\ell = y_k : \ell \in u_1,\ k \in u_2,\ |\ell \cap u_1| = |k \cap u_2| \big\}\\
    &\wedge
    \bigwedge\limits_{\ell=1}^{h(i)} \!\big[P_{i+1} (x_\ell) \wedge P_{i+1}
    (y_\ell) \wedge \Res_i^\ell(x_0) = x_\ell \big]\Big).
\end{align*}
\end{itemize}
\end{definition}

\mn
\begin{remark}\label{1.3}
1) Here, when dealing with $K^\omega_\ptr$ (or $K^\partial_{\tr(n)}$,
$K^\partial_{\tr(*)}$, $K^\partial_{\tr(h)}$; those are parallel cases),
we introduce the ``super$^*$" version, parallel to Definitions
\cite[1.1\subref{7.1}, 1.4\subref{7.2}]{Sh:331}.\footnote{
    And see more versions in \cite[1.5\subref{7.3}, 1.6\subref{7.3A}]{Sh:331}.
} 
So the easy {case}
\cite[1.9\subref{7.5}]{Sh:331} has to be redone, and
\cite[Claim 1.8(2)\subref{7.5}]{Sh:331} is no longer of any help and we should prove a
parallel. The role of $\bar\bfe$ here corresponds in the role of
$\psi_\tr$ in
\cite[\S2]{Sh:E59}, \cite[\S1]{Sh:331}.

\sn
2) Concerning Definition \ref{1.4}, the reader is advised to
concentrate on clauses (1) and (3)(A),(D).
\end{remark}

\begin{definition}\label{1.4}
Let $h : \partial \to \partial \setminus \{0\}$, and $\bar\bfe$ be a
function with domain $\partial$, with each $\bfe(i) = \bfe_i$ an equivalence relation on some subset of 
$\clP(h(i))$ satisfying
\[
u_1\ \bfe_i\ u_2\ \Rightarrow\ \vert u_1\vert = \vert u_2\vert.
\]

\mn
For this definition we identify a set (of natural numbers or ordinals)
with an increasing sequence enumerating it. 
When defining $\bfe_i$ we may ignore classes which are
singleton; see clause (5) on default values.

\sn
1) For $I \in K^\partial_{\tr(h)}$, $J \in K^\partial_{\tr(h')}$ 
and cardinals $\mu_1,\mu_2,\pinkappa$, we say $I$ is 
$(\mu_1,\mu_2,\pinkappa)$-\emph{super}-$\bar\bfe$-\emph{unembeddable\footnote{
    But if $\mu_1 = \mu_2 = \mu$, we write $(\mu,\kappa)$ instead of $(\bar\mu,\kappa)$.
} into} $J$ 
(for $K^\partial_\tr(h)$) \underline{when}:
\mn
\begin{enumerate}
    \item[$(*)^{I,J}_{\bar\mu,\pinkappa,\bar\bfe}$]  For every large enough regular cardinal $\chi$, for every $x \in \clH(\chi)$, for a fixed well-ordering $<^*_\chi$ of the set $\clH(\chi)$ and 
    $f_0 : {}^{\pinkappa>}\!I \to \mu_2$, $f_1 : \sseq(I) \to {}^{\kappa>}\!J$, there are 
    $\LL M_i,N_i : i < \partial\RR$ such that:
    \begin{enumerate}
        \item[$(i)$]  $M_i \prec N_i \prec M_j \prec (\clH(\chi),\in,<^*_\chi)$ for $i < j < \partial$.
\sn
        \item[$(ii)$]  $M_i \cap \mu_1 = N_i \cap \mu_1$ for $i < \partial$.
\sn
        \item[$(iii)$] $\pinkappa \subseteq M_0$, $\partial \subseteq M_0$.
\sn
        \item[$(iv)$]  $I,J,\mu,\kappa,\partial,f_0,f_1,h,x$ belong to $M_0$.
\sn
        \item[$(v)$]  There is $\eta \in P^I_\partial$ such that
        $\eta \rest i \in M_i$ for every $i < \partial$. Also, for all $i$ large enough, 
        ``$\Res^0_i(\eta)$, $\Res^1_i(\eta),\ldots, \Res^{h(i)-1}_i (\eta) \in N_i \setminus M_i$ realize the same Dedekind cut on
        $$
        \big\{ (\eta\rest i) \caret \LL s\RR \in I : s \in M_i \cap \{\nu(i) : \nu \in P_{i+1}^I,\ \nu \rest i = \eta \rest i\} \big\}."
        $$
        (Recall that $<^I_1$ linearly orders $\{(\eta \rest i) \caret \LL s\RR : s \in J\}$.) 
\sn
        \item[$(vi)$] In clause $(v)$, we may add: if $h(i)>1$ and $u_1\ \bfe_i\ u_2$ then 
        \begin{enumerate}
            \item[$(\alpha)$] If $\ell_1 \in u_1 \wedge \ell_2 \in u_2 \wedge 
            |u_1\cap \ell_1| = |u_2\cap\ell_2|$ then $f_0(\Res_i^{\ell_1}(\eta)) = f_0(\Res_i^{\ell_2}(\eta))$, 
            and\footnote{
                Manipulating $f_0$, the ``also" follows.
            } 
            also $f_1(\Res^{\ell_1}_i(\eta))$ and $f_1(\Res^{\ell_2}_i(\eta))$ have the same length.

            \item[$(\beta)$] For $u \subseteq h(i)$, let 
            $\bar\nu_{\eta,i,u}$ be the concatenation of the sequences $f_1(\Res^\ell_i(\eta))$ for $\ell\in u$.
            
            \underline{Then} the sequences $\bar\nu_{\eta,i,u_1},\bar\nu_{\eta,i,u_2} \in {}^{\kappa>}\!J$ 
            \begin{enumerate}[$\bullet_1$]
                \item Satisfy $f_0(\bar\nu_{\eta,i,u_1}) = f_0(\bar\nu_{\eta,i,u_2})$.
\sn
                \item They realize the same atomic type\footnote{
                    If we omit `closure' then later we will need to say more: for every $\nu \in P^J_\partial$,
\[
                    \Big(\bigcup_{i<\partial} M_i \Big) \cap \{\Res^\ell_i(\nu) : \ell < h(i),\ i < \partial\}
\]
                    is included in some $M_m$. Also, in \ref{1.6C} we need  
            $ {\partial_1} \defeq \sum\limits_{\alpha<\partial} \big( |\alpha|^{\aleph_0} \big)^+$.
                } 
                over the closure of $J\cap M_i$ in $J$. 
            \end{enumerate}      
        \end{enumerate}
\sn
        \item[$(vii)$] For every $\bar\nu \in {}^{  
        {\omega}>}\!J$, 
        for every $i < \partial$ large enough, we can replace $J \cap M_i$ by $J \cup M_i \cup \{\nu_\ell : \ell \in \lh(\bar\nu)\}$ in clause $(vi)(\beta)$ above.
    \end{enumerate}
\end{enumerate}

\sn
2) For $I,J \in K^\partial_{\tr(h)}$ and cardinals $\mu_1,\mu_2,\kappa$,
we say\footnote{
    This is helpful in constructing Boolean algebras as in \S2 in more cardinals without using Definition \ref{1.4}(1) ({or} even $\psi'_{\rt(h)}$), but this is the minor variant and the reader can ignore it.
}
that $I$ is $(\mu_1,\mu_2,\partial)$-super-$\bar\bfe$-unembeddable$'$ 
into $J$ (for $K^\partial_{\tr(h)}$) \underline{when}:
\mn
\begin{enumerate}
    \item[$(*)'_{I,J,\mu,\partial,\bar\bfe}$] 
    For every large enough $\chi$ and 
    $x \in \clH(\chi)$, for a fixed well-ordering $<^*_\chi$ of $\clH(\chi)$, 
    there exists $M$ such that: 
    \begin{enumerate}
        \item[$(i)$]  $M \prec(\clH(\chi), \in,<^*_\chi)$
\sn
        \item[$(ii)$]  $x\in M$
\sn
        \item[$(iii)$]  $\|M\| = \partial$
\sn
        \item[$(iv)$]  There is $\eta \in P^I_\partial$ such that 
        $$
        i < \partial \wedge \ell < h(i) \Rightarrow \Res_i^\ell(\eta) \in M,
        $$  
        and for every function $f\in M$ from $I$ to $\mu$, for unboundedly many $i < \partial$, we have: 
        \begin{enumerate}
            \item[$\circledast$] If $\ell'_0 < \ldots < \ell'_{k-1} < h(i)$,  $\ell^{''}_0 < \ldots < \ell^{''}_{k-1} < h(i)$, and
            $\{\ell'_0, \ldots, \ell'_{k-1}\}\ \bfe_i\ \{\ell^{''}_0,\ldots,\ell^{''}_{k-1}\}$, \underline{then} 
            $$
            f\big(\LL\Res^{\ell'_m}_i(\eta) : m < k\RR \big) = f\big(\LL \Res^{\ell^{''}_m}_i(\eta) : m < k\RR \big).
            $$
        \end{enumerate}
\sn
        \item[$(v)$]  If $\nu\in P^J_\partial$ then either $\nu \in M$ \underline{or} for some $k < \partial$ we have 
        $\nu \rest k \in M$ and $\nu \rest (k+1) \notin M$. 
    \end{enumerate}
\end{enumerate}
\mn
3) Let $\max(\bfe_i) \defeq \max\{|u| : u\in \dom(\bfe_i)\}$.
\begin{enumerate}
    \item Let $\bar\bfe^0$ be defined by 
$$
    \bfe^0_i = \bar\bfe^0(i) \defeq \big\{\big(\{\ell\},\{k\}\big) : \ell,k < h(i) \big\}.
$$ 
    \item Let $\bar\bfe^1$ be defined so that 
    $\bfe^1_i = \bar\bfe^1(i)$ 
    is the closure of 
    $$
    \big(\{0,\dots,h(i)-2\}, \{1,\dots, h(i)-1\}\big)
    $$ 
    to an equivalence relation.
\sn
    \item Similarly, $\bfe^2_i = \bar\bfe^2(i)$ 
    is the closure of  
$$
    \big(\{0,\ldots\lfloor h(i)/2\rfloor-1\},\{\lfloor h(i)/2\rfloor,\ldots, 2\lfloor h(i)/2\rfloor-1\} \big).
$$ 
    \item If $\bar\bfe = \bar\bfe^0$ we may omit {the superscript}.
\end{enumerate}

\sn
4)  For $\LL I_\xi : \xi \in W\RR$, $W$ a set of ordinals,
$I_\xi \in K^\partial_{\tr(h)}$ (standard, for simplicity) letting 
$$
\zeta_* \defeq \sup\!\big(W \cup \big\{\eta(i) : \eta \in \textstyle\bigcup\limits_{\xi \in W} I_\xi \text{ and } \lh(\eta) = i+1 < \partial \big\} \big)+1,
$$ 
we define $\sum\limits_{\xi\in W} I_\xi \in K^\partial_{\tr(h)}$ as 
$\big\{\LL\ \RR\big\} \cup \big\{\LL\xi\RR 
\mathop{\otimes}\limits_{\zeta_*} \eta : \xi \in W \text{ and }\eta\in
I_\xi\big\}$, where $\mathop{\otimes}\limits_{\zeta_*}$ is as in part (5). 

\mn
5) Recall $\xi \mathop{\otimes}\limits_{\zeta_*}\eta$ is 
$\LL\ \RR$ if $\eta=\LL\ \RR$,
and is $\LL\zeta_* \times \xi+\eta(0),\eta(1),\eta(2),\ldots\RR$ otherwise.

\mn
6) Above, if   
   $\mu_1 = \mu =\mu_2$  
we may write  
{$\mu$ instead $\mu,\mu$}
\end{definition}

\mn
\begin{definition}\label{1.5}
1) $K^\partial_{\tr(h)}$ has the $(\chi,\lambda,\mu,\pinkappa)$-\emph{super}-$\bar\bfe$-\emph{bigness property} \underline{when} there are standard $I_\zeta\in K^\partial_{\tr(h)}$
for $\zeta < \chi$ with $|I_\zeta| = \lambda$ such that $I_\zeta$ is 
$(\mu,\pinkappa)$-super-$\bar\bfe$-unembeddable into $I_\eps$ for each 
$\zeta \neq \eps < \chi$.

\mn
2) $K^\partial_{\tr(h)}$ has the \emph{full} $(\chi,\lambda,\mu,\pinkappa)$-super-$\bar\bfe$-bigness property \underline{when} there are standard 
$I_\zeta\in K^\partial_{\tr(h)}$ for $\zeta<\chi$ with 
$|I_\zeta| = \lambda$  such that $I_\zeta$ is 
$(\mu,\pinkappa)$-super-$\bar\bfe$-unembeddable into 
$J_\zeta = \sum\limits_{\substack{\eps<\chi\\ \eps\neq\zeta}}I_\eps$ for each $\zeta < \chi$.

\mn
3) We may also add superscripts to distinguish slightly different super-bigness properties: 
super$^\nr$ will be used for the properties defined in parts (1) and (2) above;
super$^\vr$ will be almost identical, but
we replace unembeddable by unembeddable$'$ (i.e.,
in Definition \ref{1.4} we replace $(*)^{I_\zeta,J_\zeta}_{\mu,\pinkappa,\bar\bfe}$ 
by $(*)'_{I_\zeta,J_\zeta,\mu,\pinkappa,\bar\bfe}$). 

\mn
3A) We may replace $\lambda$ by $\bar{\lambda} = (\lambda_0,\lambda_1)$ 
if `$|I_\xi| = \lambda$' is replaced by 
$$
`|I_\xi| = \lambda_0 \text{ and } \big|\{\eta \in I_\xi : \lh(\eta) < \partial\} \big| = \lambda_1.\text{'}
$$

\mn
4) {Whenever we state a theorem, definition, or claim which does not depend on the specific version of bigness, we will write `super$^\x$.'}
\end{definition}

\mn
\begin{remark} 
Also, $K^\omega_\tr$ can be brought into the framework above as a 
specific case (i.e. $h$ is constantly 1).
\end{remark}

\mn
\begin{claim}[\textbf{Monotonicity}]\label{1.6A}
For every $h : \partial \to \omega\setminus \{0\}$, we have:

\sn
$1)$ If $K^\partial_{\tr(h)}$ has the [full] 
$(\chi_1,\lambda_1,\mu,\kappa)$-super$^\x$ $\bar\bfe$-bigness 
properties, $\chi_2 \leq \chi_1$, and
$\lambda_2 \geq \lambda_1$, \underline{then}  $K^\partial_{\tr(h)}$ has the
[full] $(\chi_2,\lambda_2,\mu,\kappa,\theta)$-super$^\x$ $\bar\bfe$-bigness
property.

\sn
$2)$  If $K^\partial_{\tr(h)}$ has the full 
$(\chi,\lambda,\mu,\kappa)$-super$^\x$ $\bar\bfe$-bigness 
property and $\chi_1 \defeq \min\{\chi,\lambda\}$, \underline{then} 
$K^\partial_{\tr(h)}$ has the 
$(2^{\chi_1},\lambda,\mu,\kappa)$-super$^\x$ $\bar\bfe$-bigness. 
\end{claim}

\begin{PROOF}{\ref{1.6A}}
1) Straightforward.

\sn
2) Similar to 
\cite[1.9(1)\subref{7.5}]{Sh:331}, but we elaborate. 

If $\LL I_\alpha : \alpha < \chi\RR$ exemplifies ``$K^\partial_{\tr(h)}$
has the full $(\chi,\lambda,\mu,\pinkappa,\theta)$-super$^\x$ $\bar\bfe$-bigness 
property,'' $\chi_1 \defeq \min\{\chi,\lambda\}$ and $h(0) = n_*$,
then we let $J_A = \sum\limits_{\alpha\in A} I_\alpha$ for $A \subseteq \chi_1$
(see Definition \ref{1.4}(4)).

Let $\LL A_\alpha : \alpha < 2^{\chi_1}\RR$ be such that 
$A_\alpha \subseteq \lambda$ and 
$\alpha\neq\beta\ \Rightarrow\ A_\alpha \not\subseteq A_\beta$. Now\\ 
$\LL J_{A_\alpha} : \alpha < 2^{\chi_1}\RR$
exemplifies ``$K^\omega_{\tr(h)}$ has the $(2^{\chi_1},\lambda,\mu,\kappa,
\theta)$-super$^\x$ $\bar\bfe$-bigness property."
\end{PROOF}

\medskip
On the [full] strong $(\chi,\lambda,\mu,\pinkappa)$-bigness property 
(and strongly finitary version) see \cite[2.5\subref{2.3}]{Sh:E59}. 
By \ref{1.6B} below, for $\psi_{tr(h)}$ 
from Definition \ref{1.2}(2) it is a consequence of the super$^\x$ 
version, and as in \cite{Sh:E59}, \cite{Sh:331} it is useful.

\begin{claim}\label{1.6B}
If $K^\partial_{\tr(h)}$ has the [full] $(\chi,\lambda,
\mu^{<\pinkappa}, \pinkappa)$-super-bigness property and 
$\chi \leq \lambda$, \underline{then} $K^\partial_{\tr(h)}$ has the
[full] strong  $(\chi,\lambda,\mu,\pinkappa)$-bigness property for 
$\psi_{\tr(h)}$ for functions $f$
which are strongly finitary on $P_\partial$.
\end{claim}

\begin{PROOF}{\ref{1.6B}}
The result follows by the definitions and \ref{1.6C} below.
\end{PROOF}

\mn
Analogously to \cite[1.10\subref{7.5A}]{Sh:331}, we have:

\begin{claim}\label{1.6C}
Recalling \emph{\ref{1.4}(1), (3)(A),(D)}, if $(*)^{I,J}_{\mu_1,\pinkappa_1}$ (where $\mu_1 = \mu^{<\pinkappa}$, 
$\pinkappa_1 = \pinkappa$, $\{I, J\}\subseteq K^\partial_{\tr(h)}$ 
are standard\footnote{
    See Definition \ref{1.4}(1).
}) 
and $h \in {}^\partial(\omega \setminus \{0\})$,
\underline{then}  $I$ is strongly $(\mu,\pinkappa,\psi_{\tr(h)})$-unembeddable
into $J$ for embeddings which are strongly finitary on $P^I_\partial$.
\end{claim}

\begin{PROOF}{\ref{1.6C}}
Recalling \ref{1.4}(3), we have $\bar\bfe = \bar\bfe^0$.
Without loss of generality, $I$ and $J$ are subsets of 
${}^{\partial\geq}({}^{\partial>}\theta)$
for some cardinal $\theta$, and let
$<^*$ be a well-ordering of $\cM_{\mu,\pinkappa}(J)$ (which respects being
a subterm). Suppose $f$ is a function from $I$ into 
$\cM_{\mu,\pinkappa}(J)$, so for $\eta \in I$, let
\[
f(\eta) \defeq \sigma_\eta(\nu_{\eta,0},\ldots,\nu_{\eta,i},\ldots)_{i < \alpha_\eta}
\]

\mn
for some term $\sigma_\eta$, ordinal $\alpha_\eta < \kappa$, and 
$\nu_{\eta,i}\in J$.  $f$ is \emph{strongly finitary} on $P_\partial$, which means
\[
\eta \in P^I_\partial\, \Rightarrow\,  \alpha_\eta < \omega \wedge [
\sigma_\eta \text{ has finitely many subterms}].
\]

\mn
Let $\chi$ be regular large enough and $<^*$ a well-ordering of $\clH(\chi)$, and define
\[
g(\eta) = \big\{\alpha : \text{the $\alpha^\tthh$ element by $<^*$ is a subterm of } f(\eta) \big\}
\]

\sn
for $\eta\in P^I_\partial$
(so we use the ``strongly finitary'' only to ensure that $g(\eta)$ is finite for 
$\eta \in P^I_\partial$).

We define functions $f_0,f_1$ as follows.
\begin{enumerate}
    \item [$\boxplus_1$]
    \begin{enumerate}
        \item $\dom(f_1) = \dom(f_0) \defeq \bigcup\limits_{i<\partial} P_i^I$

        \item If $\eta \in P_i^I$ and 
        $$
        f(\eta) = \sigma_\eta(\ldots,\nu_{\eta,\ell},\ldots)_{\ell<\alpha_\eta} = \LL\sigma_{\eta,\iota}(\ldots,\nu_{\eta,j},\ldots) : \iota < \beta_\eta\RR
        $$ 
        then
        \begin{enumerate}
            \item $f_0(\eta) \defeq \cd(\LL \sigma_{\eta,\iota} : \iota < \beta_\eta\RR)$, where
\sn
            \item $\cd$ is the $<_\chi^*$-first one-to-one function from some subset of\\ $\{\sigma_{\nu,\iota} : \nu \in \bigcup\limits_{j<\partial} P_j^I,\ \iota < \beta_\nu\}$ into $\mu$.
\sn
            \item $f_1(\eta) \defeq \LL \nu_{\eta,\iota} : \iota < \alpha_\eta\RR$.
        \end{enumerate}
    \end{enumerate}
\end{enumerate}
Choose $\chi$ large enough (e.g.~$\chi$ strong limit such that $I,J \in \clH(\chi)$) and $<_\chi^*$ a well-ordering of $\clH(\chi)$ which respects 
$$
`a \text{ is a subterm of $b$'}\ \Rightarrow\ a <_\chi^* b.
$$
We now apply \ref{1.4}(1) with our present $I,J,f_0,f_1$, and get
$\LL M_i, N_i : i < \partial\RR$ as in the conclusion of 
Definition \ref{1.4}(1). Let $\eta\in P^I_\partial$ be as 
in clauses \ref{1.4}(1)$(v)$--$(vii)$.
Apply clause \ref{1.4}(1)$(vii)$ to the sequence $\LL\nu_{\eta,j} : j  \alpha_\eta\RR$ (so $\alpha_\eta$ is finite).

Let $i_\bullet < \partial$ be large enough (as in \ref{1.4}(1)$(vii)$.)
Let $x_0 = y_0 = \eta$ and $x_1 = (\eta \rest i_\bullet) \caret \LL \alpha_1\RR$,  
$y_1 = (\eta \rest i_\bullet) \caret \LL \alpha_{h(i_\bullet)-1}\RR$, where 
$\eta(i_\bullet) = \LL\alpha_\ell : \ell < h( i 
)\RR$. 

Now, using clause $(vi)(\beta)$+$(vii)$ of $(*)_{\mu,\pinkappa}^{I,J}$, the rest should be clear.
\end{PROOF}

\mn
\begin{lemma}\label{1.7}
Let $h \in {}^\omega(\omega\setminus \{0,1\})$. 

\mn
$1)$ $K^\partial_\ptr$ and $K^\partial_{\tr(h)}$ have the full 
$(\lambda,\lambda,\mu,\pinkappa)$-super--bigness property \underline{when} $\lambda > \mu^\partial + \lambda^{<\kappa}$.

\mn
$2)$   Above, we can deduce that $K^\partial_{\tr(h)}$ has the full 
$(\lambda,\lambda,\mu,\kappa)$-$\psi_{\tr(h)}$-bigness property. 
\end{lemma}

\begin{PROOF}{\ref{1.7}}
See \ref{1.15} and the end of \S1; we shall prove more in \S1B. (In fact, we can prove more as in \cite{Sh:331}.)
\end{PROOF}

\mn
\begin{claim}\label{1.8}
Let $h$ be as in \emph{\ref{1.7}}, and assume $\mu \ge \partial$.

\sn
$1)$  Let $I\in K^\partial_{\tr(h)}$. \underline{Then} $I$ is atomically $\mu$-stable
\underline{iff} $(A)$ and $(B)$ hold, where
\mn
\begin{enumerate}[$(A)$]
    \item  For $i < \partial$ and $\eta\in P^I_{i+1}$, the linear order 
    $$
    \big(\{\nu \in P^I_{i+1} : \nu \rest i = \eta \rest i\},{<^I_1}\big)
    $$ 
    is atomically $\mu$-stable (i.e.~for every subset of cardinality 
    $\leq\mu$, only $\leq\mu$-many Dedekind cuts are realized).
\sn
    \item  For any $I' \subseteq I$ with $|I'| \leq \mu$, the set
    $$
    \big\{\eta \in P^I_\partial : i < \partial \wedge \ell < h(i)\, \Rightarrow\, \Res^\ell_i(\eta) \in I' \big\}
    $$ 
    has cardinality $\leq\mu$. (So if $\mu = \mu^\partial$ then this certainly holds.)
\end{enumerate}
\mn
$2)$  For $\mu = \cf(\mu) > \partial$, ``atomically $({<}\,\mu)$-stable'' is
characterized similarly (for $\mu=\chi^+$, this means ``atomically
$\chi$-stable'').

\sn
$3)$  If $I\in K^\partial_{\tr(h)}$ is standard, $\mu = \cf(\mu)$, and
$$
\alpha < \mu \Rightarrow |\alpha|^\partial < \mu,
$$ 
\underline{then}  
$I$ is atomically $({<}\,\mu)$-stable.

\sn
$4)$  The family of ``atomically $({<}\,\mu)$-stable $I \in K^\partial_{\tr(h)}$\!" 
is closed under well-ordered sums.
\end{claim}

\begin{PROOF}{\ref{1.8}}
1)   Let $J\subseteq I$ be of cardinality $\leq\mu$. 
Without loss of generality 
\mn
\begin{enumerate}
    \item[$\boxtimes_1$]  $\eta\in J \wedge \lh(\eta) > i+1 \wedge \ell < h(i)\, \Rightarrow\, \Res^\ell_i(\eta) = (\eta \rest i) \caret \big\LL \eta(i)(\ell)\big\RR \in J$,
\end{enumerate}
[Why? Recall $\mu \geq \partial$.]

\mn
Let 
\begin{align*}
    J' \defeq \big\{ &\eta \rest (i+1) : \eta \in I,\ i < \lh(\eta) \big\}\\ 
    \cup\ \big\{& \eta \in J : \lh(\eta) \text{ is limit, and } j < i \wedge m < h(j)  \Rightarrow \Res_j^m(\eta) \in I \big\},
\end{align*} 
and for $\nu \in J'$ let\footnote{
    As clearly $\nu \lhd \eta \Rightarrow \lh(\eta) \ge \lh(\nu) + 1$.
} 
$$
J^*_\nu \defeq \big\{\eta \in I\setminus J : \nu \lhd \eta \text{ and } \eta \rest (\lh(\nu) + 1) \notin J \big\}.
$$
So clearly
\begin{enumerate}
    \item[$(*)$]   $\LL J^*_\nu : \nu\in J'\RR$ is a partition of $I \setminus J$ into $\leq\mu$-many sets.
\end{enumerate}
[Why? See clause (B) of the assumption.]

\mn
For $\eta\in I \setminus J$ let $\bfj(\eta) \defeq \max\{j : \eta \rest j \in J'\}$. 
It is well-defined (and $\leq\partial$) by our choice of $J'$ above, and clearly  $\eta \in J^*_{\eta\rest \bfj(\eta)}$.

We now observe:
\mn
\begin{enumerate}
    \item[$\otimes$]  If $n<\omega$, $\bar{\eta}' = \LL\eta'_\ell : \ell < n\RR$, 
    $\bar\eta'' = \LL \eta''_\ell : \ell < n \RR$, and $\eta'_\ell$, 
    $\eta''_\ell\in I$, then the following five clauses are a sufficient condition for 
    $$
    \tp_\qf(\bar\eta',J,I) = \tp_\qf(\bar{\eta}'',J,I)\text{:}
    $$
    \begin{enumerate}
        \item If $\eta'_\ell\in J$ or $\eta''_\ell\in J$ \underline{then} $\eta'_\ell = \eta''_\ell$.
\sn
        \item $\lh(\eta'_\ell) = \lh(\eta''_\ell)$
\sn
        \item If $\eta'_\ell \notin J$ (equivalently, $\eta''_\ell\notin J$) \underline{then}  $\bfj(\eta'_\ell) = \bfj(\eta''_\ell)$ --- call it $\bfj_\ell$ --- and $\eta'_\ell \rest \bfj_\ell=\eta''_\ell\rest \bfj_\ell$.
\sn
        \item for $\ell_1,\ell_2 < n < \omega$ and $j < \partial$, we have
        \begin{enumerate}
            \item[$(\alpha)$] $\eta'_{\ell_1} \rest j = \eta'_{\ell_2} \rest j \Leftrightarrow \eta''_{\ell_1}\rest j = \eta''_{\ell_2}\rest j$ \textbf{[This Follows]}
\sn
            \item[$(\beta)$] If the conditions in $(\alpha)$ are true, 
            $j < \min\!\big\{\lh(\eta'_\ell),\lh(\eta'_{\ell_2}) \big\} $,\\
            $m_1, m_2 < h(j)$, and for $\iota=1,2$ we have 
            \begin{align*}
                \big[& j+1 < \lh(\eta'_{\ell_\iota}) \wedge t'_\iota= \eta'_{\ell_\iota}(j)(m_\iota) \wedge t''_\iota = \eta''_{\ell_\iota}(j)(m_\iota) \big]\\
                \vee\ \big[& j+1 = \lh(\eta'_{\ell_\iota}) \wedge  t'_\iota = \eta'_{\ell_\iota}(j) \wedge t''_\iota = \eta''_{\ell_\iota}(j), \big],
            \end{align*}
            \underline{then}
            \begin{enumerate}[$\bullet_1$]
                \item $(\eta'_{\ell_1} \rest j) \caret \LL t'_1\RR <^I_1
            (\eta'_{\ell_1} \rest j) \caret \LL t'_2\RR \Leftrightarrow$\\
            $(\eta''_{\ell_1} \rest j) \caret \LL t''_1\RR 
            <^I_1 (\eta''_{\ell_1} \rest j) \caret \LL t''_2\RR$
\sn
                \item If $(\eta_{\ell_1}' \rest j) \caret \LL t'_2\RR \in J' \vee (\eta_{\ell_1}' \rest j) \caret \LL t''_1\RR \in J'$ \underline{then} they are equal.
            \end{enumerate}

        \end{enumerate}
\sn
        \item If $(\alpha)$ \underline{then} $(\beta)$, where:
        \begin{enumerate}
            \item[$(\alpha)$] First, $\eta'_\ell\in J^*_\nu$, $\eta'_\ell \rest j\in J'$, and $\eta'_\ell\rest (j+1) \notin J'$ (hence similarly for $\eta''_\ell$). 
            
            Second, $\nu \lhd \rho\in J$ and $m_1,m_2 < h(\lh(\nu))$.
            
            Third, for $\iota = 1,2$ we have $\bullet_1^\iota \vee \bullet_2^\iota$, where
            \begin{enumerate}[$\bullet_1$]
                \item $j+1 < \lh(\eta'_\ell)\wedge t' = \eta'_\ell(\bfj_\ell)(m_1) \wedge t'' = \eta''_\ell(\bfj_\ell)(m_1)$
\sn
                \item $j+1 = \lh(\eta'_\ell)\wedge t' = \eta'_\ell(j)\wedge t'' = \eta''_\ell(j)$.
            \end{enumerate}
\sn
            \item[$(\beta)$] $\nu \caret \LL t_{\ell_1}'\RR <^I_1 \nu \caret \LL t_{\ell_2}'\RR\ \Leftrightarrow\ \nu \caret \LL t_{\ell_1}''\RR <^I_1 \nu \caret \LL t_{\ell_2}''\RR$
\sn
        \end{enumerate}            
    \end{enumerate} 
\end{enumerate}
It is easy to check that this is true. Also, $\otimes$ defines the
equivalence relation (equality of quantifier-free types in $I$ over $J$) as various
pieces of information being the same. Now in all cases we have $\leq\mu$
choices (for clauses (d),(e) in $\otimes$, recall clause (A) in the
assumption) {so} we are done.

\mn
2) Similarly.

\mn
3) Follows, as well-orders are atomically $\mu$-stable.

\mn
4) Straightforward.
\end{PROOF}

\medskip
\begin{claim}\label{1.10}
$1)$ If $I\in K^\partial_{\tr(h)}$ is standard and $\lambda$ satisfies 
$(\forall\alpha < \lambda)\big[ |\alpha|^\partial < \lambda \big]$
\underline{then} $I$ is $\lepref{\lambda}$-atomically stable.

\mn
$2)$ In \emph{\ref{1.8}(1)} we can waive the assumption `$\mu \geq \partial$,' \underline{if} we 
  weaken  
clause \emph{\ref{1.8}(1)(B)} to 
\begin{enumerate}
    \item [$(B)^\supminus$] For all $I' \in [I]^{\leq\mu}$, the set
$$
    \big\{\eta \in I : \lh(\eta) \text{ is limit, and } \big(\exists i_* < \lh(\eta)\big) \big(\forall i > i_* \big) \big(\exists\ell < h(i) \big) \big[\Res^\ell_i(\eta) \in I'\big] \big\}
$$
    has cardinality $\leq\mu$.
\end{enumerate}

\end{claim}

\begin{PROOF}{\ref{1.10}}
1) Obvious, because well-orders are automatically $\lepref{\lambda}$-atomically stable; by \ref{1.8}(1) for $\lambda$ successor and by 
\ref{1.8}(2) for $\lambda$ a limit cardinal.

\mn
2) Easy (but this result will not be used later).
\end{PROOF}

\bigskip
\subsection{Existence Proofs} {}\ 

\medskip
In this subsection we first generalize the simple Black Box (B.B.~for short; see \cite[Lemma 1.5\subref{4.5A}]{Sh:309}) to our context, and then we prove \ref{1.7}(and more, as promised earlier) on the existence of full super-bigness properties.

\sn
\begin{lemma}\label{1.12}
\textbf{\emph{The $h$-fold simple B.B. Lemma.}}

Assume $\lambda,\partial \geq \aleph_2$, 
$h : \partial \to \omega \setminus \{0\}$, $I \in K_{\tr(h)}^\partial$ 
is as in Definition \emph{\ref{1.2}} for $J \defeq (\lambda,<)$ 
(that is, in \emph{\ref{1.2}(1)(A)} the set of elements of $I$ is equal to $\set_{\tr(h)}(\lambda)$),\footnote{
    Recalling $P_\partial^I = \{\eta \in I : \lh(\eta) = \partial\}$, $P_{<\partial}^I = \bigcup\limits_{i<\partial} P_i^I$, and $\eta(i) \in \inc_{h(i)}(J)$.
} 
and let $S \defeq \clH_{<\aleph_0}(\lambda)$.

\sn
$1)$ There are functions $f_\eta$ for $\eta \in P_\partial^I$, and pairwise disjoint 
$Y_\eps \subseteq P_\partial^I$ for $\eps < \lambda$ such that:
\begin{enumerate}[$(i)$]
    \item $\dom(f_\eta) = \big\{ (\Res_i^\ell)^I(\eta) : i < \partial,\ \ell < h(i) \big\}$. That is, 
    $$
    \dom(f_\eta) \defeq \big\{ \eta \rest j : j < i \text{ not successor} \big\} \cup \big\{ (\eta \rest j) \caret \big\LL \eta(j)(\ell) \big\RR : j + 1 < i,\ \ell < h(j) \big\}.
    $$

    \item $\rang(f_\eta) \subseteq S$
\sn
    \item If $f$ is a function from $P_{<\partial}^I$ into $S$, $g$ is a function from $P_{<\partial}^I$ into some $\gamma < \lambda$, and $\eps < \lambda$, then for some 
    $\eta \in Y_\eps \subseteq P_\partial^I$, we have: 
    \begin{enumerate}[$\bullet_1$]
        \item $f_\eta \subseteq f$
\sn
        \item $g \rest \{(\Res_i^\ell)^I(\eta) : \ell < h(i)\}$ is constant for each $i < \partial$.
    \end{enumerate}
\end{enumerate}

\sn
$2)$ Assume in addition that $\lambda \to \big(h(i) \big)_{\theta_i}^{\leq m_i}$ for $i < \partial$, where $m_i \leq h(i)$. \underline{Then} we can replace part $(1)$ with
\begin{enumerate}[$(i)_2$]
    \item $\dom(f_\eta) \defeq \big\{ (\eta \rest i) \caret \LL\eta(i) \rest u\RR : i < \partial,\ u \in \big[ h(i) \big]^{\leq m_i}\big\}.$
\sn
    \item $\rang(f_\eta) \subseteq S$
\sn
    \item ``If $(a)$ then $(b)$," where
    \begin{enumerate}[$(a)$]
        \item 
        \begin{enumerate}
            \item $f : P_{<\lambda}^I \to S$ 
\sn
            \item $g = \bigcup\limits_{i<\partial} g_i$, where $g_i$ is a function into 
            $\theta_i$, with domain 
            $$
            \bigcup\limits_{\rho \in P_{i+1}^I} \Big[ \big\{ (\rho \rest i) \caret \LL \alpha\RR : \alpha < \lambda \big\} \Big]^{\leq m_i}.
            $$
            \item $\eps < \lambda$
        \end{enumerate}
\sn
        \item For some $\eta \in Y_\eps \subseteq P_\partial^I$, we have:
        \begin{enumerate}
            \item $f_\eta \subseteq f$
\sn
            \item $g \rest \big[ \{(\Res_i^\ell)^I(\eta) : \ell < h(i)\} \big]^k$ is constant for each $i < \partial$ and $k \leq m_i$.
        \end{enumerate}
    \end{enumerate}
\end{enumerate}
\end{lemma}

\sn
\begin{remark}\label{1.13}
1) {Quoting} \ref{1.12} in \cite[Th 3.14, Def 3.13]{Sh:1232}, note that:

\begin{enumerate}[$(a)$]
    \item $\kappa,\lambda$ there correspond to $\aleph_0$ and $\lambda$ here.
\sn
    \item $\Lambda_{<\omega},\Lambda_\omega$ there correspond to $P_{<\partial}^I,P_\partial^I$ here.
\sn
    \item $g_\eta,g,f,\lambda$ there correspond to $f_\eta,f,g,\gamma$ here.
\end{enumerate}

\sn
2) We can allow finite $\lambda$, but then we would have to add the condition 
$$
i < \partial \Rightarrow (h(i)-1) \cdot \gamma < \lambda.
$$

\sn
3) Clearly we can say something about the case $\lambda = \aleph_1$ (e.g.~if we have guessing of clubs).

\sn
4) We can have $\LL Y_\eps : \eps < \lambda^\partial\RR$ and 
$|S| = \lambda^\partial$ by using a partition 
$\LL W_j : j < \partial\RR$ of $\partial$ into sets of cardinality 
$\partial$, and using $\LL \eta(i)(0) : i \in W_j\RR$ to code one member of $\lambda^\partial$ (or of $\clH_{<\partial}(\lambda)$).

\sn
5) In \ref{1.12}(2)$(iii)_2\,\bullet_2$, use $\LL g_{i,m} : i < \partial,\ m \leq m_i\RR$ with $g_{i,m} : [\lambda]^m \to \theta_{i,m}$, where $\beth_{m-1}(\theta_{i,m}) < \lambda$. (Similarly in \ref{1.15}.)
\end{remark}

\sn
\begin{PROOF}{\ref{1.12}}
1) Let $\LL W_{\bar s,\eps} : \bar s \in {}^{\omega>}S,\ \eps < \lambda\RR$ be a partition of $\lambda$
into $|{}^{\omega>}S \times \lambda|$-many sets, each of cardinality $\lambda$.
For $i \leq \partial$, let $\Lambda_i \defeq \{\eta \rest i : \eta \in P_\partial^I\}$, and choose $\LL f_\eta : \eta \in \Lambda_i\RR$ by induction on $i$ such that:
\begin{enumerate}
    \item[$(*)_1^i$] 
    \begin{enumerate}
        \item $\dom(f_\eta) \defeq$ 
        $$
        \big\{ \eta \rest j : j < i \text{ not successor} \big\} \cup \big\{ (\eta \rest j) \caret \big\LL \eta(j)(\ell) \big\RR : j + 1 < i,\ \ell < h(j) \big\}
        $$

        \item $\rang(f_\eta) \subseteq S$
\sn       
        \item If $\nu \lhd \eta$ then $f_\nu \subseteq f_\eta$.
\sn        
        \item If $j + 1 < i$ and $\ell < h(j)$ then 
        $$
        f_\eta\big( (\eta \rest j) \caret \LL \eta(j)(\ell)\RR \big) \defeq 
        \begin{cases}
            s_\ell &\text{if } \eta(j+2)(0) \in W_{\bar s,\eps} \text{ and } \ell < h(i)\\
            0 &\text{otherwise,}
        \end{cases}
        $$ 
        
        \item If $j=0$ or $j$ is limit $<i$, and $\eta(j)(0) \in W_{\bar s,\eps}$ for some $\bar s$, then $s_0$ will be $f_\eta(\eta \rest j)$ (or zero if $\varrho = \LL\ \RR$). 
    \end{enumerate}
\sn
    \item [$(*)_2$] For $\eps < \lambda$, let $Y_\eps \defeq \big\{\eta \in P_\partial^I : j < \partial \Rightarrow  \eta(j)(0) \in \bigcup\limits_{\bar s \in S} W_{\bar s,\eps} \big\}$.
\end{enumerate}
So $\bar f_i = \LL f_\eta : \eta \in P_i^I\RR$ is well-defined for $i \leq \kappa$, and obviously $\bar f_\partial$ satisfies clauses $(i),(ii)$ of the desired conclusion. What about clause $(iii)$?

{Fix $\eps$. Assume} $f : P_{<\kappa}^I \to S$ and 
$g : P_{<\kappa}^I \to \gamma$ for some $\gamma < \lambda$. We choose 
$\eta_i \in \Lambda_i$ by induction on $i$ to be a $\lhd$-increasing continuous seequence.

If $i = 0$ or $i$ is limit, we have no freedom. 

If $i = j + 1$ and $j$ is not a successor ordinal, then let $\bar s \defeq \big\LL f(\eta \rest j)\big\RR$ (so $\bar s \in {}^{\omega>} \clH_{<\aleph_0}(\lambda)$), and choose $\rho \in \inc_{h(j)}(W_{\bar s,\eps})$. Let 
$$
Y_\eps \defeq \big\{ \eta \in P_n^I : \eta(i)(\ell) \in \textstyle\bigcup\limits_{\bar s} W_{\bar s,\eps} \text{ for every $i < \partial$ and } \ell < h(i) \big\}
$$
and $\eta_i \defeq \eta_j \caret \LL\rho\RR$.

Lastly, if $i = j + 1$ and $j$ is a successor ordinal, then let 
$$
\bar s \defeq \big\LL f( \Res_{j-1}^\ell(\eta)) : \ell < h(j-1) \big\RR
$$ 
and choose $\rho \in \inc_{h(j-1)}(W_{\bar s,\eps})$ such that $g \rest \big\{ \eta \caret \LL \rho(\ell)\RR : \ell < h(1)\big\}$ is constant 
(this is possible because $W_{\bar s,\eps}$ is large enough). Let 
$\eta_i \defeq \eta_j \caret \LL\rho\RR$.

Now it is easy to check that $\eta$ satisfies clause $(iii)$.

\sn
2) Similarly.
\end{PROOF}

\mn
Next we fulfill a promise.

\begin{lemma}\label{1.15}
$K^\partial_\ptr$ and $K^\partial_{\tr(h)}$ (for 
$h \in {}^\partial(\partial\setminus \{0\})$) have the full 
$(\lambda,\lambda,\mu,\pinkappa)$-super-$\bar\bfe$-bigness
property as witnessed by  
an 
atomically $\Theta$-stable structure \underline{when} at least one of the following holds:

\mn
\textbf{\emph{Case 1:}} 

$\lambda$ is regular, 
$\lambda = \lambda^\partial> \mu_1 = \mu_2 = \mu$, and for all $i < \partial$ large enough we have 
$$
\lambda \to \big( h(i) \big)^{\max(\bfe_i)}_{\mu_*},
$$
where $\mu_* \defeq \mu^{<\kappa+\mu^\partial}$, $\max(\bfe_i)$ is as defined in \emph{\ref{1.4}(3)}, and $\Theta \defeq \{\sigma : \sigma^\partial = \sigma\}$.

\mn
\textbf{\emph{Case 2:}}

$\mu_1 = \mu_2 = \mu$, $\lambda = \cf(\lambda) > \theta = \theta^\partial = \theta^{<\kappa} > \mu > \kappa$, $\mu_*$ is as above, 
$$
\theta \to \big( h(i) \big)^{\max(\bfe_i)}_{\mu_*}
$$
for all $i < \partial$ large enough, and $\Theta \defeq [\theta,\infty)_\Card$.

\mn
\textbf{\emph{Case 3:}}

$\lambda$ is singular, $\mu_2 = \mu_2^\partial + \mu_2^{<\kappa}  < \lambda \leq 2^{\mu_2}$, 
$\kappa \leq \partial$,  
$$
\mu_2^+ \to \big( h(i) \big)^{\max(\bfe_i)}_{\mu_2}
$$
for all $i < \partial$ large enough, 
$$
\mu_1 \defeq 
\begin{cases}
    \mu_2 &\text{if } \partial = \aleph_0 \text{ and } \bigwedge\limits_i[\max(\bfe_i) = 1]\\
    2^\partial + 2^{<\kappa} &\text{otherwise}
\end{cases}
$$
and $\Theta \defeq [\mu_2^+,\infty)_\Card$.

\mn
\textbf{\emph{Case 4:}}

Like Case 2, but $\lambda$ is singular.
\end{lemma}

\mn
\begin{conclusion}\label{1.18}
In all four cases in \emph{\ref{1.15}}, we can deduce that $K^\partial_{\tr(h)}$ has the full 
$(\lambda,\lambda,\mu,\kappa)$-$\psi_{\tr(h)}$-bigness property. 
\end{conclusion}

\sn
\begin{remark}\label{1.21}
1) The reader may use only {Case 1 of} {\ref{1.15}}, because it suffices for many cardinals.

\mn
2) We may consider $M_i,N_i$ of cardinality $\mu_i$ (or $\theta_i$).

\mn
3) Similarly, we can use an increasing sequence $\bar\theta = \LL \theta_i : i < \partial\RR$ in Case 2.
 
\mn
4) We can weaken ``$\lambda \to \big( h(i) \big)^{\max(\bfe_i)}_{\mu_*}$'' as follows.
\begin{quotation}
    `If $\bfc_k : [\theta]^k \to \mu_k$ for $k \leq m \defeq \max(\bfe_i)$ \underline{then} there exists $u \in [\theta]^{h(i)}$ such that $\bfc_k \rest [u]^k$ is constant for all $k \leq m$.'
\end{quotation}
So $\mu_k \defeq \beth_{m-k}(\mu)$ will suffice by the proof of the Erd\H{o}s-Rado Theorem.

\mn
5) We can let $\partial$ be singular, but this is of doubtful interest.
\end{remark}

\begin{PROOF}{\ref{1.15}}
\textbf{Proof of Lemma \ref{1.15}.}

The proof splits into cases.

\bn
\centerline{\textbf{\underline{Case 1:}}}

First, as $\lambda > \partial$ are regular and $\lambda = \lambda^\partial$ (hence $\lambda = \lambda^{<\partial}$):
\begin{enumerate}
    \item [$(*)_0$] We can find a stationary subset $S_* \subseteq \{\delta < \lambda : \cf(\delta) = \partial\}$ which belongs to $\check I_\partial[\lambda]$.
\end{enumerate}
By this we mean we can choose $\olsi C = \LL C_\alpha : \alpha < \lambda\RR$ such that:
\begin{itemize}
    \item $C_\alpha \subseteq \alpha$
\sn
    \item $\otp(C_\alpha) \leq \partial$
\sn
    \item $\beta \in C_\alpha \Rightarrow C_\beta = C_\alpha \cap \beta$
\sn
    \item $\alpha \in S_* \Rightarrow \alpha = \sup (C_\alpha)$ (hence $\otp(C_\alpha) = \partial$).
\end{itemize}




\begin{enumerate}[$(*)_1$]
    \item We can choose a partition $\olsi S = \LL S_\zeta : \zeta < \lambda\RR$ of $S_*$ into stationary sets.
    \begin{enumerate}
        \item Without loss of generality, for each $\zeta < \lambda$, $\olsi C^\zeta \defeq \olsi C \rest S_\zeta$ guesses clubs (in the following weak sense):
        
\sn        
        \item For every club $E$ of $\lambda$, for stationarily many $\delta \in S_\zeta$, we have
        $$
        \alpha \in C_\delta \Rightarrow \min \!\big(E \setminus (\alpha+1)\big) < \min \!\big(C_\delta \setminus (\alpha+1)\big).
        $$
    \end{enumerate}
    \hspace{-1.6cm}Next,
    
\sn
    \item For $\zeta < \lambda$ and $\delta \in S_\zeta$, let $\rho_\delta,\nu_\delta \in 
    {}^\partial\!\lambda$ be defined by:
    \begin{enumerate}
        \item $\rho_\delta(i)$ is the $(2i+1)^\tthh$ member of $C_\delta$.
\sn
        \item $\nu_\delta(i)$ is the $(2i+2)^\tthh$ member of $C_\delta$.
    \end{enumerate}
\sn
    \item 
    \begin{enumerate}
        \item Let $\Lambda_\delta \defeq \big\{ \eta \in {}^\partial\delta : i < \partial \Rightarrow \eta(i) \in \inc_{h(i)}\!\big( [\rho_\delta(i),\nu_\delta(i)\big)\big) \big\}$ for 
        $\delta \in S$.
\sn
        \item Let $I_{\tr(h)}^\lambda \in K_{\tr(h)}^\omega$ be as in \ref{1.2}(1), with set of elements $\set_{\tr(h)}(\lambda)$.\footnote{
            See \ref{1.2}(1)(A).
        }
\sn
        \item Let $I_\zeta$ be the submodel of $I_{\tr(h)}^\lambda$ with set of elements 
        $$
        \bigcup_{\delta \in S_\zeta} \Lambda_\delta \cup P_{<\partial}^{I_{\tr(h)}^\lambda}.
        $$
    \end{enumerate}
\sn
    \item 
    \begin{enumerate}
        \item We will show that $\bar I = \LL I_\zeta : \zeta < \lambda\RR$ exemplifies the conclusion; this will suffice. 
\sn
        \item So let $\zeta_* < \lambda$, $I \defeq I_{\zeta_*}$, 
        $J \defeq \sum\limits_{\zeta \neq \zeta_*}I_\zeta$, and let 
        $\chi,x,f_0,f_1,<_\chi^*$ be as in \ref{1.4}(1).
    \end{enumerate}
\sn
    \item It suffices to find $\big\LL (M_i,N_i) : i < \partial\big\RR$ as in clauses $(i)$-$(vii)$ of $(*)^{I,J}_{\lambda,\mu,\partial}$ in Definition \ref{1.4}(1).

\sn
\sn
    \item[$(*)_{5.1}$] We can choose $\clB_\alpha^*$, $M_\alpha^*$ by induction on $\alpha < \lambda$ such that
    \begin{enumerate}
        \item $\clB_\alpha^*,M_\alpha^* \prec (\clH(\chi),\in,<_\chi^*)$, and $\clB_\alpha^*$ is $\prec$-increasing continuous in $\alpha$.
\sn
        \item $\|\clB_\alpha^*\| < \lambda$ and\footnote{
            We could choose $\|M_\alpha^*\| = \mu$, but there would be no gain in doing so.
        } $\|M_\alpha^*\| = \mu_*$.
\sn
        \item $[M_\alpha^*]^{\leq\partial} \subseteq M_\alpha^*$, $[M_\alpha^*]^{<\kappa} \subseteq M_\alpha^*$, and $\mu_* + 1 \subseteq M_\alpha^*$.
\sn
        \item $\LL\clB_\gamma,M_\gamma^* : \gamma \leq \beta\RR \in \clB_\alpha^*$ for all $\beta < \alpha$.
\sn
        \item If $\beta \in C_\alpha$ \underline{then} $\big\LL M_\gamma^* : \gamma \in C_\alpha \cap (\beta+1) \big\RR,\,\LL\clB_\gamma^* : \gamma \leq \beta\RR \in M_\alpha^*$.
\sn
        \item $x,f_0,f_1$, $\olsi S$, $\big\LL \olsi C^\zeta : \zeta < \lambda \big\RR$ (and hence $\big\LL(\eta_\delta,\nu_\delta) : \delta \in S_\zeta\big\RR$), $I,J$, and $\zeta_*$ all belong to 
        $\clB_\alpha^*$ {and} to $M_\alpha^*$.  
    \end{enumerate}
\end{enumerate}
[Why? Obvious.]

\sn
\begin{enumerate}
    \item[$(*)_{5.2}$]
    \begin{enumerate}
        \item The set $E \defeq \{\delta < \lambda : \delta \text{ is limit, } \clB_\delta^* \cap \lambda = \delta\}$ is a club of 
        $\lambda$.
\sn
        \item We choose $\delta_* \in E \cap S_{\zeta_*}$ such that $C_{\delta_*} \subseteq E$.
    \end{enumerate}
\end{enumerate}
[Why? Clause (a) is obvious, and (b) follows by our choice of $\olsi C^{\zeta_*}$.]

\sn
\begin{enumerate}
    \item[$(*)_{5.3}$] For $i < \delta_*$, we choose $M_i \defeq M_{\rho_{\delta_*}\!(i)}^*$ and $N_i \defeq N_{\nu_{\delta_*}\!(i)}^*$.
\sn
    \item[$(*)_{5.4}$] The set $\{\tp_\qf(\bar \nu,Y_i, t_{\zeta_*}) : \bar \nu \in {}^{\kappa>}\!I\}$ has cardinality $\leq\mu_*$,   
    {where}\\ $Y_i \defeq M_i \cap P_{<\delta}^I$. 
\end{enumerate}
[Why? Because in Case 1 we are assuming $\alpha < \lambda \Rightarrow |\alpha|^\partial < \lambda$.]

\medskip
It suffices to prove that $\LL M_i,N_i : i < \partial\RR$ is as required.
That is, we have to check clauses $(i)$-$(vii)$ of \ref{1.4}$(*)_{\lambda,\bar\mu}^{I,J}$, and then we will be done with Case 1.

\mn
\textbf{Clauses ($i$)-($iv$):} Easy.

\mn
\textbf{Clause ($v$):}

We have to find an appropriate $\eta \in P_\partial^I$. (In fact, it will be a member of $\Lambda_{\delta_*}$.)

For this, we will choose $\eta_i$ by induction on $i$ such that
\begin{enumerate}
    \item[$(*)_{5.5}$]  
    \begin{enumerate}
        \item $\eta_i \in M_i$ and $\lh(\eta_i) = i$. 
\sn
        \item If $j < i$ then
        \begin{enumerate}
            \item $\eta_i(j) \in \inc_{h(j)}(\big(\rho_{\delta_*}\!(j),\nu_{\delta_*}\!(j)\big))$
\sn
            \item $(\eta_i \rest j) \caret \big\LL\eta_i(j)(0)\big\RR$, $(\eta_i \rest j) \caret \big\LL\eta_i(j)(1)\big\RR$,  . . .\\ $(\eta_i \rest j) \caret \big\LL \eta_i(j)(h(j)-1) \big\RR$ all realize the same Dedekind cut (by $<_j^I$) on $I \cap M_j \cap \{(\eta_i \rest j) \caret \LL s\RR : s \in J\}$.
\sn
            \item $(\eta_i \rest j) \caret \big\LL\eta_i(j)(\ell)\big\RR \in I \defeq I_{\zeta_*}$.
        \end{enumerate}
    \end{enumerate}
\end{enumerate}
[Why is this possible? The $i=0$ case is trivial. 
For $i$ limit, let $\eta_i \defeq \bigcup\limits_{j<i}\eta_j$; it is of the correct form and belongs to $N_i$ by clause $(*)_{5.4}$. If $i = j+1$ \underline{then}, as $\eta_j \in N_j$ and ${}^\partial\!M_i \in N_i$, we can choose $\eta_i(j)$ as in the proof of \ref{1.12}.]

\mn
\textbf{Clause ($vi$):} 

In the `$i = j+1$' case in the previous proof, we just have to add the relevant additional demand.

Letting\footnote{
    See \ref{1.4}(3)(E).
} 
$m \defeq \max(\bfe_j)$, recalling $f_0 : I \to \mu$ and $f_1 : I \to {}^{\kappa>}\!J$, we define a two-place relation $\clE$ such that:
\begin{enumerate}
    \item [$(*)_{5.6}$]  
    \begin{enumerate}[(A)]
        \item $\clE \subseteq [\lambda]^{\leq m} \times [\lambda]^{\leq m}$
\sn
        \item $u_1\ \clE\ u_2$ \underline{iff}
        \begin{enumerate}[(a)]
            \item $|u_1| = |u_2|$
\sn
            \item If $\alpha_1 \in u_1$, $\alpha_2 \in u_2$, and 
            $|u_1 \cap \alpha_1| = |u_2 \cap \alpha_2|$ (and hence 
            $\eta_j \caret \LL\alpha_1\RR, \eta_j \caret \LL\alpha_2\RR \in I$), \underline{then}
            $$
            f_0(\eta_j \caret \LL\alpha_1\RR) = f_0(\eta_j \caret \LL\alpha_2\RR) < \mu.
            $$

            \item Letting $\varrho_1,\varrho_2$ list $u_1,u_2$ in increasing order (and hence 
            $\eta_j \caret \LL\varrho_1\RR$, $\eta_j \caret \LL\varrho_2\RR \in \sseq_j(I)$),\footnote{
                See \ref{1.2}(1B).
            }
            \underline{then} 
            \begin{enumerate}[$\bullet_1$]
                \item $f_1(\eta_j \caret \LL\varrho_1\RR), f_1(\eta_j \caret \LL\varrho_2\RR) \in {}^{\kappa>}\!J$ are of the same length ($<\kappa$) and realize the same quantifier-free type over $J \cap M_{j+1}$ inside $J$.
\sn
                \item $f_0(f_1(\eta_j \caret \LL\varrho_1\RR)) = f_0(f_1(\eta_j \caret \LL\varrho_2\RR))$
            \end{enumerate}
        \end{enumerate}
    \end{enumerate}
\sn
    \item[$(*)_{5.7}$]  
    \begin{enumerate}
        \item $\clE$ is an equivalence relation on $[\lambda]^{\leq m}$.
\sn
        \item $\clE$ has $\leq \|M_j\|^{<\kappa} + \|M_j\|^\partial = \mu_*$ equivalence classes.
    \end{enumerate}
\end{enumerate}
[Why? For clause (a) one can straightforwardly check the definition, and clause (b) holds because $\|M_\partial\| = \mu_*$ (recalling $(*)_{5.4}$).]

\mn
As we are assuming $\lambda \to \big(h(i) \big)_{\mu_*}^{\leq m}$ and as $\eta_j \in M_{\nu_{\delta_*}\!(j)} \prec (\clH(\chi),\in,<_\chi^*)$,
\begin{enumerate}
    \item [$(*)_{5.8}$] For the coloring $\bfc$, it is enough to find a sequence $\varrho \in \inc_{h(j)}(\lambda)$ such that
    \begin{enumerate}
        \item Every $\{\rho(\ell) : \ell < h(j)\}$ realizes the same Dedekind cut of $(\lambda,<) \rest M_j$ in $(\lambda,<)$.
\sn
        \item If $m \leq \max(\bfe_j)$ and $u_1,u_2 \in [h(j)]^m$, then 
        $$
        \bfc(\{\rho(\ell) : \ell \in u_1\}) = \bfc(\{\rho(\ell) : \ell \in u_2\}).
        $$
    \end{enumerate} 
\end{enumerate}

As in the proof of \ref{1.12}, there is $Y \in [\lambda]^\lambda$ such that every $\varrho \in \inc_{h(j)}(Y)$ satisfies $(*)_{5.8}$(a). As we have assumed $\lambda \to \big(h(i) \big)_{\mu_*}^{\leq m}$, we can choose $\varrho$ satisfying clause (b), so $(*)_{5.8}$ holds.

Now we have finished the proof of Case 1.

\bn
\centerline{\textbf{\underline{Case 2:}}}

We repeat the proof of Case 1, with some changes.
First, we can find a stationary $S_* \in \check I_\partial[\lambda]$ as in $(*)_0$. 

\sn
[Why? For stationarity,
we cannot use $\lambda = \cf(\lambda) = \lambda^\partial + \lambda^{<\partial}$, but we can apply \cite{Sh:420} because $\lambda = \cf(\lambda) > \partial^+ > \partial = \cf(\partial)$, because $\lambda > \theta = \theta^\partial > \partial$.]

\medskip
Then as in $(*)_1$-$(*)_2$, we can find $\olsi C$, $\LL S_\zeta : \zeta < \lambda\RR$, $\delta$, and $\rho_\delta,\nu_\delta$ as in Case 1. 

Now
\begin{enumerate}
    \item [$\circledast_3$]
    \begin{enumerate}
        \item We choose $\gA_\alpha$ by induction on $\alpha < \lambda$ such that
        \begin{enumerate}
            \item $\gA_\alpha \prec (\clH(\chi),\in,<_\chi^*)$
\sn
            \item $\|\gA_\alpha\| = \theta$
\sn
            \item $[\gA_\alpha]^{\leq \partial} \subseteq \gA_\alpha$ (recalling $\theta = \theta^\partial$).
\sn
            \item $\theta+1 \subseteq \gA_\alpha$ and $\LL \gA_\beta : \beta \in C_\alpha\RR \in \gA_\alpha$
        \end{enumerate}
\sn
        \item For $\delta \in S$, let  
        $$
        \Lambda_\delta \defeq \big\{ \eta \in {}^\partial\delta \cap \gA_\delta : i < \partial \Rightarrow \eta(i) \in \inc_{h(i)}\!\big( [\rho_\delta(i),\nu_\delta(i)\big)\big) \big\}.
        $$
        \item For $\zeta < \lambda$, let $I_\zeta$ be as in Definition \ref{1.2}, with set of elements 
        $$
        \bigcup_{\delta \in S_\zeta} \Lambda_\delta \cup \big\{ \Res_i^\ell(\eta) : \eta \in \textstyle\bigcup\limits_{\delta \in S_\zeta} \Lambda_\delta,\ i < \partial,\ \ell < h(i) \big\}.
        $$
    \end{enumerate}
\end{enumerate}

\medskip 
We then follow the proof in Case 1 again, until the point where we settle clause $(v)$ of \ref{1.4}$(*)_{\lambda,\bar\mu}^{I,J}$ by constructing a sequence $\LL \eta_i : i < \partial\RR$.
As before, $\lh(\eta_i) = i$, but now $\eta_i \in N_i \cap \gA_{\rho_{\delta_*}\!(i)}$.

(Note that $\gA_{\rho_{\delta_*}\!(i)} \in N_i$ is of cardinality $\theta$, and thus large enough.)

The rest of the proof is similar to Case 1, but we choose 
$\eta \in M_i \cap \gA_{\rho_i(\delta)}$ and use 
``$|\gA_\alpha \cap \lambda| = \theta$, and $\theta$ is large enough."

\bn
\centerline{\textbf{\underline{Case 3:}}}

The $\partial > \aleph_0$ case is 
   quite 
different from earlier proofs: 
\cite[\S3, 2.7, p.116]{Sh:136}, \cite[1.15(2)\subref{7.6}]{Sh:331} (which uses \cite{Sh:117}; or see\cite[1.17(2)\subref{a48}]{Sh:E62}) ---which dealt with $\partial = \aleph_0$. 

We are assuming $\mu_2 = \mu_2^{<\partial} < \lambda < 2^{\mu_2}$. Now we can choose:\footnote{
    See \cite[3.12\subref{4.EK}]{Sh:E62} --- i.e.~Engelking-Kar\l owicz \cite{EK65}.
}  
\begin{enumerate}[$(*)_1$]
    \item  =$A_i \subseteq \mu_2$ for $i < 2^{\mu_2}$ ($i < \lambda$ would be sufficient) such that
    if $W \in [\mu_2]^{<\partial}$ and $i \in \mu_2 \setminus W$, \underline{then} $A_i \not\subseteq \bigcup\limits_{j \in W} A_j$.
\sn
    \item  Let $\LL S^\zeta : \zeta < \mu_2\RR$ be pairwise disjoint stationary subsets of 
    $$
    \{\delta < \mu_2^+ : \cf(\delta) = \partial\}.
    $$ 
    For $\xi < 2^{\mu_2}$, let $\bfS_\xi \defeq \bigcup\limits_{\zeta \in A_\xi} S^\zeta$.
\sn
    \item  For $\xi < 2^{\mu_2}$, let $I_\xi$ be as in \ref{1.2} with set of elements $\bigcup\limits_{\delta \in \bfS_\xi} \Lambda_{\xi,\delta}$, where $\Lambda_{\xi,\delta}$ is defined as follows.

\sn
    $\eta \in \Lambda_{\xi,\delta}$ \underline{iff}
    \begin{enumerate}
        \item $\eta \in \set_{\tr(h)}(\lambda)$ (Of course, this follows from $\eta \in I_\xi$.)
\sn
        \item If $i < j$, $\lh(\eta) > j+1$, $k < h(i)$ and 
        $\ell < h(j)$, \underline{then} $\eta(i)(k) < \eta(j)(\ell)$.
\sn
        \item If $\lh(\eta) = \partial$ \underline{then} 
        $\bigcup\limits_{i<\delta}\eta(i)(0) = \delta$.
\sn
        \item If $\lh(\eta) > i+1 \wedge \ell < h(i) \wedge \eta(i)(\ell) = \alpha$ \underline{or} $\lh(\eta) = i+1 \wedge \eta(i) = \alpha$, then $\alpha \in \bfS_\xi$.
    \end{enumerate}
\sn
    \item Let $\xi_*$, $I \defeq I_{\xi_*}$, $J \defeq J_{\xi_*}$, 
    $\chi,x,f_0,f_1,<_\chi^*$ be as in the proof of Case 1.
\sn
    \item We can find $\big\LL (M_i,N_i) : i < \partial \big\RR$ and $\delta \in \bfS_{\xi_*}$ such that:
    \begin{enumerate}
        \item $M_i \prec N_i \prec (\clH(\chi),\in,<_\chi^*)$
\sn
        \item $\|M_i\| = \|N_i\| = \mu_1$ 
\sn
        \item $[M_i]^{<\partial} \subseteq M_i$, $[N_i]^{<\partial} \subseteq N_i$.
\sn
        \item $N_i \prec M_j$ for all $i < j < \partial$.
\sn
        \item $\lambda,\chi,f_0,f_1,<_\chi^*$, $\LL S^\zeta : \zeta < \mu_2\RR$, $\LL \bfS_\xi : \xi < 2^{\mu_2}\RR$, and $\LL\gA_\delta : \delta \in S_*\RR$ all belong to $M_i$.
\sn
        \item $\LL \sup(M_i \cap \mu_2^+) : i < \partial\RR$ is increasing with limit $\delta$.
\sn
        \item $\big[ \rho_\delta(i),\nu_\delta(i) \big)$ is disjoint to $M_i$ and belongs to $N_i$.
    \end{enumerate}
\end{enumerate}
[Why? Easy.]
\begin{enumerate}
    \item [$(*)_{5.1}$] $\big\LL(M_i,N_i) : i < \partial \big\RR$ satisfies clauses $(*)_{\lambda,\bar\mu}^{I,J}(i)$-$(iv)$ of \ref{1.2}.
\end{enumerate}
[Why? Straightforward.]
\begin{enumerate}
    \item [$(*)_{5.2}$] We can find $\eta \in P_i^{I_\xi} \cap \Lambda_{\xi,\delta}$ as required in $(*)_{\lambda,\bar\mu}^{I,J}(v)$-$(vii)$.
\end{enumerate}

\mn
Towards this, we choose $\eta_i$ of length $i$ by induction on $i < 
    \partial $ 
such that:
\begin{enumerate}
    \item [$(*)_6^i$] For $j < i$,
    \begin{enumerate}
        \item $\eta_i(j)$ {belongs to} $\inc_{h(j)}(N_j \cap \mu_2^+ \setminus \sup(M_j \cap \mu_2^+))$ and to\\ 
        $\inc_{h(j)}\big[ \rho_\delta(i),\nu_\delta(i) \big)$.
\sn
        \item $(\eta_i \rest j) \caret \big\LL\eta_i(j)(0)\big\RR$, $(\eta_i \rest j) \caret \big\LL\eta_i(j)(1)\big\RR$,  . . . $(\eta_i \rest j) \caret \big\LL \eta_i(j)(h(j)-1) \big\RR$ all realize the same Dedekind cut (by $<_j^I$) on 
        $$
        I \cap M_j \cap \{(\eta_i \rest j) \caret \LL s\RR : s \in J\}.
        $$
        \item $(\eta_j \rest j) \caret \LL\eta_j(j)(\ell)\RR \in I \defeq I_{\zeta_*}$ for all $\ell < h(j)$.
\sn
        \item $\big\LL\eta_j(j)(\ell) : \ell < h(j) \big\RR$ is increasing.
\sn
        \item $\eta_j(\ell) \in \bfS_\xi \setminus \bigcup \{\bfS_{\xi'} : \xi' \in M_1 \cap 2^{\mu_2}\}$.
    \end{enumerate}
\end{enumerate}
The proof is straightforward 
Clause $(vi)$ is handled similarly.

\bn
\centerline{\textbf{\underline{Case 4:}}}

We continue the proofs of Cases 1 and 2. 

\begin{enumerate}
    \item [$\boxplus_1$] We choose $\bar\lambda$ and $\lambda_*$ such that
    \begin{enumerate}
        \item $\bar\lambda = \LL \lambda_\eps : \eps < \cf(\lambda)\RR$ is an increasing sequence of regular cardinals.
\sn
        \item $\bigcup\limits_{\eps < \cf(\lambda)} \lambda_\eps = \lambda$
\sn
        \item $\lambda_* = \cf(\lambda_*) \in \big[ \cf(\lambda) + \theta^+,\lambda_0 \big)$
    \end{enumerate}
\sn
    \item [$\boxplus_2$] Let us define $\lambda_{<\eps} \defeq \sum\limits_{\zeta<\eps} \lambda_\zeta$.
    \begin{enumerate}
        \item Choose   
          $ {S_  i 
          ^
          *} \in \check I_\partial[\lambda_i]$  stationary and disjoint to $\lambda_{<i} + \lambda_{<i} + 1$ such that 
        $$
        \delta \in S_\eps^* \Rightarrow \cf(\delta) \neq \partial.
        $$ 
        Without loss of generality, $\delta \in S_\eps^* \Rightarrow \lambda_* \divides \delta$.
\sn
        \item Choose $\olsi S^\eps = \big\LL S_\xi : \xi \in [\lambda_{<i},\lambda_i)\big\RR$ to be a partition of $S_i^*$ into stationary subsets of  
        $\lambda_i$.
    \end{enumerate}
\end{enumerate}
[Why is this possible? As $\lambda_\eps = \cf(\lambda_\eps) > \partial^+$, we can choose $S_\eps^*$ as in \cite
{Sh:420}  
and $\olsi S^\eps$ by the definition of $\check I_\partial[\lambda_\eps]$.
\begin{enumerate}
    \item[$\boxplus_3$] 
    \begin{enumerate}
        \item Choose $W_* \in  \check I_\partial[\lambda_*]$ stationary such that $\delta \in W_* \Rightarrow \cf(\delta) = \partial    
            \wedge \theta^+ \divides \delta $.
\sn
        \item We choose $\LL W_\eps : \eps < \cf(\lambda)\RR$ to be a partition of $W_*$ into stationary subsets of  
    $\lambda_*$.  
    \end{enumerate}
\sn
    \item[$\boxplus_4$] For $\eps < \cf(\lambda)$, we can find $\big\LL (S_\xi^+,\olsi C^\xi) : \xi \in [\lambda_{<\eps},\lambda_\eps)\big\RR$ such that
    \begin{enumerate}
        \item $\big\LL S_\xi^+ : \xi \in [\lambda_{<\eps},\lambda_\eps)\big\RR$ is a sequence of pairwise disjoint subsets of $S_\eps^+ \subseteq [\lambda_{<\eps},\lambda_\eps)$ such that $\delta \in S_\xi^+ \Rightarrow |\delta| \divides \delta$.
\sn
        \item $\olsi C^\xi = \LL C_\alpha^\xi : \alpha \in S_\xi^+\RR$ 
\sn
        \item $C_\alpha^\xi \subseteq S_\xi^+ \cap \alpha$
\sn
        \item $C_\beta^\xi = C_\alpha^\xi \cap \beta$ for all $\alpha \in S_\xi^+$ and $\beta \in C_\alpha^\xi$.
\sn
        \item $S_\xi = \{\alpha \in S_\xi^+ : \otp(C_\alpha^\xi) = \partial\}$ 
              hence $ |C_\alpha^\xi| = \partial $.)
    \end{enumerate}
\end{enumerate}
[Why? By the properties of $\check I_\partial[\lambda_\eps]$.]
\begin{enumerate}
    \item [$\boxplus_5$] We can find $\big\LL (W_\eps^+,\olsi C_\eps) : \eps < \cf(\lambda)\big\RR$ such that 
    \begin{enumerate}
        \item $\big\LL W_\xi^+ : \xi \in [\lambda_{<\eps},\lambda_\eps)\big\RR$ is a sequence of pairwise disjoint subsets of $W_*$, and $\delta \in W_\xi^+ \Rightarrow \theta^+ \divides \delta$.
\sn
        \item $\olsi C_\eps = \LL C_{\eps,\alpha} : \alpha \in W_\eps^+\RR$
\sn
        \item $C_{\eps,\alpha} \subseteq W_\xi^+ \cap \alpha$
\sn
        \item $C_{\eps,\beta} = C_\alpha^\xi \cap \beta$ for all $\alpha \in W_\eps^+$ and $\beta \in C_{\eps,\alpha}$.
\sn
        \item $W_\eps^* = \{\alpha \in W_\xi^+ : \otp(C_{\eps,\alpha}) = \partial\}$ 
            hence $|C_{\eps,\alpha}| = \partial$.
    \end{enumerate}
\sn
    \item[$\boxplus_6$] For $\eps < \cf(\lambda)$ and $\xi \in [\lambda_{<\eps},\lambda_\eps)$, let
    \begin{enumerate}
        \item $\clU_\xi^+ \defeq \big\{ \alpha+\beta : \alpha \in S_\xi^+,\ \beta \in W_\eps^+,\ \otp(C_\alpha^\xi) = \otp(C_{\eps,\beta}) \big\}$
\sn
        \item If $\gamma \in \clU_\xi^+$ (so it is of the form $\alpha+\beta$ as above) \underline{then} let  
        $$
        C_{\xi,\gamma}^\bullet \defeq \big\{ \alpha_1'+\alpha_2' : \alpha_1' \in C_\alpha^\xi,\ \alpha_2' \in C_{\eps,\beta},\ \otp(C_{\alpha_1'}^\xi) = \otp(C_{\eps,\alpha_2'}) \big\}.
        $$
    \end{enumerate}
\sn    
    \item[$\boxplus_7$] For $\eps< \cf(\lambda)$ and $\xi \in [\lambda_{<\eps},\lambda_\eps)$, we choose $\gA_{\xi,\gamma}$ by induction on $\gamma \in S_\xi^\bullet$ such that
    \begin{enumerate}[$\bullet_1$]
        \item $\gA_{\xi,\gamma} \prec (\clH(\chi),\in,<_\chi^*)$
\sn
        \item $\|\gA_{\xi,\gamma}\| = \theta$
\sn
        \item $[\gA_{\xi,\gamma}]^{\leq \partial} \subseteq \gA_{\xi,\gamma}$ (recalling $\theta = \theta^\partial$).
\sn
        \item $\theta+1 \subseteq \gA_{\xi,\gamma}$ and $\LL \gA_\beta : \beta \in C_{\xi,\gamma}^\bullet\RR \in \gA_{\xi,\gamma}$.
    \end{enumerate}
\sn
    \item[$\boxplus_8$] For $\eps< \cf(\lambda)$, $\xi \in [\lambda_{<\eps},\lambda_\eps)$, and $\gamma \in S_\xi^\bullet$, let $\rho_\gamma, \nu_\gamma \in {}^\partial\!\lambda$ be defined as follows.
    \begin{enumerate}
        \item $\rho_\gamma(i)$ is the $(2i+1)^\tthh$ member of $C_{\xi,\gamma}^\bullet$.
\sn
        \item $\nu_\gamma(i)$ is the $(2i+2)^\tthh$ member of $C_{\xi,\gamma}^\bullet$.
    \end{enumerate}
\sn
    \item[$\boxplus_9$] For $\eps< \cf(\lambda)$ and $\xi \in [\lambda_{<\eps},\lambda_\eps)$:
    \begin{enumerate}
        \item For $\gamma \in S_\xi^\bullet$, let
        $$
        \Lambda_\gamma \defeq \big\{ \eta \in {}^\partial\gamma \cap \gA_{\xi,\gamma} : i < \partial \Rightarrow \eta(i) \in \inc_{h(i)}\!\big( [\rho_\gamma(i),\nu_\gamma(i)\big)\big) \big\}.
        $$
        \item let $I_\xi$ be as in Definition \ref{1.2}, with set of elements 
        $$
        \bigcup_{\gamma \in S_\xi^\bullet} \Lambda_\gamma \cup \big\{ \Res_i^\ell(\eta) : \eta \in \textstyle\bigcup\limits_{\gamma \in S_\xi^\bullet} \Lambda_\gamma,\ i < \partial,\ \ell < h(i) \big\}.
        $$
    \end{enumerate}
\sn
    \item[$\boxplus_{10}$] Let $\mathsf{X}$ denote the (finite) sequence of choices we have made so far: \emph{viz.} $\lambda,\mu,\partial$, $\bar\lambda$, $\big\LL S_\eps^* : \eps < \cf(\lambda) \big\RR$, . . . , $\big\LL\big\LL (\nu_\gamma,\rho_\gamma) : \gamma \in S_\xi^*\big\RR : \xi < \lambda \big\RR$.
\end{enumerate}
As in $(*)_4$+$(*)_5$ of Case 1,
\begin{enumerate}
    \item [$\boxplus_{11}$] 
    \begin{enumerate}
        \item We will show that $\bar I = \LL I_\zeta : \zeta < \lambda\RR$ exemplifies the conclusion; this will suffice. 
\sn
        \item So let $\zeta_* \in [\lambda_{<\eps},\lambda_\eps)$, $I \defeq I_{\zeta_*}$, 
        $J \defeq \sum\limits_{\zeta \neq \zeta_*}I_\zeta$, and let 
        $\chi,x,f_0,f_1,<_\chi^*$ be as in \ref{1.4}(1).
    \end{enumerate}
\sn
    \item[$\boxplus_{12}$] It suffices to find $\big\LL (M_i,N_i) : i < \partial\big\RR$ as in clauses $(i)$-$(vii)$ of $(*)^{I,J}_{\lambda,\mu,\partial}$ in Definition \ref{1.4}(1).
\end{enumerate}

\sn
Now we can choose 
\begin{enumerate}
    \item[$\boxplus_{12.1}$]  We can choose $\clB_\gamma^*$, $M_\gamma^*$ by induction on $\gamma \in S_{\zeta_*}^\bullet$ such that
    \begin{enumerate}
        \item $\clB_\gamma^*,M_\gamma^* \prec (\clH(\chi),\in,<_\chi^*)$
\sn
        \item $\|\clB_\gamma^*\| < \lambda$ and\footnote{
            We could choose $\|M_\gamma^*\| = \mu$, but there would be no gain in doing so.
        } $\|M_\gamma^*\| = \mu_*$.
\sn
        \item $[M_\gamma^*]^{\leq\partial} \subseteq M_\gamma^*$, $[M_\gamma^*]^{<\kappa} \subseteq M_\gamma^*$, and $\mu_* + 1 \subseteq M_\gamma^*$.
\sn
        \item $\LL\clB_j,M_j^* : j \leq \beta\RR \in \clB_\gamma^* \cap M_\gamma^*$ for all $\beta < \gamma$.
\sn
        \item If $\beta \in C_{\zeta_*,\gamma}^\bullet$ \underline{then} $\big\LL M_j^* : j \in C_\gamma \cap (\beta+1) \big\RR,\,\LL\clB_j^* : j \leq \beta\RR \in M_\gamma^*$.
\sn
        \item $x,f_0,f_1$, $I,J$, and $\zeta_*$ all belong to 
        $\clB_\gamma^*$ {and} to $M_\gamma^*$.  
    \end{enumerate}
\sn
    \item [$\boxplus_{12.2}$]  We will also add $\clB_\gamma^\bullet$ and $N_\gamma^\bullet$ such that if $\gamma = \alpha + \beta$ for some $(\alpha,\beta) \in S_\xi^+ \times \theta$, \underline{then}
    \begin{enumerate}
        \item $\clB_\gamma^\bullet,N_\gamma^\bullet \prec (\clH(\chi),\in,<_\chi^*)$
\sn
        \item $\|\clB_\gamma^\bullet\| < \theta$, $\clB_\gamma^\bullet \cap \theta \in \theta$, and $\|N_\gamma^\bullet\| = \mu_*$.
\sn
        \item $[N_\gamma^\bullet]^{\leq\partial} \subseteq N_\gamma^\bullet$, $[N_\gamma^\bullet]^{<\kappa} \subseteq N_\gamma^\bullet$, and $\mu_* + 1 \subseteq N_\gamma^\bullet$.
\sn
        \item $\big\LL (\clB_j,\clB_j^\bullet,M_j^*,N_j^\bullet) : j \leq \beta \big\RR \in \clB_\gamma^\bullet \cap N_\gamma^\bullet$ for all $\beta < \gamma$.
\sn
        \item As in $\boxplus_{12.1}$(e) above.
\sn
        \item As above.
    \end{enumerate}
\sn
    \item [$\boxplus_{12.3}$]  
    \begin{enumerate}
        \item $E_1 \defeq \big\{ \delta < \lambda_{\zeta_*} : \lambda_* \divides \delta \text{ and } \big(\bigcup\limits_{\gamma \in S_\delta^\bullet \cap \delta} \clB_\gamma \big) \cap \lambda_{\zeta_*} = \delta\}$ is a club of $\lambda_{\zeta_*}$.
\sn
        \item For $\delta \in E_1 \cap S_{\zeta_*}^\eps$, let $E_{2,\delta} \defeq \Big\{ \alpha < \theta : \big( \bigcup\limits_{\substack{\beta < \alpha\\\beta \text{ limit}}} \clB_{\delta+\beta}^\bullet \big) \cap \theta = \alpha \Big\}$.
    \end{enumerate}
\sn
    \item [$\boxplus_{12.4}$]  
    \begin{enumerate}
        \item Choose $\alpha_1 \in E_{2,\alpha_1}$ such that $C_{\alpha_2}^\bullet$ weakly guesses $E_{2,\alpha_1}$.
\sn
        \item Choose $\alpha_2 \in E_1 \cap S_{\zeta_*}^\eps$ such that $C_{\alpha_1}^\eps$ weakly guesses $E_1$.
\sn
        \item Now choose $\big\LL (M_i,N_i) : i < \partial\big\RR$ as follows.
        \begin{enumerate}
            \item For $\iota = 1,2$, let $\beta_{1,\iota} \in C_{\delta  
            }^\eps$ be such that $\otp(C_{\delta  
            }
            ^\eps) \cap \beta_{1,\iota} =  
            {2i} + \iota$.
\sn
            \item For $\iota = 1,2$, let $\beta_{2,\iota} \in C_{\alpha_2}^\bullet$ be such that $\otp(C_{\alpha_2}^\bullet) \cap \beta_{2,\iota} =   
        {i}$.
\sn
            \item Let $M_i \defeq N_{\beta_{1,1}+\beta_{2,1}}^\bullet$ and $N_i \defeq N_{\beta_{1,2}+\beta_{2,2}}^\bullet $.
        \end{enumerate}
    \end{enumerate}
\end{enumerate}
The rest is as in the proof of Case 1.
\end{PROOF}

\mn
We still owe the reader a proof of \ref{1.7}.

\begin{PROOF}{\ref{1.7}}
We would like to apply \ref{1.15}. Note that the case `$\max(\bfe_i) = 1$ for every $i < \partial$' simplifies matters, as $\mu_2 \to (h(i))_{\mu_2}^1$ holds.

\medskip
As we assume $\lambda > \mu^\partial + \mu^{<\kappa}$ (which implies $\lambda > \partial^+$) we know that Case 2 applies whenever $\lambda$ is regular, and use Case 4 when it is singular.
\end{PROOF}

\newpage

\section{Applications to Boolean Algebras}
\label{par2}

Here we construct some Boolean algebras with ``no non-trivial morphism.''

We shall mainly use $\BA_\tr(I)$, $I\in K^\omega_\tr$ for 
constructing mono-rigid ccc Boolean algebras; $\BA_{\tr(h)} (I)$, 
$I\in K^\omega_{\tr(h)}$, $h\in {}^\omega(\omega\setminus\{0,1,2\})$  for constructing complete
mono-rigid ccc Boolean algebras; and $\BA_\trr(I)$, $I\in K^\omega_\tr$ 
for constructing Bonnet--rigid Boolean algebras. In each case, for every $I$ from 
a relevant family (which exemplifies full bigness in the relevant case), we derive a 
Boolean algebra $\BA_\x (I)$, chosen to fit the proof of the case of rigidity we are
interested in (this is Definition \ref{2.1}). We then build a Boolean algebra $\bfB$ 
of cardinality $\lambda$, planting a copy of $\BA_\x(I_a)$ below enough elements 
$a\in\bfB$ such that $a\neq b \Rightarrow I_a\neq I_b$ (see \ref{2.4}). We mainly 
show that $\BA_{\tr(h)} (I)$ satisfies a strong version of the ccc hence the ccc 
is preserved (see \ref{2.6}), hence the outcome of the construction \ref{2.4} is as
required with respect to the ccc, completeness, and cardinality.  
We then observe the relevant weak representability results (see \ref{2.10}). 
Note that if we consider the completion of a ccc Boolean algebra $\bfB$ and 
$\bfB$ is weakly represented in 
$\cM_{\aleph_0,\aleph_0} (J)$ then its completion is weakly represented in 
$\cM_{\aleph_1,\aleph_1}(J)$. 
Next (in \ref{2.11}) we deal with deducing unembeddability of
$\BA_\x(I)$ into a Boolean algebra 
$\bfB$ which is weakly represented in 
$\cM_{\mu,\kappa} (J)$, the main case is part (2). 
We deduce as conclusions that there are mono-rigid [complete] Boolean algebras (\ref{2.13}, \ref{2.14}). We then deal with Bonnet-rigid Boolean algebras (\ref{2.15} `til the end).

\begin{definition}\label{2.1}
1) For $I\in K^\partial_\tr$, recall that $\BA_\tr(I)$ is the Boolean
algebra generated freely by $\{x_\eta : \eta \in I\}$, except that:
\begin{enumerate}
    \item[$(*)_1$] $\eta\lhd\nu\in P^I_\partial \Rightarrow\ x_\eta\ge x_\nu$.
\end{enumerate}
\mn
2)  For $I\in K^\partial_\ptr$ let $\BA_\ptr(I)$ be the Boolean algebra
freely generated by $\{x_\eta : \eta \in I\}$, except that for $\eta\in I$ with
$\lh(\eta) = \omega$, letting 
$\eta = \big\LL\LL\alpha_0,\beta_0\RR,\ldots,\LL\alpha_i,\beta_i\RR\ldots\big\RR_{i<\partial}$, 
the following holds:
\mn
\begin{enumerate}
    \item[$(*)_2$] For all $i < \partial$, $x_\eta \leq x_{\eta\rest i \caret \LL\alpha_i\RR}$ and $x_\eta \cap x_{\eta\rest i \caret \LL\beta_i\RR} = 0$.
\end{enumerate}

\mn
3) For $h\in {}^\partial(\omega\setminus\{0\})$ and $I\in
K^\partial_{\tr(h)}$, let $\BA_{\tr(h)}(I)$ be the Boolean algebra generated freely by
$\{x_\eta : \eta \in I\}$, except that for $\eta \in P^I_\partial$ and
$i < \partial$, letting $\eta(i) = \LL s_0,\ldots,s_{h(i)-1})\RR$, we have:
\begin{enumerate}
    \item[$(*)_3$] $x_\eta \leq x_{\eta\rest i \caret \LL s_0\RR}$
    and $x_\eta \cap \bigcap\limits^{h(i)-1}_{\ell=1} x_{\eta\rest i \caret \LL s_\ell \RR} = 0$.
\end{enumerate}
\mn
The second equality is trivial if $h(i)=1$, so 
usually $h\in {}^\partial(\omega\setminus \{0,1\})$. If $(\forall i) [h(i)=1]$ 
this is like the case of $I \in K^\partial_\tr$, and if $(\forall i) [h(i)=2]$
this is like the case of $I\in K^\partial_\ptr$.


\sn
4)  For $I\in K^\partial_\tr$ (or just $I$ is a set of sequences of ordinals
closed under initial segments) let $\BA_\trr(I)$ be the Boolean
algebra freely generated by $\{x_\eta : \eta \in I\}$, except that:
\begin{enumerate}
    \item  $x_{\eta \caret \LL\alpha\RR}\cap x_{\eta \caret \LL\beta\RR}= 0$ for\footnote{
        We are, of course, assuming $\eta \caret \LL \alpha \RR,\eta \caret \LL \beta \RR \in I$; similarly in other cases.
    } 
    $\alpha\ne\beta$.
\sn
    \item  $x_\eta\le x_\nu$ for $\nu\lhd\eta$.
\sn
    \item If $\eta$ has finitely many immediate successors 
    $\{\eta \caret \LL\alpha_\ell\RR : \ell < k_\eta\}$ and $k_\eta \ge 2$ \underline{then}  $x_\eta = \bigcup\{x_{\eta \caret \LL\alpha_\ell\RR} : \ell < k_\eta\}$.
\sn
    \item  If $\eta\lhd\nu$ and every $\rho$ satisfying $\eta \unlhd \rho \lhd \nu$ has a unique successor, then $x_\eta = x_\nu$.
\end{enumerate}

\mn
5)   For $I\in K^\partial_{\tr(h)}$ and $g\in {}^\partial\omega$,
$h\in {}^\partial(\omega\setminus\{0,1\})$ 
satisfying\footnote{
    I.e. $(\forall i)[g(i)\leq h(i)]$.
} 
$g\leq h$, we define $\BA_{\tr(h,g)}(I)$ as the Boolean algebra generated freely by $\LL x_\eta : \eta\in I\RR$, except that:
\mn
\begin{enumerate}
    \item[$(*)_5$] If $\eta \in I$, $\lh(\eta) = \partial$, $j < \omega$, and
    $\eta(j) = \LL\alpha_0,\ldots,\alpha_{k-1}\RR$ where $k = h(j)$,\underline{then}
    \begin{enumerate}
        \item[$(\alpha)$]  $x_\eta \leq \bigcup\limits_{m=0}^{g(j)-1} 
        x_{(\eta \rest j) \caret \LL\alpha_m\RR}$
\sn
        \item[$(\beta)$]  
        $x_\eta \cap 
        \bigcap\limits^{h(j)-1}_{m=g(j)} x_{(\eta\rest j) \caret \LL \alpha_m\RR} = 0$.     (Usually we assume $0<g<h$.)     
    \end{enumerate}
\end{enumerate}
\mn
6)  Assume that $h\in {}^\partial(\omega\setminus \{0,1\})$ and $\bar\bfe$ is a sequence of length $\partial$ with $\bfe_i = \{(u_{1,i},u_{2,i})\}$, where 
$u_{1,i}, u_{2,i}$ are disjoint non-empty subsets of $h(i)$  and $|u_{1,i}| \geq 2$.  
For $I \in K^\partial_{\tr(h)}$, we define $\BA_{\tr(h),\bar\bfe}(I)$ 
as the Boolean algebra freely generated by 
$\{x_\eta : \eta \in I\}$, except that for $\eta \in P^I_\partial$ and 
$i < \partial$, letting $\eta(i)=\LL s_0,\ldots, s_{h(i)-1}\RR$, we have 
\begin{enumerate}
    \item[$(*)'_3$]   $x_\eta \leq \bigcap\limits_{\ell\in u_{1,i}} x_{(\eta\rest i) \caret \LL s_\ell\RR}$ and $x_\eta \cap \bigcap\limits_{\ell\in u_{2,i}} x_{(\eta\rest i) \caret \LL s_\ell\RR}=0$. 
\end{enumerate}
\mn
(We have much freedom in this case. Note that $\bar\bfe$ plays a different role here than it did in \S1: compare with Definition \ref{1.4}.)
\end{definition}

\mn
\begin{notation}\label{2.2}
1)  Clearly $K^\partial_{\tr(h,g)} = K^\partial_{\tr(h)}$ for $g$ as in \ref{2.1}(3) and $h \in {}^\partial\{1\}$. 
Note that for $I \in K^\partial_{\tr(h)}$, if $g=1$ then $\BA_{\tr(h,g)}(I)$ is 
essentially $\BA_{\tr(h)}(I)$. Also, if $h=1$ then 
$K^\partial_{\tr(h)} = K^\partial_\tr$ and $\BA_{\tr(h)}(I) = \BA_\tr(I)$.

\sn
2) When we state a result that holds for $\tr,\, \ptr,\, \trr,\, \tr(h)$, or $\tr(h,g)$, we will replace the corresponding subscripts with an $x$. 

\sn
3)  Recall that when we say ``a Boolean algebra is freely generated by 
$$
X \defeq \{x_i : i \in U\},
$$
except the set equations . . ,'' we have
\textbf{0} and \textbf{1} ($= -\mathbf{0}$) in the Boolean algebra.

\sn
4) For a Boolean algebra $\bfB$ and $a \in \bfB$, $\bfB \rest a$ is the
naturally defined Boolean algebra, but $\mathbf{1}_{\bfB\rest a} = a$. Essentially, we do not consider $\mathbf{1}_\bfB$ as an {individual} constant of $\bfB$.
\end{notation}

\mn
\begin{definition}\label{2.3}
For Boolean algebras $\bfB$, $\bfB_1$ and $a^* \in \bfB_1\setminus\{\mathbf{0}_{\bfB_1}\}$, 
we define the ``$\bfB$-surgery of $\bfB_1$ at $a^*$'' or ``surgery of $\bfB_1$ 
at $a^*$ by $\bfB$'', called $\bfB_2$, as a Boolean algebra extending $\bfB_1$ 
such that $\bfB_2 = [\bfB_1 \rest (-a^*)]\times [(\bfB_1 \rest a^*) * \bfB]$, 
where $\times$ is a direct product and $*$ free product. Alternatively, 
$\bfB_2$ can be generated as follows: first make $\bfB$ disjoint to $\bfB_1$ 
(by taking an isomorphic copy) and then $\bfB_2$ is freely generated by 
$\bfB_1 \cup \bfB$, except the relations
\begin{align*}
& \mathbf{0}_{\bfB_1}=\mathbf{0}_\bfB =0,\\
&a\cap b=c\quad\text{ (for $a,b,c\in \bfB_1$ such that $a \cap b = c$ in $\bfB_1$)},\\
&a\cup b=c\quad\text{ (for $a,b,c\in\bfB_1$ such that $a\cup b=c$ in $\bfB_1$)}
\end{align*}

\begin{align*}
\mathbf{1}_{\bfB_1}-b=c\quad &\text{ (for $b,c\in \bfB_1$ such that  $\mathbf{1}_{\bfB_1}-b=c$ in
$\bfB_1$)},\\ 
a\cap b=c\qquad &\text{ (for $a,b,c\in \bfB$ such that 
$a\cap b=c$ in $\bfB$)},\\
a \cup b = c\qquad &\text{ (for $a,b,c\in \bfB$ such that $a\cup b=c$ in $\bfB$)},\\
\mathbf{1}_\bfB - b = c\quad &\text{ (for $b,c\in\bfB$ such that $\mathbf{1}_{\bfB}-b=c$)}
\end{align*}

and
\[
\mathbf{1}_{\bfB}=a^*.
\]
\end{definition}

\begin{construction}\label{2.4}
Let $x$ be one of $\{\tr,\ptr,\tr(h),\trr,\tr (h,g)\}$ and let $\lambda$ 
be a cardinal such that $\alpha < \lambda^+$ (usually $\alpha = \lambda$, always
$\alpha > 0$). The idea is to construct a Boolean algebra by defining an
increasing continuous sequence $\bfB_i$ ($i\leq\alpha$), $\bfB_0$ trivial, 
and we get $\bfB_{i+1}$ by a surgery of $\bfB_i$ at $a^*_i\in \bfB_i$ by
$\bfB^*_i= \BA_\x(I_i)$ (see Definition \emph{\ref{2.1}} and \emph{\ref{2.2}(2)}), where 
$|I_i| = \lambda$, $I_i \in K^\partial_x$, and $I_i$ is 
strongly $\psi_x$-unembeddable into $\sum\limits_{j<\alpha, j\neq i}I_j$ 
(or e.g.~super$^y$-$\bar\bfe$-unembeddable into it, for $y \in \{\nr,\vr\}$; see \emph{\ref{1.5}}). 

We denote $\bfB = \bfB_\alpha$ by $\Sur_x\LL I_i, a^*_i : i < \alpha\RR$.
Usually we would like to have
$\bfB_\alpha \setminus \{0\} = \{a^*_i : i < \beta_\alpha\}$. If there are 
$\LL I_i : i<\alpha\RR$ as above and $\alpha$ is divisible by $\lambda$ 
then this is clearly possible.
\end{construction}

\mn
\begin{definition}\label{2.5}
1)  A Boolean algebra satisfies the $\lambda$-chain condition
(or the $\lambda$-\emph{cc}) \underline{iff}  there are 
no $\lambda$ elements which form an antichain (i.e.~they are $\neq 0$ and the intersection of any two is zero).

\sn
2)  A Boolean algebra satisfies the strong $\lambda$-chain condition
or the $\lambda$-Knaster condition \underline{iff} 
among any $\lambda$ elements there are $\lambda$ which are pairwise not
disjoint.
\end{definition}

\mn
\begin{claim}\label{2.6}
Let $\x\in\{\tr,\ptr,\tr(n),\tr(h),\tr(*)\}$, $I\in K^\partial_x$, 
$\lambda > \partial$ regular.

\sn
$1)$  If $\x=\tr$ \underline{then}  $\BA_\x(I)$ satisfies the strong
$\lambda$-chain condition.

\sn
$2)$  If $\x=\ptr$ \underline{then}  $\BA_\x(I)$ satisfies the strong
$\left(2^\partial\right)^{\!+}\!\!$-chain condition.

\sn
$3)$  If $k\geq 3$ and $I \in K^\partial_{\tr(k)}$ is
standard, \underline{then}  $\bfB = \BA_{\tr(k)}(I)$ satisfies the strong
$\lambda$-chain condition; similarly for $K^\partial_{\tr(*)}$,
for $K^\partial_{\tr(h)}$ with $h\in {}^\partial(\omega\setminus 3)$, and
$K^\partial_{\tr(h,g)}$ (for $h\in {}^\partial(\omega\setminus 3)$ and 
$g\in {}^\partial \omega$ such that $g\leq h$). 

Instead of $h\in {}^\partial(\omega\setminus 3)$, we can demand 
$h\in {}^\partial (\omega\setminus 1)$ and $h(i)\geq 3$ for every large enough $i$.

\sn
$4)$ If $\x=\ptr$ then $\BA_\x(I)$ satisfies the strong $\lambda$-chain condition
provided that $I$ is atomically ($<\lambda$)-stable; for example, if
$(\forall \alpha<\lambda)\big[|\alpha|^\partial < \lambda \big]$.

\sn
$5)$ If $h,\bar\bfe$ are as in \emph{\ref{2.1}(6)} and 
$(*)^i_{\bar\bfe}$ below holds for every $i$ large enough, $\lambda > 
{\partial}$ is regular 
     uncountable 
, and 
$I \in K^\partial_{\tr(h)}$, \underline{then} 
$\BA_{\tr(h),\bar\bfe}(I)$ 
satisfies the strong $\lambda$-chain condition, where:
\begin{enumerate}
    \item[$(*)^i_{\bar\bfe}$]  $\bfe_i=\{(u^i_1,u^i_2)\}$, where 
    $u^i_1,u^i_2 \subseteq \{0,\ldots h(i)-1\}$ are disjoint and non-empty with $|u_2^i| \geq 2$.
\end{enumerate}
\end{claim}

\sn
\begin{remark} 
Clearly we can similarly phrase sufficient conditions for ``for any family of
$\lambda$ non-zero elements there is an uncountable subfamily such that any
$k$ members of the subfamily have non-zero intersection."
\end{remark}

\noindent
Before we prove \ref{2.6}, recall the well known fact: 
(Here $\bfB_0=\{0, 1\}$ is the two-element Boolean
algebra.)

\begin{fact}\label{2.6A}
1) If $\bfB$ is the Boolean algebra freely generated by $\{x_t : t\in I\}$ 
except {for a} set $\Lambda$ of equations in $\{x_t : t\in I\}$, 
(so each member of $\Lambda$ has the form 
$\sigma(x_{t_0},\ldots,x_{t_{n-1}})=0$, where $\sigma(y_0,\ldots,y_{n-1})$ 
is a Boolean term, $t_0,\ldots,t_{n-1}\in I$) \underline{then}, for a Boolean 
term $\sigma^*(x_{s_0},\ldots,x_{s_{n-1}})$, we have 
$(\alpha) \Leftrightarrow (\beta)$, where:
\mn
\begin{enumerate}
    \item[$(\alpha)$]  $\bfB\models\sigma^*(x_{s_0},\ldots,x_{s_{n-1}}) > 0$
\sn
    \item[$(\beta)$]  For some function $f : I \to \{0,1\}$, we have:
    \begin{enumerate}
        \item  $f$ respects $\Lambda$; i.e.~recalling $\bfB_0 \defeq \{0,1\}$, 
        $$\sigma (x_{t_0},\ldots,x_{t_{m-1}}) \in \Lambda \ \Rightarrow\ 
        \bfB_0\models ``0 = \sigma(f(t_0),\ldots,f(t_{m-1}))".$$

        \item $\bfB_0\models`\sigma^*(f(s_0),\ldots,f(s_{n-1}))=1$'
    \end{enumerate}
\end{enumerate}
\mn
2)  In fact, if $f :  I\to \{0,1\}$ satisfies clause (a) \underline{then} 
there is a unique homomorphism $\hat{f}$ from $\bfB$ into
 $\bfB_0$ such that $s \in I \Rightarrow \hat{f}(x_s) = f(s)$.
\end{fact}

Now we return to proving \ref{2.6}.

\begin{PROOF}{\ref{2.6}}
1)  We take $x=\tr$ and check the strong $\lambda$-chain conditions. 
Note that by \ref{2.6A} and the definition of $\BA_\tr(I)$, we have: 
\mn
\begin{enumerate}
    \item[$(*)_1$]  $x_{\eta_1} \cap \ldots \cap x_{\eta_k} \cap 
    (-x_{\nu_1}) \cap \ldots \cap (-x_{\nu _m}) = 0$\ iff 
    $$(\exists i,j)\big[\nu_i \lhd \eta_j \in P^I_\omega \vee \nu_i = \eta_j \big].$$
\end{enumerate}
\mn
[Why?  The `if' implication is trivial, recalling Definition \ref{2.1}(1). 
For proving the `only if' implication, assume that the second statement 
holds. Define $f : I \to \{0,1\}$ by $f(\eta)=1$ iff 
$(\exists \ell) [\eta=\eta_\ell \vee \eta\lhd \eta_\ell \in P^I_\partial]$; 
clearly it respects the equations in the definition of $\BA_\tr(I)$
and $\hat{f}$ maps $x_{\eta_1}\cap\ldots\cap x_{\eta_k} \cap 
(-x_{\nu_1}) \cap\ldots\cap (-x_{\nu_n})$ to 1, so by \ref{2.6A} we are done.]

Now for $u\in [I]^{<\omega}$, let $x_u=\bigcap\limits_{\eta\in u} x_\eta$ and
$x_{-u}=\bigcap\limits_{\eta\in u} (-x_\eta)$. Clearly, if $a\in
\BA_\x(I) \setminus\{0\}$ then for some $u,v\in [I]^{<\aleph_0}$,
we have $0 < x_u \cap x_{-v} \le a$ (hence $u$ and $v$ are disjoint). 
In fact, $a$ is a finite union of such elements. To check the strong 
$\lambda$-chain condition it suffices to take 
$\{(u_i,v_i) : i < \lambda\}\subseteq [I]^{<\aleph_0}\times [I]^{<\aleph_0}$ 
such that $(\forall i<\lambda)[x_{u_i}\cap x_{-v_i}\neq 0]$, and
to find $A\in [\lambda]^\lambda$ such that
\[
(\forall i,j\in A) \big[ x_{u_i} \cap x_{-v_i} \cap x_{u_j} \cap x_{-v_j} \neq 0 \big].
\]

\mn
We may assume that $\LL u_i : i\in A \RR $ and $\LL v_i : i\in A\RR$ are 
$\Delta$-systems (say with hearts $u^*$, $v^*$, respectively), so  
$u_i\cap v_i= \varnothing$ and 
$$
u_i\cap v^* = u^*\cap v_i = u^*\cap v^* = \varnothing.
$$
We may assume $i \ne j \in A$ implies $u_i\cap v_j = \varnothing$, 
$u_i \neq u_j$, and $v_i \neq v_j$. We may assume that for some non-zero 
$m,n < \omega$, for every $i \in A$, we have $|u_i| = m \wedge |v_i| = n$.
Say $u_i=\{\eta_{i,\ell} : \ell < m\}$, $v_i = \{\nu_{i,\ell} : \ell < n\}$ 
(without repetitions) and for each $\ell < m$ the sequence 
$\LL \eta_{i,\ell}:i \in A\RR$ is constant or is without
repetitions, and similarly $\LL \nu_{i,\ell}:i \in A\RR$. 
We may also assume
\begin{enumerate}
    \item[$(*)_2$]  $\LL\lh(\eta_{i,\ell}) : \ell < m\RR$, 
    $\LL\lh(\nu_{i,\ell}) : \ell < n\RR$ is the same for all $i\in A$.
\end{enumerate}

\mn
Clearly then, using the $\Delta$-system assumption,
\begin{enumerate}
    \item[$(*)_3$] For $i\in A$, $\ell<m$, $k<n$ there is at most one $j\in A$ 
    such that $\nu_{j,k} \lhd \eta_{i,\ell} \in P^I_\partial$.
\end{enumerate}

\mn
[Why?  If we have $\nu_{j,k} \unlhd \eta_{i,\ell} \in P^I_\partial$, 
note that $\nu_{i,k} \not\!\!\unlhd\, \eta_{i,\ell}$ 
by $(*)_1$, hence $\nu_{j,k} \neq \nu_{i,k}$. So $i\neq j$ and hence $\nu_{j,k}\notin v^*$, and 
$\nu_{j,k} = \eta_{i,\ell} \rest \lh(\nu_{j, k})$. Thus 
$j \neq j_1 \in A \Rightarrow \nu_{j_1,k} \neq \nu_{j,k}$ and hence 
$$
j \neq j_1 \in A \Rightarrow \nu_{j_1,k} \neq \eta_{i,\ell} \rest \lh(\nu_{j, k}) = \eta_{i,\ell} \rest \lh(\nu_{j_1,k}).
$$
Hence 
$j\neq j_1 \in A \Rightarrow \neg[\nu_{j_1,k} \unlhd \eta_{i,\ell}]$ and we have finished.]

\medskip
So for $i\in A$, the set
\[
w_i \defeq \{j:\text{for some }\ell<m,\ k<n\text{ we have }\nu_{j,k}
\lhd \eta_{i,\ell}\in P^I_\omega\}
\]
has at most $mn<\aleph_0$ members. So by $(*)_1$ it suffices to find $A'\in
[A]^\lambda$ such that $i \neq j \in A' \Rightarrow j\notin w_i$. By
the Hajnal Free Subset theorem \cite{Ha61}\footnote{
    Or see \cite[3.14\subref{4.Ha}]{Sh:E62}.
}
there is\footnote{
    Note that $(-x_{\nu_{j_1,\ell_1}}) \cap (-x_{\nu_{j_2,\ell_2}}) > 0$ always holds.
} 
such an $A'$.

\sn
2)  The case $x=\ptr$ is similar, but more
complicated. First note\footnote{
    This will also be used in the proof of part (4).
}
\mn
\begin{enumerate}
    \item[$(*)_4$] Assume $I\in K^\partial_{\ptr}$ and $\bfB=\BA_{\ptr} (I)$. If $m,n<\omega$ and $\nu_k,\eta_\ell\in I$ for $\ell<m$ ad $k<n$, \underline{then} 
    $$
    \bfB\models x_{\eta_0}\cap\ldots\cap x_{\eta_{m-1}}\cap (-x_{\nu_0}) \cap \ldots \cap (-x_{\nu_{n-1}}) = 0
    $$ 
    \underline{iff} at least one of the following conditions holds:
    \begin{enumerate}
        \item   $(\exists\ell,k < m) \big[ \lh(\eta_\ell) = \partial \wedge \Suc_R(\eta_k,\eta_\ell) \big]$
\sn
        \item  $(\exists\ell<m)(\exists k<n) \big[ \lh(\eta_\ell) = \partial \wedge \Suc_L(\nu_k,\eta_\ell) \big]$
\sn
        \item  $(\exists\ell,k<m)(\exists j<\omega)(\exists \alpha,\beta,\gamma)\big[\lh(\eta_\ell) = \lh(\eta_k) = \partial \wedge\\ \eta_\ell\rest j = \eta_k\rest j\ \wedge\ \eta_\ell(j) = \LL\alpha,\beta\RR\ \wedge\ \eta_k(j) = \LL\beta,\gamma\RR\big]$
\sn
        \item  $(\exists\ell < m)(\exists k < n)[\eta_\ell = \nu_k]$.
    \end{enumerate}
\end{enumerate}
\mn
[Why? If (a) or (b) or (c) or (d) holds then the intersection is zero by
the equations we have imposed defining $\BA_{\ptr}(I)$ in Definition 
\ref{2.1}(2), so the ``if" implication holds. Next we prove the other 
implication, so we assume (a), (b), (c), and (d) fail, and we shall use \ref{2.6A}. 
We have to define $f(\rho)$ for $\rho\in I$; we do it  by cases.

\mn
\textbf{Case 1}: $\lh(\rho)=\partial$, $\rho\in\{\eta_0,\ldots,\eta_{m-1}\}$. 

Let $f(\rho)=1$.

\mn
\textbf{Case 2}: $\lh(\rho)=\partial$ but Case 1 does not hold. 

Let $f(\rho)=0$.

\mn
\textbf{Case 3}:  $\lh(\rho)< \partial$ and for some $\ell<m$ we have $\lh (\eta_\ell)
=\partial$ and $\Suc_L(\rho,\eta_\ell)$. 

Let $f(\rho)=1$.

\mn
\textbf{Case 4}: $\lh(\rho) < \partial$ and for some $\ell<m$ we have 
$\lh(\eta_\ell) = \partial$ and $\Suc_R(\rho,\eta_\ell)$. 

Let $f(\rho)=0$.

\mn
\textbf{Case 5}: $\lh(\rho) < \partial$, $\rho\in\{\eta_\ell : \ell < m\}$. 

Let $f(\rho)=1$.

\mn
\textbf{Case 6}: No previous case applies. 

Let $f(\rho)=0$.

\mn
First, $f$ is well-defined. I.e.~there are 
no contradictions in this definition --- the less trivial cases are 
between cases 3+4, 3+5, and 4+5. 
(We are safe, as clauses (c), (b), and then (a) of $(*)_4$ would fail, 
respectively.\footnote{
    Actually, cases 3+5 cannot contradict each other.
})

Second, we show that $f$ respects 
the equations from Definition \ref{2.1}(2); that is, from $(*)_2$ there. 
If $x_\eta \leq x_{\eta\rest i \caret \LL \alpha_i\RR}$ is an instance
of $(*)_2$ of \ref{2.1}(2) and $f$ fails it (that is, $f(\eta)=1$,
$f(\eta\rest i \caret \LL\alpha_i\RR)=0$) then necessarily by 
$\lh(\eta) = \partial$ Case 1 occurs for $\eta$, hence Case 3 occurs for 
$(\eta\rest i) \caret \LL\alpha_i\RR$.   
But by Case 3 $f((\eta\rest i) \caret \LL\alpha_i\RR)=1$: a contradiction.

Similarly for the other equation in
$(*)_2$ of \ref{2.1}(2), using Case 4 instead Case 3. Third:
$f(x_{\eta_\ell})=1$ for $\ell<m$ by Cases 1, 5, and $f(\nu_k) = 0$ 
for $k<n$ as by failure of clause (d), Case 2 occurs if 
$\lh(\nu_k)=\partial$, and Case 6 occurs if 
$\lh(\nu_k) < \partial$. So by \ref{2.6A} we are done proving $(*)_4$.]

\smallskip
Let $a_\alpha\in \BA_\x(I)\setminus\{0\}$ for 
$\alpha < \lambda = (2^\partial)^+$, 
so as before without loss of generality 
$$
a_\alpha = x_{\eta_{\alpha,0}} \cap \ldots \cap x_{\eta_{\alpha,n_\alpha-1}} \cap (-x_{\eta_{\alpha,n_\alpha}}) \cap \ldots \cap (-x_{\eta_{\alpha,m_\alpha-1}}) \neq 0.
$$ 
Without loss of generality $n_\alpha=n^*$, $m_\alpha=m^*$ and $P^I_\partial \cap \{\eta_{\alpha,\ell} : \ell < m^*\} \neq \varnothing$ 
(for notational simplicity below). We can define
$\eta_{\alpha,\ell}$ (for $\ell \in [m^*, \partial)$) such that 
\[
\Suc_L(\rho,\eta_{\alpha,\ell})\vee \Suc_R (\rho,\eta_{\alpha,\ell})
\Rightarrow \rho \in \{ \eta_{\alpha,j} : j < \partial\}.
\]

\smallskip
Without loss of generality, the atomic type of 
$\LL\eta_{\alpha,\ell} : \ell < \partial\RR$ in $I$ does
not depend on $\alpha$ and they form a $\Delta$-system: i.e.
\sn
\begin{enumerate}
    \item[$(*)$]  $\eta_{\alpha,\ell_1} = \eta_{\beta,\ell_2} \wedge \alpha \neq
    \beta\ \Rightarrow\ (\forall\alpha_1,\beta_1<\lambda)[\eta_{\alpha_1,\ell_1} = \eta_{\alpha_1,\ell_2}  =\eta_{\beta_1,\ell_1} = \eta_{\beta_1,\ell_2}]$.
\end{enumerate}
\mn
Now we apply $(*)_4$: check that each case fails.

\bn
5) As part (3) is a special case of part (5), we will prove the latter.

Without loss of generality we deal with
$K^\partial_{\tr(h,g)}$. Let $a_\alpha\neq 0$ (for $\alpha<\lambda$) 
be non-zero pairwise disjoint elements, let 
$a_\alpha = \sigma_\alpha(\bar{x}_{\bar{\eta}_\alpha})$, where 
$\sigma_\alpha$ a Boolean term and $\bar{\eta}_\alpha$ a finite 
sequence from $I$ (i.e.~we write $\bar{x}_{\LL\eta_{\alpha,0},\ldots,\eta_{\alpha,k_\alpha-1}\RR}$
instead of $\LL x_{\eta_{\alpha,0}},\ldots,x_{\eta_{\alpha,k_\alpha-1}}\RR$). 
Without loss of generality  $\sigma_\alpha = \sigma$ and
$\bar{\eta}_\alpha=\LL\eta_{\alpha,0},\ldots,\eta_{\alpha,k-1}\RR$
is without repetition, and
\[
a_\alpha=\bigcap\limits_{\ell<k(0)}x_{\eta_{\alpha,\ell}}\cap
\bigcap\limits_{k(0)\leq\ell<k}\left(-x_{\eta_{\alpha,\ell}}\right).
\]
As earlier, 
\begin{enumerate}
    \item[$(*)_0$] Without loss of generality $\alpha \neq \beta\, \wedge\, \eta_{\alpha,\ell} = \eta_{\beta,k} \Rightarrow \ell = k$.
\end{enumerate}

So there is $i_\alpha = i(\alpha) <   
   {\partial}$ such that 
\begin{enumerate}
    \item[$(*)_1$] $\lh(\eta_{\alpha,\ell}) < \partial\, \Rightarrow\, \lh(\eta_{\alpha,\ell}) \leq i_\alpha$, and $\lh(\eta_{\alpha, \ell_1}) = \lh(\eta_{\alpha,\ell_2}) = \partial \wedge \ell_2 \neq \ell_2$ implies
\[
    \eta_{\alpha,\ell_1} \rest i_\alpha \neq \eta_{\alpha,\ell_2} \rest i_\alpha 
    \text{ and } (\forall i \geq i_\alpha-1) \big[ h(i) \geq 3 \big].
\]
    \item[$(*)_2$] Without loss of generality, 
    \begin{enumerate}
        \item $\   
           {u_\alpha} = u_*$, where
        \begin{align*}
            v_\alpha \defeq \big\{i < 
              i_ \alpha 
              : &\ (\exists \ell < k) \big[ \lh(\eta_{\alpha,\ell}) = i \big] \text{ \underline{or} }\\
            &\ (\exists 
              {\ell_1,\ell_2}
            < k) \big[ \max\{j < \partial : \eta_{\alpha,\ell_1} \rest j = \eta_{\alpha,\ell_2} \rest j\} = i \big] \big\}
        \end{align*}
\sn
        \item If $i \in v_*$, $\ell < k$ and $\lh(\eta_{\alpha,\ell}) > i$, \underline{then} $\Res_i^m(\eta_{\alpha,\ell}) \in \{\eta_{\alpha,m} : m < k\}$.
\sn
        \item The truth values of `$\lh(\eta_{\alpha,\ell}) = i$' and `$\Res_i^m(\eta_{\alpha,\ell}) = \eta_{\alpha,k}$' do not depend on $\alpha$.
    \end{enumerate} 
\end{enumerate}
[Why? Easy. (If necessary, we can change $\bar \eta_\alpha$ and $\sigma_\alpha$, and then uniformize $\sigma_\alpha$ and $k$ again.)]

\begin{enumerate}
    \item[$(*)_3$] Without loss of generality, for every $i \in v_*$ 
       {and $\ell$},
    \begin{enumerate}
        \item $(\forall \ell_1,\ell_2 < k)\big[\alpha \neq \beta \wedge \eta_{\alpha,\ell_1} = \eta_{\beta,\ell_2} \Rightarrow \ell_1 = \ell_2 \big]$
\sn
        \item For some $W_* \subseteq k$,
        \begin{itemize}
            \item If $\ell \in W_*$ then 
            $\LL\eta_{\alpha,\ell} : \alpha < \lambda\RR$ is 
            constant (
              {call this value} $\eta_\ell^*$).
\sn
            \item If $\ell \in k \setminus W_*$ then $\LL\eta_{\alpha,\ell} : \alpha < \lambda\RR$ is without repetition.
        \end{itemize}
    \end{enumerate}
\sn
    \item[$(*)_4$] Assume $\alpha \neq \beta < \lambda$, and we will prove that $B \models `a_\alpha \cap a_\beta\ 
        {>\, k}$'.
\end{enumerate}

\mn
Let $\bfe_i \defeq (u_1^i,u_2^i)$. We will define the function $f : I \to \{0,1\}$ as follows.
\begin{enumerate}
    \item [$(*)_{4.1}$] For $\rho \in I$, 
    $$
    f(\rho) \defeq 
    \begin{cases}
        1 &\text{if } \rho 
          \in \{ \eta_{\alpha,\ell}, \eta_{\beta,\ell} : \ell < k\} \\
        1 &\text{if } \lh(\eta_{\alpha,\ell}) = \partial,\ i < \partial,\ m \in u_1^i, \text{ and } \rho = \Res_i^m(\eta_{\alpha,\ell})\\
        0 &\text{otherwise.} 
    \end{cases}
    $$
\end{enumerate}
Clearly this is well-defined. Now we just need to show that $f$ respects the relevant equation. So assume, towards contradiction:
\begin{enumerate}
    \item [$(*)_{4.2}$] $  \rho 
       \in P_\partial^I$, $i < \gamma$, $f(\eta) = 1$, $\ell < h(i)$, and
    $$
    f(\rho) = 0,\ \rho( {\ell} 
        ) \in u_2^i \text{ 
          \underline{or} }  f(\rho) = 1 \wedge \rho( {\ell} 
           ) \in u_1^i
    $$
    \begin{enumerate}
        \item[$\bullet_1$] $f(\Res_i^\ell(\eta)) = 0$ for some $\ell \in u_1^i$,
    \end{enumerate}
    \underline{or}
    \begin{enumerate}
        \item[$\bullet_2$] $f(\Res_i^\ell(\eta)) = 1$ for every $\ell \in u_2^i$.
    \end{enumerate}  
\end{enumerate}
As $f(\eta) = 1$ and $\eta \in P_\partial^I$, necessarily $\eta \in \{ \eta_{\alpha,\ell}, \eta_{\beta,\ell} : \ell < k\}$, so $\bullet_1$ is impossible: therefore we assume $\bullet_2$.
Also by symmetry, without loss of generality we may assume $\eta \defeq \eta_{\alpha,\ell_*}$.

Now if $i \in v_*$, then this is impossible by $(*)_2$, hence $i \notin v_*$. So $f(\Res_i^\ell(\eta)) = 1$ for $\ell \in u_2^i$, hence there is $\nu_\ell \in \{ \eta_{\alpha,\ell}, \eta_{\beta,\ell} : \ell < k\} \cap P_\partial^I$ and $k_i \in u_1^i$ such that $\Res_i^{k_\ell}(\nu_\ell) = \Res_i^\ell(\eta)$. Also, we know $\nu_\ell \notin \big\{ \eta_{\alpha,m} : m  \in k \setminus \{\ell_*\} \big\}$ because $i \notin v_*$. 

All together,
$$
\nu_\ell \in \big\{ \eta_{\beta,m} : m  \in k \setminus W_*\big\}.
$$
But as $|u_2^i| \geq 2$ we can choose $\ell_1 \neq \ell_2 \in u_2^i$, so $\nu_1,\nu_2$ witness `$i \in \nu_*$!' This is a contradiction.

\vspace{7mm}

\noindent
4) Let $a_\alpha\in \BA_\x[I]\setminus\{0\}$ for 
$\alpha < \lambda$, 
so as before without loss of generality 
$$
a_\alpha = x_{\eta_{\alpha,0}} \cap \ldots \cap x_{\eta_{\alpha,n_\alpha-1}} \cap (-x_{\eta_{\alpha,n_\alpha}}) \cap \ldots \cap (-x_{\eta_{\alpha,m_\alpha-1}}) \neq 0.
$$
Without loss of generality $n_\alpha=n^*$, $m_\alpha=m^*$ and $P^I_\partial \cap \{\eta_{\alpha,\ell} : \ell < m^*\} \neq \varnothing$ 
(for notational simplicity below).

Let 
$$
\Lambda \defeq \{\eta \rest i : \eta \in P_\partial^I,\ i < \partial\} \cup (I \setminus P_\partial^I).
$$
Let $\LL \Lambda_\alpha : \alpha < \lambda\RR$ be such
that $\Lambda_\alpha \subseteq \Lambda$ with cardinality $<\lambda$, increasing 
continuously with $\alpha$ with union $\Lambda$. 

For each $\alpha < \lambda$, let 
$$
\Lambda_\alpha^+ \defeq \big\{\eta \in P_\partial^I : (\forall i < \partial) [\eta \rest i \in \Lambda_\alpha] \big\} \cup \big\{ \eta \caret \LL s\RR, \eta \caret \LL t\RR : \eta \caret \LL (s,t)\RR \in \Lambda_\alpha \big\}.
$$ 
So by our present assumptions (i.e.~$I$ is atomically stable and $|\Lambda_\alpha| < \lambda$) we can choose $\beta_\alpha < \lambda$ such that $\Lambda_\alpha^+ \subseteq \Lambda_{\beta_\alpha}^+$, and clearly 
$$
E_0 \defeq \big\{ \delta < \lambda : \delta \text{ a limit ordinal and } (\forall \alpha < \delta)[ \beta_\alpha < \delta] \big\}
$$
is a club of $\lambda$. 

Let $\LL\alpha_\eps : \eps < \lambda\RR$ list $E_0$ in increasing order. 

For $\eps < \lambda$ let 
\begin{align*}
    \base(\eps) \defeq \big\{\zeta < \eps : &\text{ for some } \ell < m_*,\ i < \partial,\ m < h(i) \text{ we have}\\
    &\ \nu_{\delta,\ell} \in \Lambda_{\alpha_{\eps+1}} \setminus \Lambda_{\alpha_\eps} \text{ \underline{or}}\\ 
    &\ \eta_{\delta,\ell} \rest i \in \Lambda_{\alpha_{\eps+1}} \setminus \Lambda_{\alpha_\eps} \text{ \underline{or}}\\
    &\ \Res_i^m(\eta_{\alpha_\eps,\ell}) \in \Lambda_{\alpha_{\eps+1}} \setminus \Lambda_{\alpha_\eps} \big\}.
\end{align*}
Now,
\begin{enumerate}
    \item [$(*)_1$] For all $\eps < \lambda$, $\otp(\base(\eps)) \leq \partial \cdot k$.
\end{enumerate}
Hence easily
\begin{enumerate}
    \item [$(*)_2$] For some stationary $S \subseteq \lambda$, we have a function $f : S \to \lambda$ such that 
    $$
    \zeta \in S \Rightarrow f(\zeta) > \zeta\ \wedge\ \zeta \notin \base(f(\zeta))\ \wedge\ `\zeta \cap \base(f(\zeta)) \text{ is a bounded subset of $\zeta$'.}
    $$
    \item [$(*)_3$] Without loss of generality, for some 
    $\xi_* < \partial$ we have 
    $$
    \zeta \in S \Rightarrow \zeta \cap \base(f(\zeta)) \subseteq \xi_*.
    $$
\end{enumerate}
The rest should be clear.
\end{PROOF}

\mn
\begin{claim}\label{2.7}
$1)$ If $\bfB_1, \bfB$ satisfy the strong $\lambda$-chain condition,
$a^* \in \bfB_1\setminus \{\mathbf{0}_{\bfB_1}\}$, and $\bfB_2$ is the result 
of a $\bfB$-surgery of $\bfB_1$ at $a^*$, \underline{then} $\bfB_2$ satisfies 
the strong $\lambda$-chain condition. If one of $\bfB_1$, $\bfB$ satisfies 
the strong $\lambda$-chain condition, and the other only the $\lambda$-chain 
condition, \underline{then} $\bfB_2$ satisfies the $\lambda$-chain condition.

\sn
$2)$  If $\bfB_2$ is the result of a $\bfB$-surgery of $\bfB_1$ at
$a^*$, \underline{then} $\bfB_1 \leqdot \bfB_2$ (i.e.~$\bfB_1$ is a
subalgebra of $\bfB_2$, and every maximal antichain of $\bfB_1$ 
is a maximal antichain of
$\bfB_2$. This is also called ``$\bfB_2$ is a regular extension of
$\bfB_1$").
\end{claim}

\begin{PROOF}{\ref{2.7}}
Well known (and easy).
\end{PROOF}

\mn
\begin{claim}\label{2.8}
The relation $\lessdot$ between Boolean algebras is a partial order, and if
a sequence $\LL \bfB_i : i < \alpha\RR$ is $\lessdot$-increasing continuous 
\underline{then} $\bfB_0\lessdot \bigcup\limits_{i<\alpha} \bfB_i$, and if each 
$\bfB_i$ satisfies the $\chi$-chain condition (for a regular $\chi$), \underline{then} so does $\bigcup\limits_{i<\alpha} \bfB_i$.

(Similarly for the strong $\chi$-chain condition.)
\end{claim}

\begin{PROOF}{\ref{2.7}}
Well known: Solovay-Tenenbaum \cite{ST} for the $\chi$-chain
condition, and Kunen-Tall \cite[p.179]{KT79} for the strong
$\chi$-chain condition. 
\end{PROOF}

\mn
\begin{claim}\label{2.9}
$1)$  In Construction \emph{\ref{2.4}}, if $|I_i|=\lambda$ (hence
$|\BA_\x(I_i)| = \lambda$ for $i < \alpha$) \underline{then} 
$\|\bfB_i\| = \lambda$ for  $0< i \leq \alpha$.

\mn
$2)$ In \emph{\ref{2.4}}, if each $\BA_\x(I_i)$ satisfies the [strong]
$\chi$-chain condition and $\chi$ is regular \underline{then} 
$\bfB = \Sur_x \LL I_i, a^*_i : i < \alpha\RR$ satisfies the
[strong] $\chi$-chain condition.

\mn
$3)$ Assume that in \emph{\ref{2.4}} we use non-trivial $\bfB_0$ and
$|I_i| = \lambda$. \underline{Then}  $\|\bfB_i\| = \lambda + \|\bfB_0\|$ for all $i \in [1,\alpha]$. 
If in addition $\bfB_0$ satisfies the $\lambda$-cc, and each $\BA_\x(I_i)$
satisfies the strong $\lambda$-chain condition, \underline{then} 
$\bfB$ satisfies the $\lambda$-cc; if in addition
$\bfB_0$ satisfies the strong $\lambda$-cc, \underline{then}  so does $\bfB$.
\end{claim}

\begin{PROOF}{\ref{2.9}}
1) Trivial.

\sn 
2)  By \ref{2.5}, \ref{2.6}, \ref{2.7}, \ref{2.8}.

\sn 
3)  Similar. 
\end{PROOF}

\mn
\begin{lemma}\label{2.10}
$1)$ For the construction in \emph{\ref{2.4}},
$\bfB_\alpha$ is weakly representable in\\
$\cM_{\aleph_0,\aleph_0}\big(\sum\limits_{i<\alpha}I_i\big)$ 
(see Definition \cite[2.1(c),(d)\subref{f2}]{Sh:E59}).

\sn
$2)$ Moreover, $\bfB_\alpha\rest (1-a^*_i)$ is weakly representable
in $\cM_{\aleph_0,\aleph_0}\big(\!\sum\limits_{\substack{j<\alpha\\j\neq i}} I_j \big)$.

\sn
$3)$ If $\bfB_\alpha$ satisfies the $\theta$-chain condition
\underline{then} $\bfB^c_\alpha$ (the completion of $\bfB_\alpha$) can 
be weakly represented in
$\cM_{\theta,\theta}(\sum\limits_{j<\alpha}I_j)$.  This
representation can extend the one from \emph{\ref{2.10}(1)}.

\sn
$4)$ Similarly for \ref{2.10}$(2)$.

\sn
$5)$ If in \emph{\ref{2.4}} we use a non-trivial $\bfB_0$, we have to adapt. 
For example, assume $\bfB_0$ is weakly representable in a relevant way 
(e.g.~for $(1)$, assume $\bfB_0$ is weakly represented in 
$\cM_{\aleph_0,\aleph_0}(J + \sum\limits_{i<\alpha}I_i)$).
\end{lemma}

\sn
\begin{remark}\label{2.10z}
1) Why do we use `weakly representable' and not just `representable?' The point is that in Construction \ref{2.4} we have to choose the order of the elements on which we do the surgery. 

\sn
2) Note that the Boolean algebra $\BA_\x(I_i)$ used in the construction usually satisfies the strong $\lambda$-cc for $\lambda$ uncountable.
\end{remark}

\sn
\begin{PROOF}{\ref{2.10}}
1) Define $f(0)=0$, $f(1)=1$. Given $b\in \bfB_\alpha$ not equal to 0 or 1, say
that $b$ first appears in $\bfB_{i+1}$. 

Say
\[
b = \big(b',\textstyle\bigcup\limits_{j<m} (c_j\cap d_j) \big)
\]
\mn
with $b'\in \bfB_i\rest (-a^*_i)$ and $c_j \in \bfB_i\rest a^*_i,d_j \in \BA_\x(I_i)$. 
Say (by induction hypothesis) $f(b')=x'$, $f(c_j)=x_j$,
$f(a^*_i)=x$, $d_j=\sigma_j(x_{\eta_0},\ldots,x_{\eta_{m-1}})$
where $\sigma_j$ is a Boolean term, and $\eta_0,\ldots,\eta_{m-1}\in I_i$. 

Then we set
\[
\begin{array}{c}
f(b)=F_k(x,x',x_0,\ldots,x_{m-1},\eta_0,\ldots,\eta_{m-1}),\\
k\ \text{ codes }\ \LL m,n,\sigma_0,\ldots,\sigma_{m-1}\RR,
\end{array}
\]
where $F_k$ is a suitable function symbol. Thus, $f(b)$ codes all the
relevant information about $b$.

\mn
2)  We may assume that $a^*_i\neq 0,1$. We go exactly as
in (1) up to $\bfB_i$. For $\alpha>i$, we use $(-a^*_i)$ in place of 1,
and working always with $\bfB_\alpha\rest (-a^*_i)$. Note that no
terms involving $I_i$ appear then.

\mn
3) For each $a\in \bfB^c_\alpha$ we can fix $\kappa<
\theta$ and a sequence $\LL b_\gamma:\gamma<\kappa\RR$ of elements
of $\bfB_\alpha$ such that $a=\bigcup\limits_{\gamma<\kappa}b_\gamma$. Then
let $f_\alpha=F(\sigma_\gamma:\gamma<\kappa)$, where $f(b_\gamma)=
\sigma_\gamma$ for all $\gamma<\kappa$.

\mn 
4-5) Similarly. 
\end{PROOF}

\mn
\begin{remark}\label{2.10Anew}
1) In \ref{2.12}-\ref{2.13} below 
we can omit the `weak' from representation and the
`strong' from unembeddability.

\sn
2)  Why weakly represented? As the order of the construction 
and the choice of the $a^*_i$ play a role in the 
definition, we can overcome this in various ways but there is  no real 
reason for doing this
\end{remark}

\mn
\begin{lemma}\label{2.11}
$1)$ Suppose $I\in K^\partial_\tr$ is strongly 
$(\partial,\aleph_0,\psi_\tr)$-unembeddable into $J\in K^\partial_\tr$, 
and $\bfB$ is a Boolean algebra weakly representable in 
$\cM_{\partial,\aleph_0}(J)$. \underline{Then} $\BA_\tr(I)$ 
is not embeddable into $\bfB$.

\sn
$2)$ Suppose $\partial = \kappa$, $I\in K^\partial_\tr$ is strongly $(\mu,\kappa,\psi_\tr)$-unembeddable 
into $J$ for embeddings which are strongly finitary on
$P^I_\partial$, and $\bfB$ is a Boolean algebra weakly represented in 
$\cM_{\mu,\kappa}(J)$. \underline{Then} $\BA_\tr(I)$ is not embeddable into $\bfB$.
\end{lemma}

\begin{PROOF}{\ref{2.11}}
1) Let $g : B \to \cM_{\partial,\aleph_0}(J)$ be
a weak representation of $\bfB$ into $\cM_{\partial,\aleph_0}(J)$ (with
the well-ordering $<^*$), and $h$ be an embedding of $\BA_\tr(I)$ into
$\bfB$. For $\eta\in I$ define $f(\eta) = g(h(x_\eta))$. As $I$ is strongly
$(\partial,\aleph_0,\psi_\tr)$-unembeddable into $J$, there are $\nu_1$,
$\nu_2$, $\eta$, $n$ such that 
$\eta \in P^I_\omega$, $\nu_1 = \eta \rest(n+1)$, $\nu_1 \rest n = \nu_2\rest n$, 
$\nu_2(n) <^J_1 \nu_1(n)$, $\lh(\nu_1) = \lh(\nu_2) = n+1$, and
\[
\big\LL f(\nu_1),f(\eta)\big\RR \approx \big\LL f(\nu_2),f(\eta)\big\RR
\mod \big(\cM^*_{\aleph_0,\aleph_0}(J),{<^*} \big).
\]

\mn
Hence (because $g$ is a weak representation)
\[
h(x_\eta) < h(x_{\nu_1}) \Leftrightarrow h(x_\eta) < h(x_{\nu_2})\quad 
\text{ (in $\bfB$).}
\]

\mn
But $h$ is an embedding, hence $x_\eta<x_{\nu_1} \Leftrightarrow x_\eta<x_{\nu_2}$ in
$\BA_\tr(I)$, contradicting the definition of $\BA_\tr(I)$.

\sn 
2) Similar. 
\end{PROOF}

\mn
\begin{lemma}\label{2.12}
$1)$ Suppose $I$, $J\in K^\partial_{\ptr}$ and $I$ is standard, strongly
$(\mu,\kappa,\psi_\ptr)$-unembeddable into $J$ by $f$ strongly finitary
on $P^I_\partial$. If $\bfB$ is a Boolean algebra weakly representable in
$\cM_{\partial,\aleph _0}(J)$ (say, by $g$), $\bfB \subseteq \bfB_1$
dense\footnote{
    E.g.~$\bfB_1$ is the completion of $\bfB$ --- the case that  interests us.
} 
in $\bfB_1$, and $g_1$ extends $g$ and is a weak representation 
of $\bfB_1$ in $\cM_{\mu,\kappa}(J)$, \underline{then} $\BA_\ptr(I)$ is not 
embeddable into $\bfB_1$.

\mn
$2)$ Analogously for $K^\omega_{\tr(h)}$, $\psi_{\tr(h)}$, $\BA_{\tr(h)}(-)$
(for $h\in {}^\omega (\omega\setminus 2)$) and $K^\omega_{\tr(h)}$,
$\psi_{\tr(h)}$, $\BA_{\tr(h,g)}(-)$.

\mn
$3)$ If $I \in K^\partial_{\tr(h)}$ is standard and 
$(\partial^{\aleph_0},\aleph_1)$-super$^\vr$ 
unembeddable into $J\in K^\partial_{\tr(h)}$, $\bfB$ is weakly 
represented in $\cM_{\partial,\aleph_0}(J)$ and satisfies the ccc 
(for example, $\rang(h)\subseteq [3,\omega)$) \underline{then}  $\BA_{\tr(h)}(I)$ is not
embeddable into the completion of $\bfB$.
\end{lemma}

\begin{PROOF}{\ref{2.12}}
1) Suppose $\bff$ is an embedding of $\BA_\ptr(I)$ into $\bfB_1$. For
$\eta\in I$, define $f(\eta)$ as follows: if $\lh(\eta) < \partial$ then 
$f(\eta) = g_1(\bff(x_\eta))$, whereas if $\lh(\eta) = \partial$, choose 
$a_\eta \in \bfB$, $0 < a_\eta \leq \bff(x_\eta)$ (possible as $\bfB$ 
is dense in $\bfB_1$) and let $f(\eta) = g(a_\eta)$.  As $I$ is strongly
$(\mu,\kappa,\psi_\ptr)$-unembeddable into $J$ by a function $f$ 
which is strongly finitary on $P^I_\partial$, there are $\nu_1$, 
$\nu_2$, $\eta$, $n$ such that $\eta\in P^I_\partial$, 
$\nu_1 = \eta \rest n \caret \LL\alpha\RR$, 
$\nu_2 = \eta \rest n \caret \LL\beta\RR$, $\eta(n) = \LL\alpha,\beta\RR$, 
$\alpha < \beta$, and
\[
\big\LL f(\nu_1),f(\eta)\big\RR \approx \big\LL f(\nu_2),f(\eta)\big\RR 
\mod (M_{\mu,\kappa}(J),<^*).
\]
Hence, as $g_1$ is a weak representation
\[
\begin{array}{lrcl}
(*)\quad&\bfB_1\models \bff(a_\eta)<\bff(x_{\nu_1})&
\Leftrightarrow& \bfB_1\models
\bff(a_\eta)<\bff(x_{\nu_2}),\\
&\bfB_1\models \bff(a_\eta)\cap \bff(x_{\nu_1})=0 &
\Leftrightarrow& \bfB_1\models
\bff(a_\eta) \cap \bff(x_{\nu_2})=0.
  \end{array}
\]
But in $\BA_\ptr(I)$, $x_{\nu_1}\geq x_\eta$, $x_{\nu_2} \cap
x_\eta=0$. Hence, as $\bff$ is an embedding,
\[
\bfB_1 \models ``\bff(x_{\nu_1}) \ge \bff(x_\eta) \wedge 
\bff(x_{\nu_2}) \cap \bff (x_\eta)=0".
\]

\mn
But $0 < a_\eta \le \bff(x_\eta)$, so $\bff(x_{\nu_1})\ge a_\eta$, 
$\bff(x_{\nu_2}) \cap a_\eta = 0$, a contradiction to $(*)$ above.

We have proved that $\BA_\ptr(I)$ is not embeddable into $\bfB_1$.

\mn 
2) Similar proof (the extra details appear in the proof of part (3)).

\mn 
3) 
Let $\partial_1 \defeq \partial^{\aleph_0}$. Assume toward contradiction that $\bff$ 
is an embedding of $\BA_{\tr(h)}(I)$ into $\bfB_1$, the completion of 
$\bfB$. Let $g : \bfB \to \cM_{\partial,\aleph_0}(J)$ be a weak 
representation (say, for the well-ordering $<^*$) of $\cM_{\partial,\aleph_0}(J)$ 
which respects subterms. So by \ref{2.4}(3) there is 
$g_1 : \bfB_1 \to \cM_{\partial_1,\aleph_1}(J)$ which extends $g$ 
and is a weak representation of $\bfB_1$ in 
$(\cM_{\partial_1,\aleph_1}(J), <^*)$.
Choose a function $f : I \to \cM_{\partial_1,\aleph_1}(J)$ as in 
the proof of part (1). Let $x = \LL h,g,g_1,f,I,J,\bfB,\bfB_1\RR$ and let
$\chi$ be large enough.

As it is assumed in part (3) that ``$I$ is $(\partial,\aleph_0)$-super$^\vr$ 
unembeddable into $J$," there are $M,\eta$ as in $(*)'$ of Definition \ref{1.4}(2). 
Let $f(\eta) = \sigma_\eta (x_{\nu_{\eta,0}},\ldots, x_{\nu_{\eta,k(\eta)-1}})$, 
where $\nu_{\eta,k}\in J$ are pairwise distinct for 
$k < k(\eta) < \omega$. 
For each $k$ let $i_k \le \partial$ be maximal such that 
$\nu_{\eta,k} \rest i_k \in M$: 
it exists by clause (v) in $(*)'$ of Definition \ref{1.4}(2). 

If $i_k < \lh(\nu_{k,\ell})$, then for each $m < h(i_k)$ let 
$\nu^*_{k,m} = (\nu_{\eta,k}\rest i_k) \caret \LL s_{k,m}\RR\in M$ be $<^J_1$-minimal 
such that $\Res^m_{i_k} (\nu_{\eta,k}) <^J_1 \nu^*_{k,m}$. 
Clearly it exists, except when $\Res^m_{i_k} (\nu_{\eta,k})$ is $<^J_1$-above 
every $\{(\nu_{\eta,k}\rest i_k) \caret \LL s\RR : s \in M\}$; in that case we 
let $s_{k,m} = \infty$ with the obvious conventions. 

{Let $\bar \nu \defeq \LL \nu_k : k < k(\eta)\RR$.} We define
$$Y^* = \big\{\bar\nu : \bar\nu\text{ is similar in $J$ to } \LL \nu_{\eta,0},\ldots, \nu_{\eta,k(\eta)-1}\RR \text{ over } Z^*\big\}$$
where $Z^* = \big\{\nu_{\eta,k} : \nu_{\eta,k}\in M\big\} \cup 
\big\{\nu^*_{k,m} : k < k(\eta),\ m < h(k(\eta))\big\}$.
Clearly $Z^*$ is a finite subset of $M$. We define a filter $D$ on $Y^*$: 
$Y\in D$ \underline{iff} there are $\nu'_{k,m} <^J_1 \nu^*_{k,m}$ for all 
relevant $k,m$ such that if $\LL \nu''_k : k < k(\eta)\RR$ satisfies 
$\nu'_{k,m} \leq^J_1 \nu''_{k,m}$ for all relevant $k$ and $m$
then $\LL \nu''_k : k < k(\eta)\RR \in Y$. 

Clearly $(Y^*,D)\in M$, and by weak 
representability the following function $f_1$ belongs to $M$:
$$
\dom (f_1) = \{\varrho\in I : \lh(\varrho) < \partial\},\quad \rang (f_1) \subseteq \{0,1\}, \ \ \text{ and}
$$ 
$$
f_1 (\varrho)=\begin{cases} 1 & 
    \big\{\bar\nu \in Y^* : \BA_{\tr(h)}(J) \models \bff(x_\varrho)\geq \sigma_\eta(x_{\nu_0}, \ldots,x_{\nu_{k(\eta)-1}}) \big\} \in D\\ 
    &\text{ 
    {\underline{iff} that set} is} \neq \varnothing \mod D\\
    0 & \text{otherwise.}
\end{cases}
$$

Recall that $\sigma_\eta$ is a $\tau_{\partial,\aleph_0}$-term, hence it is
$\in M$. So by the choice of $M$ and $\eta$, 
for unboundedly many $i$, (as $\tau=\tau_0$; see Definition \ref{1.4}), 
we have that the truth values of 
$$
\BA_{\tr(h)} (J) \models \bff(x_{\Res^\ell_i(\eta)}) 
\geq \sigma_\eta (x_{\nu_{\eta,0}},\ldots, x_{\nu_{\eta,k(\eta)-1}})
$$
are the same for all $\ell<h(i)$. As $\bff$ is an embedding, 
$$
\bfB_1 \models \bff(x_\eta)\geq \sigma_\eta (x_{\nu_{\eta,0}},\ldots, x_{\nu_{\eta,k(\eta)-1}}) > 0,
$$
and $\BA_{\tr(h)}(I) \models x_{\Res^0_i(\eta)}\geq x_\eta$, we have 
$$
\bfB_1 \models ``\bff(x_{\Res^0_i(\eta)}) \geq \bff(x_\eta) \geq 
f(\eta) = \sigma_\eta(x_{\nu_{\eta,0}},\ldots x_{\nu_\eta,k(\eta)-1}) > 0."
$$
So $f_1(\Res^0_i(\eta)) = 1$, hence by the choice of $n$ we have 
$\ell < h(i) \Rightarrow f_1(\Res^\ell_i(\eta)) = 1$. 
So $\bfB_1 \models ``\bigcap\limits_{\ell<h(i)} 
\bff (x_{\Res^\ell_i(\eta)}) \cap f(x_\eta) > 0"$, but $\bff$ is an 
embedding and 
$$
\BA_{\tr(h)} (J) \models ``0<f(\eta) \leq \bff (x_\eta)"
$$
hence 
$\BA_{\tr(h)} (I) \models \bigcap\limits_{\ell<i} x_{\Res^\ell_i(\eta)} 
\cap x_\eta > 0$, contradicting the definition of $\BA_{\tr(h)}(I)$. 
\end{PROOF}

\mn
\begin{conclusion}\label{2.13}
Suppose $\lambda>\aleph_0$. \underline{Then}:
\mn
\begin{enumerate}[$(A)$]
    \item  There is a rigid Boolean algebra $\bfB$ satisfying the $\aleph_1$-chain condition $\lambda$.
\sn
    \item Moreover, if $a,b \in \bfB$ are $\neq 0$ with $a-b \neq 0$, \underline{then} 
    $\bfB \rest a$ cannot be embedded into $\bfB\rest b$ (hence $\bfB$ has no one-to-one endomorphism $\neq\id$).
\sn
    \item  Moreover, we can find such $\bfB_i$ (for $i<2^\lambda$) with 
    $|\bfB_i|= \lambda$; and if $a\in \bfB_i$, $b \in \bfB_j$ with $i\neq j$ or 
    $a-b\neq 0$ then $\bfB_i \rest a$ cannot be embedded into $\bfB_j\rest b$.
\end{enumerate}
\end{conclusion}

\begin{PROOF}{\ref{2.13}}
We leave it to the reader as the next proof is similar (but here we should use 
$(\lambda,\lambda,\aleph_0,\aleph_0)$-$\psi_\tr$-bigness, Theorem 
\cite[2.20\subref{7.11}]{Sh:331}, and $\x=\tr$ instead of 
$(\lambda,\lambda,2^{\aleph_0},\aleph_1)$-$\psi_{\tr(h)}$-bigness,  
\cite[1.11\subref{7.6}]{Sh:331}, and $\x=\tr(h)$ there, respectively. 
(Also, we have dealt with it in \cite[2.16\subref{2.7}]{Sh:E59}).
\end{PROOF}

\mn
\begin{conclusion}\label{2.14}
$1)$ There is a complete Boolean algebra $\bfB$ satisfying the ccc,
having density $\lambda$ 
\[
\text{ (in fact, } a \in \bfB \setminus \{0\}\ \Rightarrow\ \bfB \rest a
\text{ has density }\lambda, \text{ so } |\bfB| = \lambda^{\aleph_0}),
\]
and monorigid (i.e.~every one-to-one
endomorphism is the identity) p\underline{rovided that}:
\begin{enumerate}
    \item[$(*)_1$]  $K^\partial_\ptr$ has the full strong 
    $(\lambda,\lambda,2^{\aleph_0},\aleph_1)$- $\psi_{\ptr}$-bigness property for $f$ strong finitary on $P_\omega$, by standard atomically $({<}\,\aleph_1)$-stable 
    $I \in K^\partial_\ptr$.
\end{enumerate}

\mn
$2)$ We can replace $(*)_1$ by $(*)_2 \vee (*)_3 \vee (*)_4$, where for some $h\in {}^\omega(\omega\setminus 3)$:
\mn
\begin{enumerate}
    \item[$(*)_2$] $\lambda$ is as in \ref{1.7}$(1)$ \underline{or}
\sn 
    \item[$(*)_3$] $K^\partial _{\tr(h)}$ has the full strong 
    $(\lambda,\lambda,2^{\aleph_0},\aleph_1)$-$\psi_{\tr(h)}$-bigness property \underline{or}
\sn
    \item[$(*)_4$] $K^\partial_{\tr (h)}$ has the full super$^\vr (\lambda,\lambda,2^{\aleph_0},\aleph_1)$-bigness property.
\end{enumerate}

\mn
$3)$ Moreover, we can find such $\bfB_i$ (for $i<2^\lambda$)
satisfying the following: if $a \in \bfB_i\setminus\{0\}$, $b\in \bfB_j\setminus \{0\}$,
$[i\neq j \vee (i = j\wedge a-b \neq \mathbf{0}_{\bfB_i})]$, 
then $\bfB_i\rest a$ cannot
be embedded into $\bfB_j \rest b$.
\end{conclusion}

\begin{PROOF}{\ref{2.14}}
We first prove parts (1) and (2). For part (1) let $h\in {}^\omega\omega$ be
constantly 2. First note that if $f$ is a one-to-one endomorphism
$\ne \id$ of any Boolean algebra $\bfB$, then there is an element
$a \ne 0$ with $a\cap f(a) = 0$. First, choose $x$ with $x\neq f(x)$. 
If $x \cap -f(x)\neq 0$ we can take $a = x \cap -f(x)$; if 
$-x \cap f(x) \ne 0$ we can take $a=-x\cap f(x)$. 
Hence for {(1)} and {(2)} we only need to find $\bfB$ of power $\lambda$
such that if $a,b\in \bfB$ are non-zero and $a-b\neq 0$ (and even $a\cap b=0$), 
then $\bfB \rest a$ cannot be embedded in $\bfB\rest b$.

Now let $\LL I_\alpha : \alpha < \lambda\RR$ exemplify the full strong
$(\lambda,\lambda,2^{\aleph_0},\aleph_1)$-$\psi_{\tr(h)}$-bigness property
for $f$ strongly finitary on $P_\omega$; such a sequence exists by  $(*)_1$ or 
$(*)_2$ or $(*)_3$ or $(*)_4$ by \ref{1.7}(1), \ref{1.6B} for any 
$h\in {}^\omega(\omega\setminus 3)$. Let $\bfB=\Sur_x\LL I_\alpha$,
$a^*_\alpha : \alpha < \lambda\RR$ be as in the construction \ref{2.4}
for $x = \tr(h)$, such that $\bfB \setminus \{0\} = \{a^*_\alpha : \alpha < \lambda\}$. Then by \ref{2.9}(1), $|\bfB|=\lambda$. By \ref{2.6}(3), \ref{2.6}(4), each 
$\BA_{\tr (h)} (I_\alpha)$ satisfies the strong $\aleph_1$-cc, hence by 
\ref{2.9} the Boolean algebra $\bfB$ satisfies the $\aleph_1$-chain condition. 
Let $\bfB^*$ be its completion.  Now let $a,b \in \bfB^*$ be non-zero, with 
$c = a-b \neq 0$.  Toward contradiction, suppose $f$ is an embedding of 
$\bfB^* \rest a$ into $\bfB^* \rest b$. Then $f(c) \cap c = 0$, and 
$f \rest (\bfB\rest c)$ is an embedding of $\bfB^* \rest c$ into $\bfB^* \rest f(c)$. 
But $\bfB$ is dense in $\bfB^*$ hence $a^*_\alpha \le c$ for some $\alpha$, 
hence $\BA_{\tr(h)}(I_\alpha)$ is embeddable into $\bfB^* \rest c$, hence
into $\bfB^* \rest f(c)$, hence into 
$\bfB^* \rest (-c) = \bfB^* \rest (-a^*_\alpha)$. But by \ref{2.10}(3), 
$\bfB\rest (-a^*_\alpha)$ is weakly representable in 
$\cM^*_{2^{\aleph_0},\aleph_1} \big(\sum\limits_{\beta\neq\alpha,\beta<\lambda} I_\beta\big)$. This
contradicts \ref{2.12} when we assume $(*)_4$.

For part (3) let $\LL I_{\alpha,\beta} : \alpha, \beta < \lambda\RR$
rename $\LL I_\alpha : \alpha < \lambda\RR$. We shall choose, for
$\xi < 2^\lambda$, functions $f_\xi$, $g_\xi$ from $\lambda$ to $\lambda$ and
$A_\xi \in [\lambda]^\lambda$ such that $g_\xi$ is one-to-one, $\rang(f_\xi) = A_\xi$,
$(\forall\alpha \in A_\xi)(\exists^\lambda \beta < \lambda)[f_\xi(\beta) = 1]$, 
and $\xi_1 \neq \xi_2 \Rightarrow\ A_{\xi_1} \not\subseteq A_{\xi_2}$. For 
$\xi < 2^\lambda$, let $\bfB^\xi$ be constructed as 
$\Sur_x\LL I_{f_\xi(\alpha),g_\xi(\alpha)},a^\xi_\alpha : \alpha < \lambda\RR$. 
For simplicity, assume that for some $\xi$, for every $a\in \bfB^\xi \setminus \{0\}$
and $\zeta\in A_\xi$, we have $a^\xi_\alpha = a$ and $f_\xi(\alpha) = \zeta$.  
Let $\bfB^{\xi,*}$ be the completion of $\bfB^\xi$. As $g_\xi$ is one-to-one clearly
$\bfB^\xi$ satisfies the demand in (2), and as 
$\xi \neq \zeta < 2^\lambda\ \Rightarrow\ A_\xi \not\subseteq A_\zeta$ 
the demands in (3) also hold.
\end{PROOF}

\begin{conclusion}\label{2.15}
$1)$ For $\lambda > \aleph_0$, there is a Boolean algebra $\bfB$ of
cardinality $\lambda$ with no non-trivial endomorphism onto itself. 
Moreover, it is Bonnet-rigid (defined below).

\sn
$2)$  We can find such $\bfB_i$ (for $i<2^\lambda$) such that
for $i$, $j < 2^\lambda$, $a\in \bfB_i\setminus\{0\}$, $b \in
\bfB_j\setminus\{0\}$ there is no embedding of $\bfB_i\rest a$ into a
homomorphic image of
$\bfB_j \rest b$ except when $i = j \wedge a \leq b$.
\end{conclusion}

\noindent
We prove it later.

\begin{remark}\label{2.15Anew}
We shall use Boolean algebras built from cases of $\BA_\trr(I)$ (see Definition
\ref{2.1}(4)) hence {they have} no long chains. We can go in the inverse direction
using Boolean algebras built from orders --- using, for example, $\LO(I)$ the
linear order with elements $\{x_\eta,y_\eta : \eta \in I\}$ such that:
\mn
\begin{enumerate}
    \item  $\lh(\eta) < \omega$ implies $x_\eta < y_\eta$,
    $y_{\eta \caret \LL\alpha\RR} < x_{\eta \caret \LL\beta\RR}$ for 
    $\alpha < \beta$, and $x_{\eta\rest n} < x_\eta < y_\eta < y_{\eta\rest n}$ 
    for $n < \lh(\eta)$.
\sn
    \item $\lh(\eta) = \omega$ implies $x_{\eta\rest n} < x_\eta = 
    y_\eta < y_{\eta\rest n}$ for $n < \omega$.
\end{enumerate}
In such cases we need a parallel to Lemma \ref{2.19}, which is true.
\end{remark}

\noindent
We make some preparations to the proof of \ref{2.15}.

\begin{definition}\label{2.16}
A Boolean algebra $\bfB$ is called \emph{Bonnet--rigid} \underline{iff}  there are no
Boolean algebra $\bfB'$ and homomorphisms $\bff_\ell : \bfB \to \bfB'$ 
(for $\ell=0,1$) such that $\bff_0$ is one-to-one and $\bff_1$ is onto 
$\bfB'$, except when $\bff_0=\bff_1$.
\end{definition}

\begin{observation}\label{2.17}
$1)$ If $\bfB$ is Bonnet--rigid \underline{then} it has no onto endomorphism
$\neq\id_\bfB$.

\sn
$2)$  A Boolean algebra $\bfB$ is Bonnet--rigid \underline{if}:
\begin{enumerate}
    \item[$(*)$]  For no disjoint non-zero $a,b\in \bfB$ is there an 
    embedding of $\bfB \rest a$ into $a$ homeomorphic image of $\bfB\rest b$.
\end{enumerate}
\end{observation}

\begin{PROOF}{\ref{2.17}}
1) Otherwise choose $\bfB'=\bfB$, $\bff_0$ the identity, and
$\bff_1$ the given endomorphism.

\sn
2) Towards contradiction, assume $\bff_\ell : \bfB \to \bfB'$
(for $\ell=0,1$) contradict Bonnet--rigidity. First, suppose $\bff_1$ is not
one-to-one, so for some $a \in \bfB$, $a \neq 0$, $\bff_1(a)=0$.

For any $b\in \bfB$, $\bff_1(b-a) = \bff_1(b)-\bff_1(a) = \bff_1(b)$. 
So $\bfB'$ is a homomorphic image of $\bfB \rest (1-a)$ and 
$\bfB\rest a$ can be embedded into it, so we are finished.

Second, assume $\bff_1$ is one-to-one. Then $\bff_1$ is an 
isomorphism from $\bfB$ onto $\bfB'$
hence $\bff_1^{-1}\bff_0:\bfB\to \bfB$ is an embedding 
(well-defined as $f_1$ is one-to-one and onto). It is not the
identity (otherwise $\bff_0=\bff_1$) so for some $a\in 
\bfB$, the elements $a,\bff^{-1}_1 \bff_0(a)$ are 
disjoint and non-zero; choose $b=\bff^{-1}_1 \bff_0(a)$.
\end{PROOF}

\noindent
To prove \ref{2.15}, we shall use $\BA_\trr(I)$ (see Definition \ref{2.1}(4)). Note:

\begin{claim}\label{2.18}
$1)$  The only atoms of $\BA_\trr(I)$ are $x_\eta$, where $\eta$ ($\in I$)
has no immediate successor, or at least
\begin{enumerate}
    \item[$(*)$] For all $\nu_1,\nu_2 \in I$, we have $\eta \lhd \nu_1 \wedge \eta \lhd \nu_2\ \Rightarrow\ \nu_1, \nu_2\text{ are $\lhd$-comparable.}$
\end{enumerate}

\mn
$2)$ The set $\{x_\eta : \eta \in I\}$ is a dense subset of $\BA_\trr(I)$.
\end{claim}

\begin{PROOF}{\ref{2.18}}
Check. 
\end{PROOF}

\begin{lemma}\label{2.19}
If $\bfB$ is a homomorphic image of $\bfB_0=\BA_\trr(I)$, \underline{then}  $\bfB$ is
isomorphic to some $\BA_\trr(J)$, $J$ weakly representable in $\cM_{\aleph_0,
\aleph_0}(I)$ hence $\bfB$ is weakly representable in
$\cM_{\aleph_0,\aleph_0} (I)$.
\end{lemma}

\begin{PROOF}{\ref{2.19}}
So let $\bfJ$ be an ideal of $\bfB_0$ such that $\bfB$ is
isomorphic to $\bfB_0/\bfJ$. Let
\[
I_1 = \{\eta \in I : x_\eta \notin \bfJ\};
\]
$I_1$ is an approximation to $J$. (Clearly $I_1$ is closed under initial
segments by \ref{2.1}(4)(b).) Let
\begin{align*}
A_0 =  \big\{\eta\in I_1 : &\ \eta\text{ has $<\aleph_0$ immediate successors in }I_1, \text{ say }\\
&\ \eta \caret \LL\alpha_\ell\RR\ \text{ for } \ell < m,\ \text{ and
}\big(x_\eta - \textstyle\bigcup\limits_\ell x_{\eta \caret \LL\alpha_\ell\RR} \big) \in
\bfJ \big\},\\  
A_1 = \big\{\eta\in I_1 : &\ \eta \text{ has $<\aleph_0$ immediate successors in }I_1, 
\text{ say }\\
&\ \eta \caret \LL\alpha_\ell\RR\ \text{ for } \ell<m,\ \text{ and
} \big(x_\eta - \textstyle\bigcup\limits_\ell x_{\eta \caret \LL\alpha_\ell\RR} \big) \notin \bfJ \big\},\\
A_3 = \big\{(\eta,\nu) : &\ \eta\in A_0,\ \eta\lhd\nu\in I_1,\ \lh(\nu) \text{ is limit},\ x_\eta-x_{\nu\rest i}\in \bfJ,\\
&\text{ when }\lh(\eta)\leq i < \lh(\nu)\text{ and for no }\eta'
 \lhd\eta\text{ does }(\eta',\nu)\\
&\text{ have those properties}\big\},\\
\end{align*}
and let
$A_4 = \big\{(\eta,\nu)\in A_3 :  x_\eta - x_\nu \notin \bfJ \big\}$.

\mn
Now for $\eta\in I$ let $\alpha_\eta = \min\{\alpha : \eta \caret
\LL\alpha\RR\notin I\}$.

Put
\[
J = I_1 \cup\ \{\eta \caret \LL\alpha_\eta\RR:\eta\in A_1\} \cup\ 
\{\eta \caret \LL\alpha_\eta+1\RR: (\eta,\nu)\in A_4\}.
\]

\mn
Now $\BA_\trr(J)$ is isomorphic to $\bfB$, and the lemma should be clear.
\end{PROOF}

\noindent
Now we can turn to

\begin{PROOF}{\ref{2.15}}
\textbf{Proof of \ref{2.15}}:

\noindent
1) Let $\LL I_\alpha:\alpha<\lambda\RR$ exemplify that
$K^\partial_\tr$ has the full strong
$(\lambda,\lambda,\aleph_0,\aleph_0)$-bigness property, 
$I_\alpha$ standard. 

Without loss of generality:
\begin{enumerate}
    \item[$(*)_1$] $\alpha \neq \beta\ \Rightarrow\ I_\alpha \cap I_\beta = 
    \{\LL\ \RR\}$
\sn
    \item[$(*)_2$]  If $\nu\in I_\alpha$ then for some $\eta$ we have 
    $\nu \unlhd \eta \in I_\alpha$ and $\lh(\eta) = \omega$.
\end{enumerate}
\mn
We construct as in \ref{2.4}, using $\BA_\trr(I_\alpha)$
(i.e.~$x=\trr$ there) but making the surgeries on atoms only, 
getting $\bfB = \Sur \LL I_\alpha,a^*_\alpha : \alpha < \lambda\RR$. 
Looking at the construction, it is clear that $\bfB = \BA_\trr(I^*)$, where
\begin{align*}
    I^* = \big\{\eta_1 \caret \eta_2 \caret \ldots \caret \eta_n :&\ n < \omega,\ \eta_\ell \in I_{\alpha_\ell} \text{ for some $\alpha_\ell < \lambda$, and for } \ell < n\\ 
    &\text{ we have } \lh(\eta_\ell)=\omega \text{ and $a^*_{\alpha_\ell + 1}$ is } x_{\eta_\ell} \big\}.
\end{align*}

By \ref{2.17}(2), it suffices to prove:
\mn
\begin{enumerate}
    \item[$(**)$]  If $a,b$ are disjoint non-zero and $\bfB'$ is a homomorphic image of $\bfB\rest b$ then $\bfB\rest a$ cannot be embedded into $\bfB'$.
\end{enumerate}
\mn
Suppose $(**)$ fails and $a,b$, $\bfB'$ exemplify this. By Claim
\ref{2.18} and $(*)$, there is $\eta\in I^*$ with $x_\eta \leq a$ and
$\lh(\eta)$ limit, and let $\alpha$ be such that $a^*_\alpha = x_\eta$. 
Clearly $\bfB'$ is also a homomorphic image of $\bfB\rest (1-x_\eta)$, 
hence by \ref{2.19} it is weakly representable in 
$\cM^*_{\aleph_0,\aleph_0} \big(\sum\limits_{j<\lambda,j\neq\alpha} I_j\big)$ 
and $\bfB'\cong \BA_\trr(I^+)$ for some $I^+$ weakly representable in
$\cM_{\aleph_0,\aleph_0} \big(\sum\limits_{j<\lambda,j\neq\alpha}I_j\big)$.

We can conclude:
\mn
\begin{enumerate}
    \item[$(***)$]  $\BA_\trr(I_\alpha)$ is weakly representable in
$\cM_{\aleph_0,\aleph_0}$ ($\sum\limits_{j<\lambda,j\neq\alpha}I_j$).
\end{enumerate} 
\mn
But from this the contradiction is trivial (we could avoid the ``weakly").

\sn 
2) No new point.
\end{PROOF}

\newpage

\section{Arbitrary length of a Boolean Algebra with no small infinite
homomorphic image}

We recall the definition of the length 
(and length$^+$) of a Boolean algebra (Definition \ref{5.0A}).
Our aim is to construct a Boolean algebra of cardinality 
continuum with no infinite homomorphic image of smaller cardinality. 
Toward this, for a Boolean algebra $\bfB$ of cardinality 
$\leq 2^{\aleph_0}$ satisfying the ccc, an $\omega$-sequence 
$\LL a_n : n < \omega\RR$ of pairwise disjoint members of 
$\bfB \setminus \{\mathbf{0}_\bfB\}$, and $I\in K^\omega_{\tr(h)}$, we define in
Definition \ref{5.1} an extension $\bfB' = \ba [\bfB,\bar a,I]$ of $\bfB$. 
We shall use it for $h$ with $\LL h(n) : n < \omega\RR$ going to infinity. 
The properties we need are that $\bfB \lessdot \bfB'$, 
$\|\bfB'\| \le 2^{\aleph_0}$, and $\bfB'$ satisfies the ccc.\footnote{
    See \ref{5.2}(1),(3), \ref{5.3}, and inside the proof of \ref{5.4}.
} 
Moreover, a stronger version of $\bfB\lessdot \bfB'$ 
holds (see \ref{5.2}(5)). 

Also, if $\bff$ is a homomorphism from $\bfB'$ into any 
Boolean algebra $\bfB''$ satisfying  $n<\omega\Rightarrow \bff(a_n)>0$ 
(in $\bfB''$) then $\bfB''$ has at least $2^{\aleph_0}$ elements 
(see inside the proof of \ref{5.4}). Theorem \ref{5.4} is the main result: 
if $\mu\in [\aleph_1,2^{\aleph_0}]$ then some ccc  Boolean algebra $\bfB$ 
of cardinality $2^{\aleph_0}$ and length $\mu$ has no infinite homomorphic 
image of cardinality $<2^{\aleph_0}$. For this we take care of 
every antichain $\LL a_n : n < \omega\RR$ by an extension 
$\ba[-,\bar a,I]$.  We start with a ccc 
Boolean algebra of length and cardinality $\mu$. 
In this framework we need to show that the length 
has not increased by the construction. 
For this we prove, by induction on the length of the construction, that for 
any family of $\mu^+$ finite sequences from the Boolean algebra and 
$m < \omega$, there is a subfamily of $\mu^+$ finite sequences which is an 
indiscernible set of pairwise distinct elements. 

We may like to consider a limit $\mu \in [\aleph_1,2^{\aleph_0})$ 
and ask above that its length is $\mu$ but the supremum is 
not obtained; by a similar construction (of length $2^{\aleph_0}\times \mu$)
we get such a Boolean algebra, provided that $\cf(\mu)$ is uncountable 
(see \ref{5.17}). If $\cf(\mu) = \aleph_0$ this is impossible
(see \ref{5.7bis}). We then generalize the results, replacing 
$\aleph_0$ by any strong limit $\kappa$ of cofinality $\aleph_0$.

\begin{convention}\label{5.0}
In this section, $h$ will be from ${}^\omega (\omega\setminus\{0\})$ 
(and for simplicity ${}^\omega (\omega\setminus\{0,1\})$). 
The constant function $h=2$ actually suffices,\footnote{
    That is, using $K^\omega_\ptr$ as in the proof of \ref{5.4}.
} 
\underline{but} if we would like to have the ccc, we'd better use $h \geq 3$.
\end{convention}

\begin{definition}\label{5.0A}
For a Boolean algebra $\bfB$ let
\[
\begin{array}{rcl}
\mathrm{length}(\bfB)&=&\sup\{|A| : A\subseteq \bfB,\ A\text{ is linearly ordered
by } {<_\bfB}\},\\
\mathrm{length}^+(\bfB)&=&\sup\{|A|^+ : A\subseteq \bfB,\ A\text{ is linearly
ordered by }{<_\bfB}\}.
\end{array}
\]
\end{definition}

\begin{definition}\label{5.1}
For a Boolean Algebra $\bfB^*$, $\bar{a} = \LL a_n : n < \omega\RR 
\subseteq \bfB^* \setminus \{\mathbf{0}_{\bfB^*}\}$ such that
$\bigwedge\limits_{n<m} a_n \cap a_m = 0$, and $I \in K^\omega_{\tr(h)}$, we
define a Boolean Algebra $\ba[\bfB^*,\bar{a},I]$ as follows.

It is freely generated by $\bfB^* \cup \{x_\eta : \eta \in I\}$, except for the
following equalities:
\mn
\begin{enumerate}[(a)]
    \item   All the equalities which $\bfB^*$ satisfies, and $x_\eta\leq \mathbf{1}_{\bfB^*}$.
\sn
    \item  If $n<\omega$ is even, $k=h(n)-1$, $\eta\in P^I_\omega$, 
    $\nu = \eta \rest n$, and $\eta(n) = \LL\alpha_0,\alpha_1,\alpha_2,\ldots,\alpha_{k-1}\RR$, then
    \[
        a_n - \bigcup\limits_{\ell<\frac{k-1}{2}}\left((x_{\nu\caret \LL\alpha_{2\ell}\RR}- x_{\nu \caret \LL\alpha_{2\ell+1}\RR}\right)\leq x_\eta.
    \]
\sn
    \item  If $n<\omega$ is odd, $k = h(n)-1$, $\eta\in P^I_\omega$, 
    $\nu = \eta \rest n$, and $\eta(n) = \LL\alpha_0,\alpha_1,\alpha_2,\ldots,\alpha_{k-1}\RR$ then
    \[
        \Big(a_n \cap \bigcap\limits_{\ell<\frac{k-1}{2}} \left(1-(x_{\nu \caret \LL\alpha_{2\ell} \RR} - x_{\nu \caret \LL\alpha_{2\ell+1}\RR})\right)\Big) \cap x_\eta = 0.
    \]
\sn
    \item  $x_{\LL\ \RR_I}=0$ ($\LL\ \RR_{I}$ is the root of $I$).
\end{enumerate}
\end{definition}

\begin{claim}\label{5.2}
$1)$ For $\bfB^*,\bar a$, $I$ as in Definition \ref{5.1},
$\ba[\bfB^*,\bar{a},I]$ is an extension of $\bfB^*$ 
(so the equalities do not cause
members of $\bfB^*$ to become identified, and of course 
$\mathbf{1}_{\ba[\bfB^{*},\bar a,I]} = \mathbf{1}_{\bfB^*})$.

\sn
$2)$ For $I_1,I_2\in K^\omega_{\tr (h)}$, $\bfB^*,\bar a$ as in
Definition \ref{5.1}, if $I_1\subseteq I_2$ \underline{then}  $\ba[\bfB^*,
\bar a,I_1]$ is a subalgebra of $\ba[\bfB^*,\bar a,I_2]$.

\sn
$3)$ In $(1)$, $\bfB^* \lessdot \ba[\bfB^*,\bar a,I]$.

\sn
$4)$ In $(2)$, if also $I_1\subseteq^* I_2$ (which means that 
$I_1\subseteq I_2$ and\\
$\eta\in P^{I_2}_\omega\setminus I_1\ \Rightarrow\ \bigvee\limits_{n,
\ell}\Res^\ell_n(\eta)\notin I_1$)
\underline{then}  $\ba[\bfB^*,\bar a,I_1]\lessdot\ba[\bfB^*,\bar a,I_2]$.

\sn
$5)$ In $(4)$, for every non-zero $c \in \ba [\bfB^*,\bar a,I_2]$ there is
$d^*$ such that:
\sn
\begin{enumerate}
    \item[$(i)$]   $c \leq d^*\in\ba[\bfB^*,\bar a,I_1]$
\sn
    \item[$(ii)$]  If $0 < b \le d^*$ and $b \in \ba[\bfB^*,\bar a,I_1]$ 
    then $c \cap b \neq 0$.
\end{enumerate}
\end{claim}

\begin{PROOF}{\ref{5.2}}
1) It is a particular case of (2) for $I_1 = \{\LL\ \RR\},\ I_2 = I$.

\sn
2) Let $d^* \in \ba [\bfB^*,\bar a,I_1]\setminus \{0\}$.
We would like to prove that $\ba[\bfB^*,\bar a,I_2] \models d^* \neq 0$; by the definition of these two Boolean algebras (see \ref{5.1}), this suffices.
Clearly, without loss of generality, for some $\alpha_* \leq \omega$ we have:
\begin{enumerate}
    \item[$(*)_1$] $\alpha_* < \omega \wedge d^*\leq a_{\alpha_*}\text{ or } \alpha_* = \omega \wedge (\forall n) [d^* \cap a_n = 0].$
\end{enumerate}

\mn
Now we shall define a function $\bff : \bfB^* \cup \{x_\eta : \eta \in I_2\}
\to \ba[\bfB^*,\bar a,I_1] \rest d^*$, which will map
all the equations appearing in the definition of $\ba[\bfB^*,\bar a,I_2]$ to
ones satisfied in $\ba[\bfB^*,\bar a,I_1] \rest d^*$ 
and maps $d^*$ to itself; this suffices.

Now we define $\bff=\bff^{d^*}$ as follows:
\mn
\begin{enumerate}
    \item[$(*)_2$]
    \begin{enumerate}[(A)]
        \item For $b \in \bfB^*$, $\bff(b) = b \cap d^*$ (or more exactly, the interpretation of $b \cap d^*$ in $\ba[\bfB^*,\bar{a},I_1]$).
\sn
        \item  For $\eta \in I_1$, $\bff(x_\eta) = x_\eta \cap d^*$.
\sn
        \item  If $\eta\in P_\omega^{I_2}$, $\eta\notin I_1$, let
    $$
        \bff(x_\eta)= 
        \begin{cases}
            d^* &\text{ if }\alpha_*\text{ is even (including $\alpha_* = \omega$),}\\
            0 &\text{ if }\alpha_*\text{ is odd (and $<\omega$).}
        \end{cases}
    $$
\sn
        \item  For $\eta \in I_2 \setminus I_1$ such that (C) does not apply, let $\bff(x_\eta) = 0$.
    \end{enumerate}
\end{enumerate}
\mn
Now check: the main point being that the equations in clauses (b)+(c) of
Definition \ref{5.1} hold trivially by the present choice in clause $(*)_2$(C).


\sn 
3) Again, it suffices to prove this for the context of (2); i.e.~to prove (4).

\sn 
4) By part (5).

\sn
5) We can find $\Lambda$ and $I_2'$ such that 
\begin{enumerate}
    \item [$\boxdot_1$]
    \begin{enumerate}
        \item $\Lambda$ is a finite subset of $I_2 \setminus I_1$.
\sn
        \item $I_2' \defeq I_1 \cup \Lambda \cup \big\{ \Res_i^\ell(\eta) : i + 1 < \lh(\eta),\ \ell < h(n), \text{ and } \eta \in \Lambda\big\}$
\sn
        \item $ 
          {c} \in \ba[\bfB^*,\bar a,I_2']$.
    \end{enumerate}
\end{enumerate}
Note that
\begin{enumerate}
    \item [$\boxdot_2$]
    \begin{enumerate}
        \item $\ba[\bfB^*,\bar a,I_1] \subseteq \ba[\bfB^*,\bar a,I_2'] \subseteq \ba[\bfB^*,\bar a,I_2]$
\sn
        \item It suffices to prove the existence of $d^*$ in the case `$I_2' = I_2$.'
    \end{enumerate}
\end{enumerate}
Now
\begin{enumerate}
    \item [$\boxdot_3$]
    \begin{enumerate}
        \item Without loss of generality, if $\eta \in \Lambda \setminus P_\omega^{I_2}$, $i +1 < \lh(\eta)$, and $\ell < h(i)$ \underline{then} $\Res_n^\ell(\eta) \in \Lambda$.
\sn
        \item Let $n_* < \omega$ be such that:
        \begin{enumerate}
            \item If $\eta \in P_n^{I_2} \cap \Lambda$ and $n < \omega$, 
            then $n < n_*$.
\sn
            \item If $\eta \neq \nu \in P_\omega^{I_2} \cap \Lambda$ then $\lh(\eta \cap \nu) < n_*$.
\sn
            \item $c$ belongs to the subalgebra of $\ba[\bfB^*,\bar a,I_2]$ generated by 
            $$\{x_\eta : \eta \in \Lambda\} \cup \ba[\bfB^*,\bar a,I_1].$$
        \end{enumerate}
\sn
        \item Without loss of generality, if $\eta \in P_\omega^{I_2} \cap \Lambda$ and $n < n_*$, $\ell < h(n)$ \underline{then} $\Res_n^\ell(\eta) \in \Lambda$.
    \end{enumerate}
\end{enumerate}
[Why? Just check.]

\sn
Next,
\begin{enumerate}
    \item [$\boxdot_4$] {Let} $\LL \eta_i : i < m\RR$ list $\Lambda$ without repetition such that $\eta_i = \Res_n^\ell(\eta_j) \Rightarrow i < j$.
\end{enumerate}
Now
\begin{enumerate}
    \item [$\boxdot_5$]
    \begin{enumerate}
        \item For $k < m$, let $I^*_{2,k} \defeq I_1 \cup \{\eta_k\} \cup \{\Res^\ell_i (\eta_k) :\Res^\ell_i (\eta_k) \text{ is well-defined}\}$.
        \item For $k \leq m$, let $I_{2,k} \defeq I_1 \cup \bigcup\limits_{\ell < k} I^*_{2,k}$.
    \end{enumerate}
\end{enumerate}
Clearly
\begin{enumerate}
    \item [$\boxdot_6$] $I_1 = I_{2,0} \subseteq^* I_{2,1} \subseteq^* \ldots \subseteq^* I_{2,m} \subseteq I_2$
\end{enumerate}
As $\lessdot$ is transitive (and by part (2))
\begin{enumerate}
    \item [$\boxdot_7$] Without loss of generality $m=1$, and one of the following occurs:
    \begin{enumerate}
        \item [\textbf{A})] $I_2 \setminus I_1 = \{\eta_0\}$ and $\lh(\eta_0) < \omega$.
\sn
        \item [\textbf{B})] $I_2 \setminus I_1 \subseteq \big\{ \eta_0,\Res^\ell_n(\eta_0) : n_0 \leq n < \omega,\ \ell < h(n) \big\}$ and $\lh(\eta_0)=\omega$.
    \end{enumerate}
\end{enumerate}
Now note
\begin{enumerate}
    \item [$\boxdot_8$] If $c = c_1 \cup c_2$ with $c_1,c_2 \neq 0$ in $\ba[\bfB^*,\bar a,I_2]$ {then} it suffices to prove \ref{5.2}(5) for $c_1$ and for $c_2$.
\end{enumerate}
By $\boxdot_7$+$\boxdot_8$, noting that if $c \in \ba[\bfB^*,\bar a,I_1]$ then the conclusion is trivial, without loss of generality we have
\begin{enumerate}
    \item [$\boxdot_9$] $c = d^* \cap x_{\eta_0}$ or $c = d^* - x_{\eta_0}$, where $d^* \in \ba[\bfB^*,\bar a,I_1] \setminus \{0\}$.
\end{enumerate}

Lastly, 
\begin{enumerate}
    \item [$\boxdot_{10}$] It suffices to construct a function $\bff$ as there  such that that the homomorphism $\hat\bff$ from $\ba[\bfB^*,\bar{a},I_2]$ into $\ba[\bfB^*,\bar{a},I_1] \rest d^*$ induced by $\bff$ (which is the identity on the latter by its definition) will satisfy $\hat\bff(c) \ge d^*$.
\end{enumerate}

The proof now splits into cases (A) and (B) from $\boxdot_7$.

\medskip
In case (A), let $\bff(x_{\eta_0})$ be $d^*$ if $m_0=1$, and let
$\bff(x_{\eta_0})$ be $0$ if $m_0=0$ (and $\bff(b)=b\cap d^*$ 
if $b\in \ba [\bfB^*,\bar a,I_1]$). In case (B), if $\alpha_* = \omega$ 
then act similarly; i.e. define $\bff(x_\nu)=d^*$ for 
$\nu \in I_2 \setminus I_1$ if $m_0 = 1$, 
and $0$ if $m_0 = 0$ for $\eta\in I_2\setminus I_1$. So in
case (B), without loss of generality $\alpha_* < \omega$, and so by the assumptions of \ref{5.4}(4) we have 
\begin{enumerate}
    \item[$\boxdot_7$] For some $n_* < \omega$, we have (\textbf{A}) and (\textbf{B}), where 
    \begin{enumerate}
        \item [\textbf{A})] $\alpha_* < n_* < \omega$
\sn
        \item [\textbf{B})] If $n \geq n_*$ and $\ell < h(n)$, then $\Res_n^\ell(\eta_0) \notin I_1$.
    \end{enumerate}
\end{enumerate}
Now by repeated use of case (A),
without loss of generality\footnote{
    But we do not have to use it.
}
\[
\big(\forall n\leq\alpha_*\big) \big(\forall\ell<h(n)\big)
\big[\Res^\ell_n(\eta_0)\in I_1\big].
\]

\mn
Let
\[
\begin{array}{ll}
\bff(b)=b\cap d^*   &\text{ for }b \in \ba[\bfB^*,\bar{a},I_1],\\
\bff(x_{\eta_0})=d^*&\text{ if } m_0 = 1,\\
\bff(x_{\eta_0})=0  &\text{ if } m_0 = 0,\\
\bff(x_{\Res^\ell_n\! (\eta_0)}) = 0   & \text{ if } n \in [n_*,\omega) \text{ and } \ell < h(n) \text{ is odd,}\\
\bff(x_{\Res^\ell_n\! (\eta_0)}) = d^* & \text{ if } n \in [n_*,\omega) \text{ and } \ell < h(n) \text{ is even.}
\end{array}
\]

\mn
Now check.
\end{PROOF}

\begin{claim}\label{5.3}
Assume $h \geq 3$ or just $h(n) \geq 3$ for $n$ large enough. 
If $\bfB^*,\bar{a}, I$ are as in Definition \ref{5.1}, $I$ standard, 
$\lambda = \cf(\lambda) > \aleph_0$ and $\bfB^*$ satisfies the
[strong] $\lambda$-cc, \underline{then}  $\ba[\bfB^*,\bar{a},I]$
satisfies the [strong] $\lambda$-cc.
\end{claim}

\begin{PROOF}{\ref{5.3}}
Let $c_i\in\ba[\bfB^*,\bar a,I]$ for $i<\lambda$, $c_i \ne 0$. 
Without loss of generality  $c_i$ has the form
\[
c_i = d_i \cap \bigcap\limits_{\ell<m_{i,0}} x_{\eta_{i,\ell}} \cap
\bigcap \limits_{\ell\in [m_{i,0},m_{i,1})} \bigl(1-x_{\eta_{i,\ell}}\bigr),
\]
where $\eta_{i,\ell} \in I$, $d_i \in \bfB^* \setminus \{0\}$. By \ref{5.4}(4),
without loss of generality $d_i \leq a_{n_i}$ for some $n_i < \omega$ or 
$n_i = \omega \wedge \bigwedge\limits_{n<\omega} [d_i \cap a_n = 0]$. 
Without loss of generality $m_{i,0} = m_0$, $m_{i,1} = m_1$, 
$\lh(\eta_{i,\ell}) = n_\ell$, $n_i = n^*$, and
$\LL\eta_{i,\ell} : \ell < m_1\RR$ is without repetition for every $i$. 

Also letting $k_i < \omega$ be the minimal $k$ such that
\mn
\begin{enumerate}
    \item[$(*)$] 
    \begin{enumerate}
        \item $\lh(\eta_{i,\ell}) < \omega \Rightarrow\ \lh(\eta_{i,\ell}) \leq k$
\sn
        \item $n^* < \omega \Rightarrow  n^* < k$
\sn        
        \item $\ell_1 < \ell_2 < m_1 \Rightarrow \eta_{i,\ell_1} \rest k \ne 
        \eta_{i,\ell_2}\rest k$
\sn
        \item $(\forall n)[n \geq k \Rightarrow h(n)\geq 3]$
    \end{enumerate}
\end{enumerate}
\mn
and without loss of generality $k_i = k^*$; if $\lh(\eta_{i,\ell}) = \omega$,
$k < k^*$, $\ell < h(k)$ then $\Res^\ell_k(\eta_{i,\ell}) \in \{\eta_{i,m} : m < m_{i,1}\}$.

By the $\Delta$-system argument, without loss of generality 
\mn 
\begin{enumerate}
    \item[$(*)$]  If $i \ne j < \lambda$, $k \leq k^*+1$, and $m',m'' < m_1$, $\ell',\ell'' < h(k)$ and $\Res^{\ell'}_k(\eta_{i,m'}) = \Res^{\ell''}_k(\eta_{j,m'',j})$, \underline{then} for every 
    $\alpha,\beta < \lambda$ we have
    \[
        \Res^{\ell'}_k(\eta_{\alpha,m'}) = \Res^{\ell''}_k(\eta_{\alpha,m''}) =
        \Res^{\ell'}_k(\eta_{\beta,m'}) = \Res^{\ell''}_k(\eta_{\beta,m''}).
    \]
\end{enumerate}
\mn
We can now check (similarly to the proof of \ref{2.6}).
\end{PROOF}

\begin{theorem}\label{5.4}
Let $\aleph_1 \leq \mu < 2^{\aleph_0}$. There is a Boolean Algebra $\bfB$ such
that:
\mn
\begin{enumerate}[$(A)$]
    \item  $\bfB$ has cardinality $2^{\aleph_0}$ and satisfies the ccc 
    (and even the strong $\lambda$-cc if $\lambda = \cf(\lambda) > \aleph_0$).
\sn
    \item  $\bfB$ has length $\mu$ (i.e. there is in $\bfB$ a chain of 
    length $\mu$ but no chain of length $\mu^+$).
\end{enumerate}
\mn
Moreover: 
\mn
\begin{enumerate}
    \item[$(B)^+$]  If $n,m < \omega$ and $\bar c^\zeta \in {}^m\bfB$ for 
    $\zeta < \mu^+$\underline{then}  for some $Y\in [\mu^+]^{\mu^+}$ (i.e. 
    $Y \subseteq \mu^+$ of cardinality $\mu^+$), the sequence 
    $\LL\bar c^\zeta : \zeta \in Y\RR$ is a $(\qf,n)$-indiscernible set in the Boolean algebra $\bfB$ (see \ref{5.4A}(2) below). 
\sn
    \item[$(C)$]  Every infinite homomorphic image of $\bfB$ has cardinality $2^{\aleph_0}$.
\end{enumerate}
\end{theorem}

\begin{remark}\label{5.4A}
0) Recall that the \emph{length} of a BA refers to the size of a linearly ordered subset, not necessarily well-ordered. 

\sn
1) Note that $(B)^+\, \Rightarrow\, (B)$; for it $m=1$ suffices, 
for this constant $h$ is OK below, but the proof here is simpler.

\sn
2) Let $ \bar \bfc = \LL{\bar c}^\zeta : \zeta \in Y\RR$ be a sequence of 
$m$-tuples from a model $M$ (for example, a Boolean algebra) and $\Delta$ 
a set of formulas in $\bbL(\tau_M)$. We say $\bar\bfc$ is an
$(\Delta,n)$-indiscernible set \underline{iff} for any $\zeta_0,\ldots,
\zeta_{n-1}$ from $Y$ with no repetitions and $\xi_0,\ldots,\xi_{n-1}$ from
$Y$ with no repetitions, the $\Delta$-type of 
$\bar c^{\zeta_0} \caret \ldots \caret \bar c^{\zeta_{n-1}}$ in $M$ is equal 
to the $\Delta$-type of $\bar c^{\xi_0}\caret \ldots \caret \bar c^{\xi_{n-1}}$ 
in $M$. For $\Delta$ the set of quantifier free formulas we write $\qf$.
\end{remark}

\begin{PROOF}{\ref{5.4}}
Let $h : \omega \to \omega$ be, for example, $h(n) = 2n+2$.

Let $I_\beta \in K^\omega_{\tr(h)}$ be standard for $\beta < 2^{\aleph_0}$, 
have cardinality continuum, and be such that:
\mn
\begin{enumerate}
    \item[$(*)_{I_\beta}$]  For every $f : I_\beta \to \theta$,
    $\theta < 2^{\aleph_0}$, for some $\eta \in P^{I_\beta}_\omega$, 
    for every $n < \omega$,
    \[
        \big(\forall\ell < h(n)\big) \big[f(\Res^0_n(\eta)) = 
        f(\Res^\ell_n(\eta))\big]
    \]
    (i.e. $\eta(m) = \big\LL \alpha_\ell : \ell < h(n) \big\RR \Rightarrow 
    \big|\big\{f(\eta \rest n \caret \LL\alpha_\ell\RR) : \ell < h(n)\big\}\big| = 1$.)
\end{enumerate}
\mn
[Why do such $I$-s exist? The full tree will serve; that is, we let
\begin{align*}
I_\beta = \big\{\LL\bar{\alpha}^\ell : \ell < \gamma\RR : &\ \gamma\leq\omega,\
\bar{\alpha}^\ell\text{ an increasing sequence of length }h(\ell)\\
&\text{ from }2^{\aleph_0},\text{ except in the case }0 < \gamma < \omega \wedge
\ell = \gamma-1;\\
&\text{ then we demand }\bar\alpha^\ell\text{ is just an ordinal} <
2^{\aleph_0}\big\}.
\end{align*}

This is as required, as for any $f : I_\beta \to \theta$ we can choose a 
sequence $\bar\alpha_\ell = \LL\beta_{\ell,0},\ldots,\beta_{\ell,h(\ell)-1}\RR$
by induction on $\ell < \omega$, where 
$\beta_{\ell,0} < \ldots < \beta_{\ell,h(\ell)-1} < 2^{\aleph_0}$ 
and $f(\LL\bar\alpha^0,\ldots,\bar\alpha^{\ell-1},\beta_{\ell,i}\RR)$ 
does not depend on $i < h(\ell)$. This is possible as $2^{\aleph_0} > |\rang(f)|$.  
So $I_\beta$-s as required in $(*)_{I_\beta}$ indeed exist.]

\medskip
We shall now construct Boolean algebras $\bfB_\alpha$ (for $\alpha \leq 2^{\aleph_0}$) 
and $\bar a^\alpha = \LL a^\alpha_n : n < \omega \RR$ such that:
\sn
\begin{enumerate}
    \item[(I)]  
    \begin{enumerate}
        \item $\bfB_0$ is a subalgebra of $\clP(\omega)$ of cardinality 
        $\mu$ with a chain of cardinality $\mu$ satisfying the ccc 
        (even the strong $\lambda$-cc, when $\lambda = \cf(\lambda) > \aleph_0$).

\sn
        [E.g.~let $A$ be a set of $\mu$ reals, let $h$ be a one-to-one function from $\omega$ onto the rationals and $\bfB$ is the Boolean algebra of subsets of $\omega$ generated by $\big\{\{n : h(n) < a\} : a \in A\big\}$.  Clearly $\gB$ has a linearly ordered subset of cardinality $\mu$ (e.g.~its set of generators). Of course, its length is not $>\mu$ as its cardinality is $\mu$. Lastly, it satisfies the ccc because the set of nonempty rational intervals is dense in it.]
\sn
        \item $\bfB_\alpha$ is increasing continuous, of cardinality 
        $2^{\aleph_0}$ if $\alpha > 0$.
\sn
        \item $\bar a^\alpha$ is an $\omega$-sequence of pairwise 
        disjoint non-zero elements of $\bfB_\alpha$.
\sn
        \item If $\alpha < 2^{\aleph_0}$, $a_n \in \bfB_\alpha \setminus 
        \{\mathbf{0}_{\bfB_\alpha} \}$, and $\bigwedge\limits_{n\neq m} [a_n \cap a_m = 0]$ then for $2^{\aleph_0}$ many ordinals $\beta$, we have $\bigwedge\limits_{n<\omega} [a_n = a^\beta_n]$.

        [You can demand that $\{a^\alpha_n : n < \omega\}$ is a maximal antichain; it does not matter.]
\sn
        \item $\bfB_{\alpha+1} = \ba[\bfB_\alpha,\bar a^\alpha,I_\alpha]$ 
        (We denote the $x_\eta$ by $x^\alpha_\eta$ for $\eta\in I_\alpha$.)
    \end{enumerate}
\end{enumerate}
\mn
There is no problem to do the bookkeeping, and 
$\bfB_\alpha \subseteq \bfB_{\alpha +1}$ by \ref{5.2}(1). We shall show 
that $\bfB \defeq \bfB_{2^{\aleph_0}}$ is as required. Obviously $\bfB$ 
has cardinality $2^{\aleph_0}$.

By \ref{5.2}(3) clearly $\bfB_\alpha \lessdot \bfB_{\alpha+1}$, 
so we can prove by induction on $\alpha$ that 
$$\beta < \alpha\, \Rightarrow\, \bfB_\beta \lessdot \bfB_\alpha,$$
by \ref{2.7} and \ref{2.8}. We can also prove by induction on $\alpha$
that $\bfB_\alpha$ satisfies the $\aleph_1$-cc
(even the strong $\lambda$-cc when $\lambda = \cf(\lambda) > \aleph_0$): 
the successor stage is proved by \ref{5.3}, the limits steps by \ref{2.8}. 
So demand (A) from \ref{5.4} holds. If $\bff$ is a homomorphism from
$\bfB$ onto some $\bfB'$ with $\aleph_0 \leq \|\bfB'\| < 2^{\aleph_0}$ 
then there are $b_n \in \bfB'\setminus \{0\}$ pairwise disjoint. 
Now for some $a_n\in \bfB$, $\bff(a_n)= b_n$ and without loss of generality 
$\bigwedge\limits_{n\neq m} [a_n \cap a_m = 0]$ (otherwise use 
$a'_n = a_n \setminus \bigcup\limits_{m<n} a_m$).
Hence for every infinite co-infinite $Y \subseteq \omega$. for some
$\alpha=\alpha_Y$:
\[
\{a^\alpha_{2n} : n < \omega\} = \{a_n : n \in Y\}\quad\text{and}\quad\{
a^\alpha_{2n+1} : n < \omega\} = \{a_n : n \in \omega \setminus Y\}.
\]

\mn
Now define 
$g : I_\alpha \to \bfB'$ by $g(\eta) = \bff(x^\alpha_\eta)$,
so by the choice of the $I_\alpha$-s (i.e. by $(*)_{I_\beta}$) 
for some $\eta = \eta_Y \in P^{I_\alpha}_\omega$, for every $n$, 
letting $\eta(n) = \LL\alpha_0,\alpha_i,\ldots,\alpha_{h(n)-1}\RR$, we have
\[
\bigwedge\limits_{\ell<h(n)} \big[ \bff(x^\alpha_{\eta\rest n \caret
\LL\alpha_\ell\RR}) = \bff(x^\alpha_{\eta\rest n \caret \LL\alpha_0\RR}) \big].
\]

\sn
Hence $\bff \bigl(x^\alpha_{\eta \rest n \caret \LL\alpha_\ell\RR} - x^\alpha_{\eta\rest n \caret\LL\alpha_{\ell+1}\RR}\bigr) = \mathbf{0}_{\bfB'}$ 
for $\ell < h(n)-1$ and hence
\[
\bff\bigl(a^\alpha_n \cap \bigcap\limits_{\ell<\frac{h(n)-1}{2}} \bigl(1 -
(x^\alpha_{\eta\rest n \caret \LL\alpha_{2\ell}\RR} -
x^\alpha_{\eta \rest n \caret \LL\alpha_{2\ell+1}\RR}) \bigr)\bigr) = 
\bff(a^\alpha_n)
\]

\sn
and  
\[
\bff\bigl(a^\alpha_n-\bigcup\limits_{\ell<\frac{h(n)-1}{2}}
\bigl(x^\alpha_{\eta \rest n \caret 
\LL\alpha_{2\ell}\RR}-x^\alpha_{\eta\rest n \caret
\LL\alpha_{2\ell+1}\RR}\bigr)\bigr)=\bff(a^\alpha_n).
\]

\sn
Hence (see Definition \ref{5.1})
\[
\begin{array}{rcl}
n\text{ is even}&\Rightarrow& \bfB' \models 
\bff(a^\alpha_n)\leq \bff(x^\alpha_\eta),\\
n\text{ is odd} &\Rightarrow& \bfB' \models \bff(a^\alpha_n)\cap 
\bff(x^\alpha_\eta)=0.
\end{array}
\]

\mn
Therefore,
\[
\begin{array}{rcccl}
m\in Y&\Rightarrow&\text{ for some even $n$, }a^\alpha_n=a_m&\Rightarrow&
\bfB' \models b_m \leq \bff(x^\alpha_\eta),\\
m\in\omega\setminus Y&\Rightarrow&\text{ for some odd $n$, }a^\alpha_n=a_m&
\Rightarrow& \bfB' \models b_m\cap \bff(x^\alpha_\eta)=0.
\end{array}
\]

\mn
As this occurs for every infinite co-infinite $Y \subseteq \omega$, 
for some $\alpha = \alpha(Y)$ and $\eta = \eta_Y$, clearly the 
$\bff(x_{\eta_\alpha}^{\alpha(Y)})$ determine $Y$ using $\bfB'$ 
{and} $\LL a_n : n < \omega\RR$.
So clearly we get $2^{\aleph_0}$-many distinct members of $\bfB'$ 
(simply put, the $\bff(x_{\eta_Y})$), 
a contradiction. So demand \textbf{(C)} of {\ref{5.4}} holds.

What about the length, i.e, clauses \textbf{(B)} and \textbf{(B)}$^+$? For \textbf{(B)}, first note that $\bfB_0$ has a chain of cardinality $\mu$ and hence so
does $\bfB$. If $J\subseteq \bfB$ is a chain, $|J|=\mu^+$, then 
$(\bfB)^+$ gives a contradiction and even the ``weakly indiscernible sequence'' 
version does because as $\bfB\models$ ccc, it has 
no subset of order type $\mu^+$ or $(\mu^+)^*$; but 
the variant of $(\bfB)^+$ implies just this ($m=1$ suffices).

So it suffices to prove that clause $\textbf{(B)}^+$ holds for $\bfB_\alpha$ by
induction on $\alpha$.

\mn
\textbf{Case 1}:  $\alpha=0$.

Trivial (we can get $\bar c^\zeta$ constant on some $Y \in [\mu^+]^{\mu^+}$).

\mn
\textbf{Case 2}:  $\alpha$ is limit, $\cf(\alpha) \ne \mu^+$.

For some $\beta < \alpha$,
\[
Y_1 = \{\zeta < \mu^+ : \bar c^\zeta \subseteq \bfB_\beta\} \in
[\mu^+]^{\mu^+},
\]

\mn
(note that if $\cf(\alpha) < \mu^+$, then we can get $Y_1 = \mu^+$)
and use the induction hypothesis.

\mn
\textbf{Case 3}: $\cf(\alpha) = \mu^+$.

Let $\LL\beta_\eps : \eps < \mu^+\RR$ be an increasing
continuous sequence with limit $\alpha$. Let $n,m,\LL\bar c^\zeta : \zeta < \mu^+\RR$ be given. Without loss of generality  $\bar c^\zeta = \LL
c^\zeta_\ell : \ell < m\RR$ is a partition of $\mathbf{1}_{\bfB_\alpha}$ 
(i.e.~$\ell_1\neq\ell_2 \Rightarrow c^\zeta_{\ell_1} \cap 
c^\zeta_{\ell_2}=0$ and $\mathbf{1}_{\bfB_\alpha}=
\bigcup\limits_{\ell<m}c^\zeta_\ell$). For each $\zeta < \mu^+$, we can find
$a^\zeta_\ell,b^\zeta_\ell \in \bfB_{\beta_\zeta}$ (for $\ell < m$) such that:
\mn
\begin{enumerate}[(a)]
    \item  $a^\zeta_\ell \le c^\zeta_\ell\leq b^\zeta_\ell$
\sn
    \item  $(0 < x \le b^\zeta_\ell - a^\zeta_\ell) \wedge x \in \bfB_{\beta_\zeta} \Rightarrow  (x\cap c^\zeta_\ell - a^\zeta_\ell \ne 0) \wedge 
    (x - c_\ell^\zeta \neq 0)$.
\end{enumerate}
\mn
[Why? By use of \ref{5.2}(5) applied to $c_\ell^\zeta$ and to $\mathbf{1}_{\bfB_{\beta_\zeta}} - c_\ell^\zeta$).  

If $\zeta$ is limit then for some $f(\zeta) < \zeta$ we have 
$\{a^\zeta_\ell,b^\zeta_\ell : \ell < m\} \subseteq \bfB_{f(\beta_\zeta)}$. 
By Fodor's lemma, for some $\eps(*) < \mu^+$ and a stationary set 
$S \subseteq \mu^+$, we have $\bigwedge\limits_{\zeta\in S} [f(\zeta) = \eps(*)]$. 

So
\mn
\begin{enumerate}
    \item[(c)]  $\eps\in S \Rightarrow \{a^\eps_\ell,b^\eps_\ell : \ell < m\} \subseteq \bfB_{\beta_{\dot\eps(*)}}$.
\end{enumerate}
\mn
Also without loss of generality 
\mn
\begin{enumerate}
    \item[(d)]  If $\eps < \zeta \in S$ then $\{c^\eps_\ell : \ell < m\} \subseteq \bfB_{\beta_\zeta}$.
\end{enumerate}
\mn
Now apply the induction hypothesis on $\bfB_{\beta_{\dot\eps(*)}}$ and 
$\LL \bar a^\zeta \caret \bar b^\zeta :\zeta \in S\RR$,
where\\ $\bar a^\zeta = \LL a^\zeta_\ell : \ell < m\RR$ and
$\bar b^\zeta \defeq  \LL b^\zeta_\ell : \ell < m\RR$.

So there is $Y \in [S]^{\mu^+}$ such that $\LL \bar b^\zeta:\zeta
\in Y\RR$ is an $(n,\qf)$-indiscernible set.  So let $\zeta_0 <
\ldots < \zeta_{n-1}$ be from $S$ and for $k \le n$ let $\bfB'_k$ be
the subalgebra of $\bfB$ generated by $X_k \defeq \bfB_{\beta_{\eps(*)}} \cup \{ c^{\zeta_i}_\ell : i < k,\ \ell < m\}$.
We now prove, by induction on $k \leq n$, that
\begin{enumerate}
    \item[(e)$_n$] $\bfB_k''$ is freely generated by $X_k$, except for:
    \begin{enumerate}[$\bullet_1$]
        \item The equations satisfied by $\bfB_{\beta_{\eps(*)}}$.
\sn
        \item $a^\zeta_\ell \le c^\zeta_\ell\leq b^\zeta_\ell$ for $\ell < m$ and $\zeta \in \{\zeta_i : i < k\}$.
\sn
        \item $\LL c_\ell^\zeta : \ell < \omega\RR$ is a partition of $\mathbf{1}$.
    \end{enumerate}
\end{enumerate}
For $k = 0$ this is trivial. For $k+1$ we use clauses (c) and (d) above. Lastly, for $k = n$ we get the desired conclusion.

\mn
\textbf{Case 4}:  $\alpha = \beta+1$.

Let $n,m<\omega$ and $\bar c^\zeta\in {}^m(\bfB_{\beta+1})$ for $\zeta<\mu^+$
be given, $\bar c^\zeta=\LL c^\zeta_\ell:\ell<m\RR$. So there are
$k_{\zeta,0} = k(\zeta,0) < \omega$, $k_{\zeta,1} = k(\zeta,1) < \omega$ and $b^\zeta_0,\ldots,b^\zeta_{k_{\zeta,0}-1} \in \bfB_\beta$, 
$\eta^\zeta_0,\ldots,\eta^\zeta_{k_{\zeta,1}-1} \in I_\beta$ 
and Boolean terms $\sigma^\zeta_\ell$ (for $\ell < m)$ such that
\[
c^\zeta_\ell = \sigma^\zeta_\ell\bigl(b^\zeta_0,\ldots,b^\zeta_{k(\zeta,0)
-1},x_{\eta^\zeta_0},\ldots,x_{\eta^\zeta_{k(\zeta,1)-1}}\bigr).
\]

\mn
Without loss of generality $\LL\eta^\zeta_\ell : \ell < k_{\zeta,1}\RR$
is a $\Delta$-system.

Without loss of generality   $k_{\zeta,0}=k_0$, $k_{\zeta,1}=k_1$,
$\sigma^\zeta_\ell=\sigma_\ell$ and $\lh(\eta^\zeta_\ell) =
m_\ell\leq \omega$ for every $\zeta < \mu^+$.   

Also, there is $k_{\zeta,2}<\omega$ such that:
\mn
\begin{enumerate}
    \item[(II)] 
    \begin{enumerate}
        \item[$(\alpha)$] $\lh(\eta^\zeta_\ell) < \omega \Rightarrow \lh(\eta^\zeta_\ell) < k_{\zeta,2}$
\sn
        \item[$(\beta)$] $\eta^\zeta_{\ell_1} \neq \eta^\zeta_{\ell_2} 
        \Rightarrow \eta^\zeta_{\ell_1} \rest k_{\zeta,2} \neq 
        \eta^\zeta_{\ell_2} \rest k_{\zeta,2}$
\sn
        \item[$(\gamma)$] $2n+2 < k_{\zeta,2}$.
    \end{enumerate}
\end{enumerate}
\mn
Without loss of generality  $\bigwedge\limits_\zeta k_{\zeta,2}=k_2$.

Without loss of generality  the statement $(*)$ (with $k_2$ here for $k^*$
there and is $>n$) from the proof of \ref{5.3} holds (essentially being a
$\Delta$-system), i.e.
\mn
\begin{enumerate}
    \item[$(*)$]  If $i \neq j < \lambda$, $k \leq k_2+1$, and $m',m'' < m_1$, $\ell',\ell'' < h(k)$ and $\Res^{\ell'}_k(\eta_{m',i}) = \Res^{\ell''}_k(\eta_{m'',j})$, then for every $\alpha,\beta < \lambda$ we have:
    \[
        \Res^{\ell'}_k(\eta_{m',\alpha}) = \Res^{\ell''}_k(\eta_{m'',\alpha}) =
        \Res^{\ell'}_k(\eta_{m',\beta}) = \Res^{\ell''}_k(\eta_{m'',\beta}).
    \]
\end{enumerate}
\mn
Let $\bar b^\zeta = \LL b^\zeta_\ell : \ell < k_0\RR$. By the induction
hypothesis, without loss of generality $\big\LL\bar b^\zeta \caret 
\LL a^\beta_\ell : \ell \leq k_2\RR : \zeta < \mu^+ \big\RR$ is 
$(\qf,n)$-indiscernible and without loss of generality the sequence 
$\big\LL \LL\eta^\zeta_\ell \rest (k_2+1) : \ell < k_1\RR : \zeta < \mu^+\RR$ 
is indiscernible (sequence of finite sequences of ordinals).

\mn
To finish the proof of \ref{5.4} it suffices to observe \ref{5.4F} below.
\end{PROOF}

\begin{observation}\label{5.4F}
If $\bfB^* = \ba[\bfB,\bar{a},I]$, $n^* < \omega$, 
$I^0 = \{\eta\in I : \lh(\eta) \leq n^*\}$, $Z \subseteq P^I_\omega$, 
and for every $\nu \in Z$ and $n \geq n^*$ the set 
$$\{\nu' \rest (n+1) :  \nu' \in Z,\ \nu' \rest n = \nu \rest n\}$$ 
has $\leq \lfloor \frac{h(n)-1}{4}\rfloor$ elements, \underline{then} 
$\{x_\eta : \eta \in Z\}$ is independent in $\bfB^*$ over 
$\bfB^0 \defeq \ba[\bfB,\bar{a},I^0]$, except the equations 
$c^+_\eta \leq x_\eta \wedge c^-_\eta \cap x_\eta = 0$ for $\eta \in Z$, where
\begin{align*}
    c^+_\eta \defeq \bigcup_{2n < n^*} &\Big( a_{2n} - \bigcup \big\{ x_{\Res^{2\ell}_{2n}(\eta)} - x_{\Res^{2\ell+1}_{2n}(\eta)} : \ell < h(n)/2 \big\} \Big) \\
    c^-_\eta \defeq \bigcup_{2n+1<n^*} &\Big( a_{2n+1} - \bigcup \big\{ 
    x_{\Res^{2\ell}_{2n+1}(\eta)} - x_{\Res^{2\ell+1}_{2n+1}(\eta)} : 2\ell+1 < h(2n+1) \big\} \Big) 
\end{align*}
(Note: $\eta_1 \rest n^* = \eta_2 \rest n^*\ \Rightarrow\ 
(c^+_{\eta_1},c^-_{\eta_1}) = (c^+_{\eta_2},c^-_{\eta_2})$.)
\end{observation}

\begin{PROOF}{\ref{5.4F}}
Let $\bff_0$ be any function with domain $X = \{x_\eta : \eta \in Z\}$ 
such that 
$$
\bff_0(x_\eta) \in \big\{c \in \ba[\bfB,\bar a,I] : c^+_\eta \leq c \leq \mathbf{1}_\bfB - c^-_\eta \big\},
$$ 
and let
\[
J^1 = I^0 \cup X \cup\{\Res^\ell_n(\nu) : \nu \in Z,\ \ell < h(n),\ n < \omega\}.
\]

\mn
Clearly, by \ref{5.2}(2) it suffices to find a homomorphism from
$\bfB_1 \defeq \ba[\bfB,\bar{a},J^1]$ into $\bfB^0$ extending 
$\id_{\bfB^0} \cup \bff_0$. For this it suffices to find a mapping $\bff$ 
from $\bfB \cup \{x_\eta : \eta \in J^1\}$ into $\bfB^0$ extending 
$\id_{\bfB^0},\bff_0$, and $\id_{\{x_\eta : \eta \in I^0\}}$, and preserving 
the equations defining $\ba[\bfB,\bar{a},J^1]$. As $\bff \rest \bfB$,
$\bff \rest \{x_\eta\in I : \lh(\eta) \leq n^*\}$, and 
$\bff \rest \{x_{(\eta,\varrho)} : \eta \in Z\}$ are defined, and
\[
J^1 = \bigcup\limits_{n\in [n^*,\omega)} \{x_\eta : \eta \in Z_n\} \cup X \cup Z
\]
where $Z_n = \{\eta \in J^1 : \lh(\eta) = n+1\}$, it will suffice to choose 
$\bff \rest \{x_\eta : \eta \in Z_n\}$ for each $n \in [n^*,\omega)$ to finish 
the definition of $\bff$. 

Let $Y_n = \{\nu \rest n : \nu \in Z_n\}$, and for $\eta \in Y_n$ let 
$X_{n,\eta} = \{\nu \in Z_n : \nu \rest n = \eta\}$. Clearly $\LL X_{n,\eta} : \eta \in Y_n \RR$ is a partition of $Z_n$.

For $\eta \in Y_n$ let 
$\cP_{n,\eta} = \{\nu \rest(n+1) : \nu \in Z,\ \nu \rest n = \eta\}$ and 
$$
\clS_{n,\eta} = \big\{(\rho\rest n) \caret \LL\rho(n)(\ell) \RR : \ell < h(n)\big\}.
$$

By the assumption on $Z$, for every $\eta \in Y_n$ 
the set $\cP_{n,\eta}$ has $< \frac{h(n)-1}{4}$ elements. Now:
\begin{enumerate}
    \item[$(*)$] For $\eta \in Y_n$ there is a function 
    $f_\eta : \clS_{n,\eta} \to \{\mathbf{0}_{\bfB^*},\mathbf{1}_{\bfB^*}\}$ 
    such that if $\nu \in \cP_{n,\eta}$ is equal to 
    $\eta \caret \LL \alpha_0,\ldots,\alpha_{h(n)-1}\RR$ 
    then for some $\ell < (h(n)-1) / 2$ we have 
    $f_\eta(x_{\eta\caret\LL \alpha_{2\ell}\RR}) = \mathbf{1}_{\bfB^*}$ and 
    $f_\eta(x_{\eta\caret\LL \alpha_{2\ell+1}\RR}) = \mathbf{0}_{\bfB^*}$.
\end{enumerate}
[Why is this possible? Let 
$$
\big\LL \nu_k = \eta \caret \LL \bar \alpha^k\RR : k < |\clS_{n,\eta}| \big\RR
$$ 
list $\clS_{n,\eta}$, where $\bar \alpha^k = \LL \alpha_0^k, \ldots, \alpha_{h(n)-1}^k \RR$. By induction on $k < |\clS_{n,\eta}|$, we can choose $\ell(k) < \frac{h(n)-1}{2}$ such that 
$$
\alpha_{2\ell(k)}^k, \alpha_{2\ell(k)+1}^k \notin \{ \alpha_{2\ell(i)}^i, \alpha_{2\ell(i)+1}^i : i < k\}.
$$
Upon arriving to $k$, $\alpha_{2\ell(i)}^i$ will disqualify at most one candidate for $\ell(k)$, and $\alpha_{2\ell(i)+1}^i$ will disqualify at most one more. So at most, $2k$ candidates are disqualified. As $k < \frac{h(n)-1}{4}$, not all $\ell < \frac{h(n)-1}{2}$ have been disqualified.]

\smallskip
Now define $\bff \rest Z_n$ as follows: if $\nu \in Z_n$ then $\nu \in \cP_{n,\eta}$ for some $\eta \in Y_n$, and so we let $\bff(x_\nu) = f_\eta(x_\nu)$.

\mn
Now check. 
\end{PROOF}

\begin{discussion}\label{5.16}
1) In the proof of clause $(B)^+$, the successor case we use the
fact that $h(n)$ converges to $\infty$, as when the level increases we need
more $\eta\in P^{I_\beta}_\omega$ to see non-freeness.

\sn
2) The proof there for limit $\alpha$ uses just 
``$\LL \bfB_i : i \leq 2^{\aleph_0}\RR$ is $\lessdot$-increasing continuous with 
projections'' (i.e. \ref{5.2}(5)), and the induction hypothesis.

\sn
3)  We can vary the construction in some ways. We can demand that each
$\bar{a}^\alpha$ is a maximal antichain --- no difference so far. We may like
to use
$\LL I_\beta:\beta<2^{\aleph_0}\RR$ such that $I_\beta$ is not super
unembeddable into $\sum\limits_{\gamma\neq\beta} I_\gamma$. 
We can construct our Boolean algebra to be monorigid (i.e.~with no one-to-one endomorphism), and even get $2^{2^{\aleph_0}}$ such
Boolean algebras, no one embeddable to another: even restricting to 
non-zero elements, even not embeddable into the completion of another. 
To carry this out we need the following for $\lambda = 2^{\aleph_0}$: 
there is $\bar I = \LL I_\alpha : \alpha < \lambda\RR$ exemplifying that 
$K^\omega_{\tr(h)}$ has the full $(\lambda,\lambda,\aleph_1,\aleph_1)$-super 
bigness property, such that for at least one $\beta$,
$I_\beta$ satisfies $(*)_{I_\beta}$ from the beginning of the proof of
\ref{5.4}. Now such a $\bar I$ does exist (with $(*)_{I_\beta}$ for every
$\beta$); this may be elaborated elsewhere.

\sn
4)  Of course the proof works for $\mu= 2^{\aleph_0}$.

\sn
5)  We can separate some parts of the proof to independent claims. We can
ask for ``$\bfB$ has length $\mu$, but no chain of cardinality $\mu$'' (i.e.
the supremum is not obtained) for $\mu$ limit. It is natural to demand
$\cf(\mu) >\aleph_0$. Next, we address this.
\end{discussion}

\begin{claim}\label{5.17}
$1)$ Assume $2^{\aleph_0} \ge \mu$ and $\aleph_0 < \kappa = \cf(\mu) < \mu$.
\underline{Then} there is a Boolean algebra $\bfB$ such that $|\bfB| = 2^{\aleph_0}$,
$\bfB$ has no homomorphic image of cardinality $\in [\aleph_0,2^{\aleph_0})$,
and $\mathrm{length(\bfB)} = \mu$, but the supremum is not obtained (i.e.~$\mathrm{length}^+(\bfB) = \mu$).

\sn
$2)$ Similarly, but slightly modifying the assumption to $\aleph_0 < \kappa = \cf(\mu) = \mu$.
\end{claim}

\begin{PROOF}{\ref{5.17}}
Like \ref{5.4}.

\sn
1) Let $\mu = \sum\limits_{i<\kappa} \mu_i$ with $\LL\mu_i : i < \kappa\RR$
be increasing continuous and $\kappa < \mu_i < \mu$. For $\eps < \kappa$, 
let $\bfB^\eps$ be a subalgebra of $\clP(\omega)$ of cardinality
$\mu_\eps$ and length $\mu_\eps$. We stipulate that $\bfB_\eps$ satisfies the ccc (and moreover, the strong $\lambda$-cc for $\lambda = \cf(\lambda) > \aleph_0$).

\sn
[Why does this exist? As in (I)(a) in the proof of Theorem \ref{5.4}.]

\smallskip
Let $\LL I_\alpha : \alpha < 2^{\aleph_0} \times \kappa\RR$ 
be as in the proof of \ref{5.4}. We define $\bfB_\alpha$ 
(for $\alpha \leq 2^{\aleph_0} \times \kappa$) similarly 
to the proof of \ref{5.4}. Specifically:
\[
\begin{array}{l}
\bfB_0 = \bfB^0, \text{ the trivial Boolean algebra,}\\
\bfB_{2^{\aleph_0}\times\eps+1}\text{ is the free product
}\bfB_{2^{\aleph_0}
\times\eps}* \bfB^\eps,\\
\bfB_\alpha\text{ is increasing continuous in }\alpha,\\
\bfB_{\alpha+1}=\ba[\bfB_\alpha,\bar a^\alpha,I_\alpha]\text{ for }\alpha<
2^{\aleph_0}\times\kappa,\ \alpha\notin\bigl\{2^{\aleph_0}\times\eps:
\eps<\kappa\bigr\},
\end{array}
\]
\mn
where $\bar a_\alpha$ is $\LL a_{\alpha,n} : n < \omega\RR$, 
$a_{\alpha,n} \in \bfB_\alpha,a_{\alpha,n} > 0$, 
$n_1 \ne n_2 \Rightarrow a_{\alpha,n_1} \cap a_{\alpha,n_2} = 0$.  The choice
of the $\bar a_\alpha$-s (i.e. the bookkeeping) is as in the proof 
of \ref{5.4} above.

So, by the proof of \ref{5.4}:
\mn
\begin{enumerate}
\item[$(*)$]  If $\alpha \leq 2^{\aleph_0} \times \kappa$ \underline{then}
\sn 
    \begin{enumerate} 
        \item[$(\alpha)$]  $\bfB_\alpha$ is $\lessdot$-increasing continuous, and satisfies the strong $\lambda$-cc if $\lambda = \cf(\lambda) > \aleph_0$.
\sn
        \item[$(\beta)$]  $\bfB_{1+\alpha}$ has cardinality $2^{\aleph_0}$ and length 
        $\mu_0 + \sum\limits_{2^{\aleph_0} \times \eps < \alpha} \mu_\eps$, which is $< \mu$ when $\alpha < 2^{\aleph_0} \times \kappa$.
\sn
        \item[$(\gamma)$] If $\alpha=2^{\aleph_0} \times \eps$ with $\eps$ a successor ordinal, then $\bfB_\alpha$ has no homomorphic image of cardinality $\in [\aleph_0,2^{\aleph_0})$. 
\sn
        \item[$(\delta)$]  If $\alpha < \beta \leq 2^{\aleph_0} \times \kappa$ 
        and $b \in \bfB_\beta \setminus \{\mathbf{0}_{\bfB_\beta}\}$ then for some 
        $a \in \bfB_\alpha$ we have $\bfB_\beta \models b \leq a$ and if 
        $\bfB_\alpha\models 0 < a' \le a$ then $a' \cap b > \mathbf{0}_{\bfB_\beta}$.
    \end{enumerate}
\end{enumerate}
\mn
[Note: for clause $(\beta)$ we use the proof of 
(B)$^+$ of \ref{5.4}. For $\alpha=2^{\aleph_0}\times\eps+1$ for 
clause $(\delta)$ we have a new clause, but easy one].

It follows that
\mn
\begin{enumerate}
    \item[$(**)$] $\bfB=\bfB_{2^{\aleph_0}\times\kappa}$:
    \begin{enumerate}
        \item Is a ccc Boolean algebra of cardinality continuum.
\sn
        \item Has length $\mu$.
\sn
        \item Has no infinite homomorphic image of cardinality $<2^{\aleph_0}$.
    \end{enumerate}
\end{enumerate}
[Why? $\bfB$ is a ccc Boolean algebra by $(*)(\alpha)$ and is cardinality continuum by $(*)(\beta)$, so clause (a) holds. The length of $\bfB$ is $\mu$ 
by {$(*)(\beta)$}. Lastly, if $f$ is a homomorphism from $\bfB$
onto some infinite $\bfB'$, then (as $\kappa = \cf(\kappa)> \aleph_0$) for 
some $\eps < \kappa$ we know $\rang(f \rest \bfB_{2^{\aleph_0} \times (\eps+1)})$ is infinite, and hence has cardinality continuum. Therefore $\bfB'$ has cardinality continuum, and so clause (c) holds as well.]

\mn
Now to finish, we just need to show
\mn
\begin{enumerate}
    \item[$(***)$]  For $J \subseteq \bfB$ a chain of cardinality $\mu$, we get a contradiction.
\end{enumerate}
\mn
Let $\bfB^*_\eps=\bfB_{2^{\aleph_0}\times\eps}$. Let $c_\alpha \in J$ 
(for $\alpha < \mu$) be pairwise distinct.

By clause $(*)(\delta)$, for each $\eps < \kappa$ and 
$\alpha < \mu^+_\eps$ we can find $b^\eps_\alpha \in \bfB^*_\eps$ such that:
\mn
\begin{enumerate}
    \item[(a)]  $c_\alpha\leq b^\eps_\alpha$ in $\bfB$.
\sn
    \item[(b)]  $0 < x \leq b^\eps_\alpha \wedge x\in \bfB^*_\eps\ \Rightarrow\ x\cap c_\alpha\neq 0$
\end{enumerate}
\mn
Note:
\mn
\begin{enumerate}
    \item[(c)]  $b^\eps_\alpha$ is unique, and
\sn
    \item[(d)]  $c_\alpha\leq c_\beta\ \Rightarrow\ b^\eps_\alpha
\leq b^\eps_\beta$.
\end{enumerate}
\mn
As $\bfB^*_\eps$ has length $\leq \mu_\eps$ and $J$ is a chain,
necessarily for some $Y_\eps \subseteq \mu^+_\eps$ with 
$|Y_\eps| = \mu^+_\eps$ we have
\mn
\begin{enumerate}
    \item[(e)]  $b^\eps_\alpha = b^\eps$ for $\alpha \in Y_\eps$.
\end{enumerate}
\mn
We can apply clause $(*)(\delta)$ to $-c_\alpha$ (for $\alpha \in Y_\eps$ and
$\bfB^*_\eps$, and possibly shrinking $Y_\eps$) 
to get $a^\eps_\alpha \in \bfB^*_\eps$ such that:
\mn
\begin{enumerate}
    \item[(f)] $(-c_\alpha)\leq a^\eps_\alpha$ and $0<x\leq a^\eps_\alpha
 \wedge x\in \bfB^*_\eps \Rightarrow x\cap(-c_\alpha)\neq 0$. 
\end{enumerate}
\mn
As above, without loss of generality, shrinking $Y_\eps$ further we get
\mn
\begin{enumerate}
    \item[(g)]  $a^\eps_\alpha = a^\eps$ for $\alpha\in Y_\eps$.
\end{enumerate}
\mn
As the length of $\bfB^*_\eps$ is $\leq \mu_\eps < \mu_\eps^+ = |Y_\eps|$, for
some $\alpha\in Y_\eps$ we have $c_\alpha\notin \bfB^*_\eps$; as
\begin{itemize}
    \item $a^\eps_\alpha \geq -c_\alpha$, $b^\eps_\alpha \geq c_\alpha$, 
    $b^\eps_\alpha \in \bfB^*_\eps$, and $a^\eps_\alpha\in \bfB^*_\eps$,
\end{itemize}
necessarily $a^\eps_\alpha \cap b^\eps_\alpha \neq 0$. 

Hence:
\mn
\begin{enumerate}
    \item[(h)]  $b^\eps \cap a^\eps \neq 0$.
\end{enumerate}
\mn
Let $g(\eps) = \min\{\zeta \leq \eps : a^\eps, b^\eps \in \bfB^*_\zeta\}$, 
so for limit $\eps$, $g(\eps) < \eps$. Hence on some stationary 
$S \subseteq \kappa$ and $\zeta_* < \kappa$ the function $g \rest S$ is
constantly $\zeta_*$, and without loss of generality 
\[
\zeta < \eps \in S \Rightarrow \big| \{\alpha\in Y_\zeta : c_\alpha \in
\bfB^*_\eps\} \big| = \mu^+_\zeta.
\]

\mn
As $\bfB$ satisfies the ccc we can find $\eps_1<\eps_2$ in $S$ such
that
\[
b^{\eps_1} \cap a^{\eps_1} \cap b^{\eps_2} \cap a^{\eps_2} \ne 0.
\]
Choose $\alpha\in Y_{\eps_1}$ such that $c_\alpha\in \bfB_{\eps_1}$
and $\beta\in Y_{\eps_2}$. Now $\{c_\alpha,c_\beta\}$ is independent: 
a contradiction. 

\mn
2) Similarly.
\end{PROOF}

\begin{remark}\label{5.6p}
We may further ask: is the restriction ``$\cf(\mu) > \aleph_0$" 
in (\ref{5.17}) necessary?
\end{remark}

\begin{observation}\label{5.7bis}
Assume that the infinite Boolean algebra $\bfB$ has the length $\mu$,
$\cf(\mu) = \aleph_0$. \underline{Then}  the length is obtained.
\end{observation}

\begin{PROOF}{\ref{5.7bis}}
Let $\dcI = \{b \in \bfB : \mathrm{length}(\bfB\rest b) < \mu\}$.

Easily
\[
b_1 \leq b_2\wedge b_2 \in \dcI \Rightarrow b_1 \in \dcI.
\]

\mn
Also clearly $\dcI$ is closed under unions. 

\sn
[Why? If $b_1,b_2\in \dcI$, 
$b = b_1 \cup b_2 \notin \dcI$ then there is a chain $\LL c_t : t \in J\RR$, 
$J$ a linear order of cardinality $\mu$, 
$[s <_J t\, \Rightarrow\, c_s <_\bfB c_t]$, and $c_t\leq b$. 

Let
\[
\bfE_l \defeq \{(t,s) \in J \times J : c_t \cap b_l = c_s\cap b_l\}.
\]
Then $\bfE_l$ is a convex equivalence relation on $J$; if $|J/\bfE_l| = \mu$ then
$\{c_t \cap b_l : t \in J\}$ exemplifies $b_l \notin \dcI$, a contradiction. 
So $|J/\bfE_l| < \mu$. Hence $\bfE = \bfE_1 \cap \bfE_2$ is a convex
equivalence relation with $\leq |J/\bfE_1| \times |J/\bfE_2| < \mu$ classes, 
but as $b = b_1 \cup b_2$ it is the equality.]

\smallskip
If $\bfB/\dcI$ is infinite then we can find $\LL a_n/\dcI : n < \omega\RR$ 
pairwise disjoint non-zero. Now $b_n \defeq a_n - \bigcup\limits_{\ell<n} a_\ell$ 
are pairwise disjoint members of $\bfB$ not in $\dcI$.  Let 
$\mu = \sum\limits_{n<\omega} \mu_n$, $\mu_n < \mu$. Let $\LL c^n_t : t \in J_n\RR$ 
be an increasing chain in $\bfB \rest b_n$, $|J_n| \geq \mu_n$
(note that we can invert $J_n$). Let $J = \sum\limits_{n<\omega}J_n$
(without loss of generality, $n < m\ \Rightarrow\ J_n \cap J_m = \varnothing$) 
and for $t\in J_n$ let $c^*_t = b_0 \cup \ldots \cup b_{n-1} \cup c^n_t$. Now 
$\LL c^*_t : t \in J\RR$ exemplifies that the length is obtained. So $\bfB/\dcI$ 
is finite, so without loss of generality  $\dcI$ is a maximal ideal. Try to choose 
$a_n \in \dcI$ satisfying $\bigwedge\limits_{\ell<n} a_n \cap a_\ell = 0$ such that
$\mathrm{length}(\bfB \rest a_n) > \mu_n$. {If we succeed, then we may repeat the
proof for the case} ``$\bfB/\dcI$ is infinite,'' hence we necessarily fail. 
Hence for some $n$ (replacing $\bfB$ by $\bfB\rest -(a_0\cup \ldots \cup a_{n-1})$) 
we have 
\[
b \in \dcI \ \Rightarrow \ \mathrm{length}(\bfB\rest b)\leq
\mu_n.
\]

Let $J \subseteq \bfB$ be linearly ordered, $|J| > \mu^+_n$. 
Possibly shrinking $J$, without loss of generality 
$J \subseteq \dcI \vee J \subseteq \bfB \setminus \dcI$.  As we
can replace $J$ by $\{\mathbf{1}_{\bfB}-b : b \in J\}$ without loss of generality  
$J\subseteq \dcI$, so for some $b \in J$ we have $|\{c \in J : c \le b\}| \ge \mu_n^+$, 
and hence $\mathrm{length}(\bfB \rest b) \ge \mu_n^+$, a contradiction.
\end{PROOF}

\begin{remark}
We may wonder if we can replace $\aleph_0$ in \ref{5.17} by another cardinals. 
Most natural are $\kappa$ strong limit of cofinality $\omega$.
\end{remark}

\begin{claim}\label{5.8}
Assume $\kappa\leq\mu<2^\kappa$, $\kappa$ strong limit and $\cf(\kappa)=
\aleph_0$. \underline{Then} there is a Boolean algebra $\bfB$ such that:
\mn
\begin{enumerate}
    \item[$(\alpha)$]  $|\bfB|=2^\kappa$
\sn
    \item[$(\beta)$]  $\bfB\models \kappa^+$-cc (and even the strong $\lambda$-cc, when $\lambda = \cf(\lambda) > {\kappa}$). 
\sn
    \item[$(\gamma)$]  $\bfB$ has length $\mu$ (and satisfies clause $(B)^+$ of \ref{5.4})
\sn
    \item[$(\delta)$]  $\bfB$ has no homomorphic image $\bfB'$ with 
    $|\bfB'| \in [\kappa, 2^\kappa)$.
\end{enumerate}
\end{claim}

\begin{PROOF}{\ref{5.8}}
Let $h\in {}^\omega\omega$ be $h(n) = 2(n+1)$. Let $\bfB_0 \subseteq \clP(\kappa)$ 
have cardinality $\mu$ and length $\mu$, with a dense $\bfB \subseteq \bfB_0 \setminus \{\mathbf{0}\}$ of cardinality $\kappa$.

\sn
[Why does such a $\bfB_0$ exist? Let $\kappa = \sum\limits_{n<\omega} \kappa_n$ with $\kappa_n < \kappa_{n+1}$, and let $\nu_\alpha \in \prod \kappa_n$ be pairwise distinct for $\alpha < \mu$. Let $\eta_\alpha \in {}^\omega\kappa$ be defined by 
$$
\eta_\alpha(n) \defeq (\kappa_n)^n \eta_\alpha(0) + (\kappa_n)^{n-1} \eta_\alpha(1) + \ldots + (\kappa_n)^1 \eta_\alpha(n-1) + (\kappa_n)^0 \eta_\alpha(n)
$$
(where $(\kappa_n)^\ell$ denotes ordinary ordinal exponentiation).
\begin{itemize}
    \item $\{\eta_\alpha : \alpha < \mu\}$ is linearly ordered under the lexicographic order $<_\lex$.
\end{itemize}
Let $A_\alpha \defeq \{\kappa_n^+ + \eta_\alpha(n) : \alpha < \mu,\ n < \omega\}$ and $\bfB_0$ be the Boolean subalgebra $\clP(\kappa)$ generated by $[\kappa]^{<\aleph_0} \cup \{A_\alpha : \alpha < \mu\}$. All our assertions are easy ro check.]

Let
\begin{align*}
I^0_\alpha = \big\{\eta : &\ \eta\text{ is an $\omega$-sequence, $\eta(n)$ is an
increasing }\ \\
&\text{ sequence of ordinals $<2^\kappa$ of length }h(n) \big\},
\end{align*}

\mn
and
\[
I_\alpha = I^0_\alpha \cup \big\{\Res^\ell_n(\eta) : n < \omega,\ \ell < h(n),\ 
\eta \in I^0_\eta \big\}
\]

\mn
so $|I_\alpha| = 2^\kappa$. Let 
$\bfB_{\alpha+1} = \ba[\bfB_\alpha,\bar a_\alpha, I_\alpha]$, $\bfB_\alpha$ 
increasing continuous for $\alpha\leq 2^\kappa$. (Again, $\bar a_\alpha$ is an $\omega$-sequence of pairwise disjoint non-zero elements of $\bfB_\alpha$ 
satisfying `each such sequence appears $2^\kappa$ times.')

Again, for $\alpha < \beta$, $\bfB_\alpha \lessdot \bfB_\beta$ (and even the
conclusion of \ref{5.2}(5) holds). The proof that $\bfB \defeq \bfB_{(2^\kappa)}$
satisfies \ref{5.8}$(\alpha)$, $(\beta)$, $(\gamma)$ is as in the proof of
\ref{5.4}.

\noindent
For $(\delta)$ we need \ref{5.19} below.
\end{PROOF}

\begin{observation}\label{5.19}
Assume that $\kappa$ is a strong limit cardinal of countable cofinality.

\sn
$1)$ If $\bfB'$ is a Boolean Algebra of cardinality $\geq \kappa$ but
$<2^\kappa$ \underline{then}:
\mn
\begin{enumerate}
    \item[$(a)$]  There are pairwise disjoint non-zero $b_n$ (for $n < \omega$) in $\bfB'$ such that
    \begin{enumerate}
        \item[$(*)$]  for no $c\in \bfB'$ do we have
        $\bigwedge\limits_{n<\omega} [b_{2n}\leq c] \wedge \bigwedge\limits_{n<\omega} [b_{2n+1} \cap c = 0]$.
    \end{enumerate}
\end{enumerate}

\sn
$2)$  For a Boolean algebra $\bfB'$, a sufficient condition for $\bfB'$ to
satisfy $(a)$ (i.e.~the existence of a sequence 
$\LL b_n : n < \omega\RR$ of pairwise disjoint elements of 
$\bfB'$ satisfying $(*)$ above) is:
\mn
\begin{enumerate}
    \item[$(b)$]   $\bfB'$ has cardinality $< 2^\kappa$ and there are $b_n \in \bfB'$ such that $\bigwedge\limits_{n<m} [b_n \cap b_m = 0]$ and 
    $\kappa = \liminf\limits_n |\bfB' \rest b_n|$.
\end{enumerate}
\end{observation}

\noindent
We first prove that \ref{5.19} suffices (for finishing the proof of \ref{5.8}).
Toward contradiction assume that $\bfB'$ is a Boolean 
algebra of cardinality $< 2^\kappa$ but $\ge \kappa$, 
and $\bfB'$ is a homomorphic image of $\bfB$. 
If clause (a) is satisfied by $\bfB'$, then the proof is very similar to the
earlier proof of \ref{5.4}: for a homomorphism $\bff : \bfB \to \bfB'$
from $\bfB$ onto $\bfB'$ we can find pairwise disjoint $a_n\in \bfB$ (for
$n<\omega$) such that $f(a_n) = b_n$. So, for some $\alpha$ we have 
$\bar a_\alpha = \LL a_n : n < \omega\RR$, and we repeat the relevant 
part of \ref{5.4}. Using clauses (b),(c) of Definition \ref{5.1} we get 
a contradiction.  We are left with proving \ref{5.19}. First, the second part.

\begin{PROOF}{\ref{5.19}}
\textbf{Proof of Observation \ref{5.19}(2)}:

We can find $\bar c^\zeta = \LL c^\zeta_n : n < \omega\RR$, $c^\zeta_n\in \bfB'$,
$c^\zeta_n\leq b_n$ for $\zeta < 2^\kappa$ such that the sequences 
$\LL c^\zeta_n : n < \omega\RR$ are pairwise distinct for $\zeta < 2^\kappa$. 
For each $\zeta$ let $b^\zeta_{2n} = c^\zeta_n$, $b^\zeta_{2n+1} = b_n - c^\zeta_n$, 
so if clause \textbf{(a)} fails then for every $\zeta < 2^\kappa$ there is 
$y_\zeta\in \bfB'$ such that for every $n < \omega$ we have
\[
b^\zeta_{2n} \leq y_\zeta,\quad b^\zeta_{2n+1} \cap y_\zeta = 0.
\]
So $\bigwedge\limits_n [y_\zeta \cap b_n = c^\zeta_n]$ and hence 
$\zeta < \xi < 2^\kappa \Rightarrow\, y_\zeta \neq y_\xi$, which contradicts
$|\bfB'| < 2^\kappa$.

\bn
\textbf{Proof of Observation \ref{5.19}(1)}:

Assume that the conclusion fails. For a cardinal $\mu$, let
\[
\dcI_\mu = \dcI_\mu[\bfB'] \defeq \{b \in \bfB' : \bfB' \rest b
\text{ has cardinality}<\mu\}.
\]

\mn
Clearly it is an ideal of $\bfB'$ increasing with $\mu$ and 
$\mathbf{1}_{\bfB'} \in \dcI_\mu \Leftrightarrow \mu > |\bfB|$. 
If $\bfB'/\dcI_\kappa[\bfB']$ is infinite then we can easily get condition 
\textbf{(B)} of part (2), and we are done. If it is finite, but 
$\dcI_\mu[\bfB] \ne \dcI_\kappa[\bfB']$ for every $\mu < \kappa$, then let 
$\kappa = \sum\limits_{n<\omega} \mu_n$, $\mu_n < \mu_{n+1}$, and choose 
$b_n \in \dcI_\kappa[\bfB'] \setminus \dcI_{\mu_n}[\bfB']$. But 
$\dcI_\kappa[\bfB'] = \bigcup\limits_{\mu<\kappa} \dcI_\mu[\bfB']$, so 
$\bigwedge\limits_n \bigvee\limits_m b_n \in \dcI_{\mu_m}[\bfB']$. So without 
loss of generality $b_n \in \dcI_{\mu_{n+1}}[\bfB'] \setminus \dcI_{\mu_n}[\bfB']$ 
and hence $\LL b_n - \bigcup\limits_{l<n} b_l : n < \omega\RR$ are as required. 
We are left with the case that for some $\mu(*) < \kappa$,
\[
\dcI \defeq \dcI_{\mu(*)}[\bfB'] = \dcI_\kappa[\bfB']
\]

\mn
and without loss of generality  $\dcI = \dcI_{\mu(*)}[\bfB']$ is a maximal ideal.

Without loss of generality  $2^{\mu(*)} < \mu_n < \mu_{n+1}$ for $n < \omega$. 
Let $b_i \in \dcI$ (for $i < \kappa$) be distinct (these exist as
$|\bfB'| \geq \kappa$ and $\dcI$ is a maximal ideal of $\bfB'$). 
By the proof of Erd\H{o}s--Tarski theorem, without loss of generality  
$\LL b_i : i < \kappa\RR$ are non-zero pairwise disjoint.

\mn
[Why? For example, apply the $\Delta$-system lemma to
\[
\big\{ \{x : x \leq b_i\} : i < (2^{\mu_n})^+ \big\},
\]

\mn
and get $Y_n\subseteq (2^{\mu_n})^+$ of cardinality $(2^{\mu_n})^+$ and a
set $A_n$ of cardinality $\leq 2^{\mu(*)}$ such that
\[
i,j \in Y_n \wedge i \neq j\ \Rightarrow\ \{x : x \leq b_i\} \cap 
\{x : x \leq b_j\} = A_n.
\]

\mn
So $|A_n| \leq \mu(*)$. Pick $Y'_n \subseteq Y_n$ of cardinality
$(2^{\mu_n})^+$ such that
\[
i,j \in Y_n' \wedge i \neq j\, \Rightarrow\, \{x : x \leq b_i\} \cap
\bigcup_{m<n} A_m' = \{x : x \leq b_j\} \cap \bigcup_{m<n} A_m',
\]

\mn
where $A'_n = \big\{x : (\exists i\in Y_n')[x \leq b_i] \big\}$. Let 
$i(n) = \min(Y_n')$. Then 
\[
X_n = \{x_i - x_{i(n)} : i \in Y_n,\ i > i(n)\} \subseteq \bfB \setminus \{0\}
\]

\mn
is an antichain, and $\bigcup\limits_n X_n$ is as required.]

Let
\[
\cP_0 = \big\{Y \in [\kappa]^{\aleph_0} : \text{there is } b \in \dcI 
\text{ such that }(\forall i\in Y)[b_i\leq b]\big\}.
\]

\mn
This is a subset of $[\kappa]^{\aleph_0}$ of cardinality $\leq |\dcI|
\cdot\mu(*)^{\aleph_0}\leq |\bfB'|+\kappa=|\bfB'|$, but\\
$[\kappa]^{\aleph_0}=  2^\kappa>|\bfB'|$, so there is $Y_0\in
[\kappa]^{\aleph_0}\setminus \cP_0$. 

Let
\[
\cP_1 = \big\{Y\in [\kappa]^{\aleph_0} : Y \subseteq \kappa \setminus Y_0 
\text{ and }(\exists b \in \dcI)(\forall i\in Y)[b_i\leq b] \big\}.
\]

\mn
By cardinality considerations as above there is 
$Y_1 \in [\kappa]^{\aleph_0} \setminus \cP_1$ disjoint to $Y_0$. 
By assumption above (i.e.~clause \textbf{(a)} fails) there is $b\in \bfB'$ 
such that $\bigwedge\limits_{i\in Y_0} [b_i \leq b]$ and 
$\bigwedge\limits_{i\in Y_1} \big[b_i \leq (1-b) \big]$. If $b\in \dcI$ we
get a contradiction to the choice of $Y_0$; if not, then $\mathbf{1}_\bfB - b \in \dcI$
contradicts the choice of $Y_1$. Hence the observation holds and hence the
Observation \ref{5.19} is proven.  Hence Claim \ref{5.8} is proven.
\end{PROOF}

\begin{remark}\label{5.19A}
In other words \ref{5.19} says
\mn
\begin{enumerate}
    \item[$(*)$]  If $\kappa$ is strong limit, $\cf(\kappa) = \aleph_0$ and $\bfB$ 
    is a Boolean algebra of cardinality $\geq\kappa$ with $\aleph_1$-separation 
    (i.e.~\textbf{(a)} of the observation fails) then $|\bfB| \ge 2^\kappa$.
\end{enumerate}
\end{remark}

\newpage

\section {Using subtrees of $({}^{\omega\geq} 2,\lhd)$ and
theories unstable in $\aleph_0$}
\sectionmark{Subtrees...}
\label{par3}

\begin{theorem}\label{3.1}
Suppose $T \subseteq T_1$ are first order theories, $T_1$ is countable, $T$
is complete, superstable but $\aleph_0$-unstable. \underline{Then} for
$\lambda > \aleph_0$ we have
\[
\numbIE(\lambda,T_1,T) \ge \min\{2^\lambda,\beth_2\}.
\]
\end{theorem}

\begin{remark}
The reader is not required to know
anything on superstable theories, just to believe a result quoted below. So
we can just assume $(*)$ from the proof.
\end{remark}

\begin{PROOF}{\ref{3.1}}
The assumption that the theory is superstable and not totally transcendental 
($=\aleph_0$-stable) is used to obtain $m_a, m_b < \omega$ and a countable 
set of definable (without parameters) equivalence relations 
$\{\varphi_n(\bar x;\bar y) : n < \omega\} \subseteq \bbL(\tau_T)$ 
such that:\footnote{
    We may write $\bar x\; \varphi_n\; \bar y$ instead of $\varphi_n(\bar{x},\bar{y})$.
}
\mn
\begin{enumerate}
    \item[$(*)$]  
    \begin{enumerate}[$(i)$]
        \item $\lh(\bar x) = \lh(\bar y) = m_a + m_b$

        \item If $M$ is a model of $T$ and $\bar a \in {}^{m_a}|M|$ \underline{then} the set $\{\bar a \caret \bar b/\varphi_n : \bar b \in {}^{m_b}|M|\}$ is finite.

        \item  If for $\ell=1,2$, $\lh(\bar a_\ell) = m_a$, $\lh(\bar b_\ell) =   m_b$, 
        and $(\bar a_1 \caret \bar b_1)\ \varphi_n\ (\bar a_2 \caret \bar b_2)$ then
        $\bar a_1 = \bar a_2$.

        \item $\varphi_{n+1}$ refines $\varphi_n$: i.e. for every $n < \omega$,
        $\bar x\ \varphi_{n+1}\ \bar y$ implies $\bar x\ \varphi_n\ \bar y$.

        \item There are (in some model $M$ of $T$)  $\bar c_\eta$, for 
        $\eta \in {}^{\omega>}2$, such that:
        \begin{enumerate}
            \item $\lh(\eta)\geq n\wedge \lh(\nu)\geq n \Rightarrow 
            \big[\bar c_\eta\; \varphi_n\; \bar c_\nu \Leftrightarrow 
            \eta \rest n = \nu \rest n \big]$
\sn
            \item $\bar c_\eta \rest m_a = \bar c_\nu \rest m_a$ and
            $\lh(\bar c_\eta) = m_a + m_b$.
        \end{enumerate}
    \end{enumerate}
\end{enumerate}
\sn
The existence of this set of equivalence relations was proved in Chapter III, 
5.1-5.3 of both \cite{Sh:a} and \cite{Sh:c}.

Clearly, without loss of generality we may expand the theory $T_1$. Let
$$
\{c_\ell : \ell < m_a\} \cup \{c_{\eta,\ell} : \ell \in [m_a,m_a + m_b]
\text{ and } \eta \in {}^{\omega>}2\}
$$ 
be new individual constants in $T_1$. We let 
$$
\bar{c}_\eta = \big\LL c_\ell : \ell < m_a \big\RR \caret \big\LL c_{\eta,\ell} : \ell \in [m_a, m_a + m_b) \big\RR,
$$ 
and suppose
\begin{align*}
(**) \qquad  T_1 \supseteq \big\{\ \ (\bar c_\eta\ \varphi_n\ \bar c_\nu) : &\ \lh(\eta), \lh(\nu) \geq n,\ \eta \rest n = \nu \rest n \big\}\ \cup \\
\big\{ \neg (\bar c_\eta\ \varphi_n\ \bar c_\nu) : &\ \lh(\eta),\lh(\nu)\geq n,\  
\eta \rest n \neq \nu \rest n\big\}. 
\end{align*}

\mn
Also without loss of generality, suppose that $T_1$ has Skolem functions 
(so the axioms saying it has Skolem functions belong to $T_1$).

We will use the following fact. [For a sequence $\bar\eta$ let 
$\bar\eta = \big\LL\bar\eta[\ell] : \ell < \lh(\bar\eta) \big\RR$ and
$\bar{a}_{\bar{\eta}} = \bar{a}_{\bar\eta[0]} \caret
\bar a_{\bar\eta [1]} \caret \bar a_{\bar\eta [2]}\ldots$.] 
\end{PROOF}

\begin{fact}\label{3.1A}
1)  Suppose
\begin{enumerate}
    \item  $T\subseteq T_1$ are first order theories, $T$ complete and superstable, unstable in $|T_1|$, $\tau = \tau(T)$ and $\tau_1 = \tau(T_1)$, and $T_1$ has Skolem functions.
\sn
    \item  $\tau_1$ is countable, or at least $\MA_\mu$ holds for $\mu = |T_1|$. 
\sn
    \item  $\varphi_n$ (for $n < \omega$), $m_a, m_b$ are as in $(*)$ above, and 
    $m_* \defeq m_a + m_b$.
\sn
    \item $\varphi_n \in \tau$ is a $(2m_*)$-place predicate, 
    $$
    \Delta = \{\varphi_n : n < \omega\},\quad \tau^+_1 = \tau_1 \cup \{d_n : n < \omega\},
    $$ 
    $\tau \subseteq \tau_1$, and $|\tau_1| \leq \mu < 2^{\aleph_0}$.
\sn
    \item $T_1$ satisfies $(**)$ above.
\end{enumerate}
\mn
\underline{Then} there are $M_1$ and $\bar a_\eta$ (for $\eta\in {}^\omega2$) such that:
\mn
\begin{enumerate}
    \item[$(\alpha)$]  $M_1$ is a model of $T_1$ and $\varphi_n^{M_1}$ is an equivalence relation such that $\varphi^{M_1}_{n+1}$ refines $\varphi_n^{M_1}$.
\sn
    \item[$(\beta)$]  $\tau(M_1) = \tau_1^+$, 
    $\{\bar a_\eta : \eta \in {}^\omega 2\} \subseteq {}^{m_*}|M_1|$, and
    $$\lh(\eta) \geq n \wedge \lh(\nu) \geq n\ \Rightarrow\ \big[\eta \rest n = \nu\rest n \Leftrightarrow (\bar a_\eta\ \varphi_n\ \bar a_\nu)\big].$$
    For $\bar \eta \in {}^m({{}^\omega 2})$, let $\bar a_{\bar \eta} = \bar a_{\eta_0} \caret \ldots \caret \bar a_{\eta_{m-1}}$.
\sn
    \item[$(\beta)_1$]  $\bar a_\eta \rest m_a = \bar a_\nu \rest m_a = 
    \LL c^{M_1}_\ell : \ell < m_a\RR$, $\lh(\bar a_\eta) = m_*$, 
    and if $n < \omega$, $\lh(\bar a) = m_a < m_*$ \underline{then} 
    $\big|\big\{ {\bar a} \caret {\bar b}/\varphi_m : \bar b \in {}^{m_b}(M_1) \big\} \big| < k_m$.
\sn
    \item[$(\gamma)$]  For every formula $\varphi(\bar x)$ from $\bbL(\tau_1)$ such that $m_*$ {divides} $\lh(\bar x)$, there is $\eta_\varphi$ such that for $n \in [\eta_\varphi,\omega)$:
    \begin{enumerate}
        \item[$(*)^1_{\varphi,n}$]  If $\bar\eta$, $\bar\nu\in {}^m({}^\omega2)$,
        $\lh(\bar\eta) = \lh(\bar\nu) = m = \dfrac{\lh(\bar x)}{m_a + m_b}$ 
        (so $\lh(\bar a_{\bar\eta}) = \lh(\bar x)$), and
        $$\big\LL\eta_\ell \rest \eta : \ell < \lh(\bar\eta)\big\RR = \big\LL \nu_\ell \rest n : \ell < \lh(\bar\nu) \big\RR$$ is without repetitions, then 
        $M_1 \models \varphi[\bar a_{\bar\eta}] = \varphi[\bar a_{\bar\nu}]$.
    \end{enumerate}
\sn
    \item[$(\delta)$]  $\LL d_n : n < \omega\RR$ is an indiscernible sequence 
    over $\{\bar a_\eta : \eta \in {}^\omega 2\}$ in $M_1 \rest \tau_1$.
\sn
    \item[$(\delta)^+$]  $d_n \ne d_m$ for $n\neq m$.
\end{enumerate}
\mn
2)  If $M_1,\tau,\tau_1,\tau_1^+,m_a,m_b,\LL \varphi_n : n < \omega \RR$ are as in $(\alpha)$, 
$(\beta)$, $(\beta)_1$, $(\gamma),(\delta),(\delta)^+$ above and $\mu = \aleph_0$ (or at least $\MA_\mu$) \underline{then},
replacing ${}^\omega 2$ by a subtree, replacing $\LL \varphi_n : n < \omega\RR$ 
by a sub-sequence and renaming, decreasing $M_1$, we can add to part (1):
\mn
\begin{enumerate}
    \item[$(\gamma)^+$] For every sequence of terms $\bar\sigma(\bar x)$ from 
    $\tau^+_1$, if $m \times (m_a+m_b) = \lh(\bar x)$, $m_a + m_b = \lh(\bar\sigma)$,
    $\bar\sigma (\bar x) \rest m_a = (\bar \sigma \rest m_a)(\bar x \rest m_d)$, 
    $m_e < m_a$, $m_d = m_e \times (m_a+m_b)$, [i.e. 
    $\bar\sigma(\bar a_{\bar\eta})\rest m_a = 
    (\bar\sigma\rest m_a)(\bar a_{\bar\eta\rest m_e})$ for 
    $\bar\eta\in {}^m(^\omega 2)$], \underline{then} there exists 
    $n_{\bar\sigma} < \omega$ such that:
\sn
    \begin{enumerate}
        \item  For $n \geq n_{\bar\sigma}$ and $\bar\eta,\bar\nu \in {}^m({}^\omega 2)$ with no repetitions, $\bar\eta \rest m_e = \bar\nu \rest m_e$, we have:
        \begin{itemize}
            \item If $\ell \neq k \Rightarrow \bar\eta[\ell]\rest n \neq \bar\eta[k] \rest n$ and $(\forall\ell < m)\big[\bar\eta[\ell] \rest n = \bar\nu[\ell] \rest n\big]$ 
            \underline{then} for every $\bar\rho\in {}^m({}^\omega 2)$, 
            $\bar\rho\rest m_e = \bar\eta \rest m_e$ implies
        \[
            \big(\bar\sigma(\bar a_{\bar\eta})\ \varphi_n\ \bar\sigma (\bar a_{\bar\rho}) \big) \Leftrightarrow \big(\bar\sigma(\bar a_\nu)\ \varphi_n\ \bar\sigma (\bar a_{\bar\rho}) \big).
        \]
        \end{itemize}        
        
        \item  For $n\geq\eta_{\bar\sigma}$ and $\bar\eta,\bar\nu \in {}^m({}^n2)$ 
        each with no repetition and 
        $$\bar\eta \rest m_e = \bar\nu \rest m_e,$$ we have:
        \begin{itemize}
            \item If there are $k\geq n$ and $\bar\eta_1,\bar\nu_1 \in {}^m({}^\omega 2)$ such that $\neg \varphi_k(\bar\sigma(\bar a_{\bar\eta_1}),\bar\sigma(\bar a_{\bar\nu_1}))$, for $\ell < m$, $\bar\eta_1[\ell] \rest n = \bar\eta[\ell]$, $\bar\nu_1[\ell]\rest n=\bar\nu [\ell]$, and
        $$
            (\forall\ell,i < m)\big[\bar\eta_1[\ell] = \bar\nu_1[i] \Leftrightarrow \bar\eta[\ell] = \bar\nu[i]\big],
        $$
            \underline{then} for every $\bar\eta^*,\bar\nu^* \in {}^m({}^\omega 2)$
            satisfying $\bar\eta^*[\ell]\rest n = \bar\eta[\ell]$, 
            $\bar\nu^*[\ell] \rest n=\bar\nu[\ell]$ (for each $\ell < m$) and
        $$
            (\forall\ell,i < m)\big[\bar\eta^*[\ell] = \bar\nu^*[i] \Leftrightarrow \bar\eta[\ell] = \bar\nu[i]\big]
        $$
            we have $\neg \big[\bar\sigma(\bar a_{\bar\eta^*})\ \varphi_n\ \bar\sigma(\bar a_{\bar\nu^*}) \big]$.
        \end{itemize}        
    \end{enumerate}
\end{enumerate}
\end{fact}

\begin{remark}
1) This is the only place where countability (or $\MA_{\vert\tau_1\vert}$) is used.

\sn
2)  For alternative proof see \ref{3.10}.
\end{remark}

\begin{PROOF}{\ref{3.1A}}
1) If we ignore $(\delta)^+$ (so can have $d_n = d_0$) use Theorem 
\cite[Ch.VII,3.7]{Sh:a}. In general, use \cite[Ch.VII,Ex.3.1]{Sh:a}.  
What if $T_1$ is uncountable but $\MA_\mu$? (The reader may ignore this
proof or see the proof of \ref{3.10}.)

Let $\bbP$ be the forcing notion of adding $\lambda = \beth_{(2^\mu)^+}$
Cohen reals, $\LL\name{\eta}_i : i < \lambda\RR$, $\name{\eta}_i \in {}^\omega 2$. 
Let $\chi = (2^\lambda)^+$ and let 
$$
\Vdash_\bbP ``\Name{M} \text{ is a model of $T_1$, the Skolem hull of } 
\{x_i : i < \lambda\},\ \bar x_i\ \varphi_m\ \bar c_{\name{\eta}_i\rest m}".
$$ 
By the Omitting Type Theorem\footnote{
    See (e.g.) \cite[Ch.VII,\S5]{Sh:c}.
} 
there are $\gB_1 \prec \gB_2$ with $\gB_1 \prec (\clH(\chi),\in,<^*_\chi)$ and 
$\|\gB_1\| = \mu$ such that $T_1$, $P$, $\Name{M}$, $\LL x_i : i < \lambda\RR$ 
belong to $\gB_1$. Also in $\gB_2$, $\LL a_\rho : \rho \in {}^\omega 2\RR$ 
is an indiscernible sequence over $\gB_1$, and $\gB_2 \models ``\name{a}_i$ 
is an ordinal $=\lambda$".

Note that any set which $\gB_2$ considers a maximal antichain of $\bbP^{\gB_2}$ really is so. 
Now we can naturally apply $\MA_\mu$.

\sn 
2) Satisfy requirement (a) by letting 
$\varphi^\ell_n(\bar x \caret \bar z) \defeq E_n\big(\bar x \caret \bar z,
F_\ell(\bar z) \caret \bar z \big)$ for $\ell < \ell^*_n < \omega$, 
where $F_\ell\in \tau^+$ are such that $\{F_\ell(\bar x) : \ell < \ell^*_n\}$ 
is a complete set of representatives for $\{\bar x \caret \bar z/\varphi_n : \bar x\}$, 
possibly with repetitions. (Remember $T_1$ has Skolem functions and there is 
$\ell^*_n$ which does not depend on $\bar z$ by compactness).  
Requirement (b) is fulfilled by trimming the perfect tree and renaming. 
\end{PROOF}

\begin{claim}\label{3.2}
For $M_1$, $\bar a_\eta$ ($\eta\in {}^\omega 2$), $\varphi_n$ as in the
conclusion of \ref{3.1A} we can conclude:
\begin{enumerate}
    \item[$\otimes$]  \underline{If} $\nu \neq \rho$ are from ${}^\omega 2$,
     $\bar\eta_\nu = \LL\eta_{\nu,\ell} : \ell < \ell(*)\RR$,
    $\bar\eta_\rho = \LL\eta_{\rho,\ell} : \ell < \ell(*)\RR$,\\
    $\bar x = \LL x_\ell : \ell < \ell(*)\RR$, 
    $\bar\sigma(\bar x) = \LL\sigma_m(\bar x) : m < m(*)\RR$, 
    $\nu\rest k = \rho \rest k$,\\ $\eta_{\nu,\ell} \rest k = \eta_{\rho,\ell} \rest k$, 
    $\LL \eta_{\nu,\ell} : \ell < \ell(*)\RR$ with no repetitions, $k > n_{\bar\sigma}$, and
    \[
        \bigwedge\limits_{n<\omega} \big[\bar a_\nu\ \varphi_n\ \bar a_\rho \Leftrightarrow \bar\sigma(\bar a_{\bar\eta_\nu})\ \varphi_n\ \bar\sigma (\bar a_{\bar\eta_\rho}) \big]
    \]
    (moreover, the $\Delta$-type of $\bar a_\nu \caret \bar a_\rho$ and 
    $\bar\sigma(\bar a_{\bar\eta_\nu}) \caret \bar\sigma (\bar a_{\bar\eta_\rho})$ 
    (in $M$) are equal for every $n$) \underline{then} 
    $\lh(\nu\cap\rho) \in \big\{\lh(\eta_{\nu,\ell} \cap \eta_{\rho,\ell}) : \ell < \ell(*) \big\}$.
\end{enumerate}
\end{claim}

\begin{PROOF}{\ref{3.2}}
Assume not. 

Let $n = \lh(\rho\cap\nu)$. Then $\varphi_n(\bar a_\rho,\bar a_\nu) \wedge 
\neg \varphi_{n+1}(\bar a_\rho,\bar a_\nu)$. We suppose first (for didactic reasons) 
for the sake of contradiction that for every $\ell< n_0$ we have
\[
\bar\eta_\nu[\ell] \neq \bar\eta_\rho[\ell]\ \Rightarrow\ 
\lh(\bar\eta _\nu[\ell] \cap \bar\eta_\rho[\ell]) < n.
\]

\mn
By the equality of types $\neg \varphi_{n+1}\bigl(\bar\sigma (\bar
a_{\bar\eta_\rho}),\bar\sigma (\bar a_{\bar\eta_\nu})\bigr)$, now we can
deduce by Fact \ref{3.1A}(2) and the assumption that the conclusion of
$(\otimes)$ fails, that $\neg \varphi_{n+1}\bigl(\bar\sigma(\bar
a_{\bar\eta_\rho}),\bar\sigma (\bar a_{\bar\eta_\nu})\bigr)$. Again, by the
equality of types $\neg \varphi_n(\bar a_\rho,\bar a_\nu)$,
a contradiction to $\varphi_n(\bar a_\rho,\bar a_\nu)$.

Now we deal with the general case, i.e.~we assume
\mn
\begin{enumerate}
    \item[$(*)$]  $(\forall\ell < n_0)\big[\lh(\bar\eta_\nu[\ell] \cap 
    \bar\eta_\rho[\ell]) \neq n \big]$.
\end{enumerate}
\mn
We shall derive a contradiction.

Define $\bar\eta\in {}^{n_0}({}^\omega 2)$:
$$
\bar\eta[\ell]=\begin{cases}
   {\bar\eta}_\rho[\ell]&\text{ if }\bar\eta_\nu [\ell]\rest n\neq
   \bar\eta_\rho [\ell]\rest n,\\
   {\bar\eta}_\nu[\ell]& \text{ otherwise.}
\end{cases}
$$

\mn
Clearly $\bar\sigma(\bar a_\eta)\rest m_a = \bar\sigma(\bar a_{\eta_\rho}) 
\rest m_a = \bar\sigma(\bar a_{\eta_\nu}) \rest m_a$
and $\bar\eta\rest m_e = \bar\eta_\nu\rest m_e = \bar\eta_\rho \rest m_e$, 
and also $\bar\eta$ is with no repetition and 
$\LL \bar\eta[\ell]\rest n : \ell < n_0\RR$ are pairwise distinct.

Since, by the definition of $\bar\eta$, for which 
$\bar\eta[\ell] \rest n = \bar\eta_\rho[\ell] \rest n$, using $(*)$ we obtain
\[
\eta[\ell] \rest (n+1) = \eta_\rho[\ell] \rest (n+1).
\]

\mn
Let $\bar b=\bar\sigma(\bar a_{\bar\eta})$. By reflexivity of the
equivalence relation we have
\[
\bar\sigma(\bar a_{\bar\eta_\rho})\ \varphi_{n+1}\ \bar\sigma
(\bar a_{\bar\eta_\rho}).
\]

\mn
By Fact \ref{3.1A}(1), $\bar\sigma(\bar a_{\bar\eta})\ \varphi_{n+1}\ 
\bar\sigma (\bar a_{\bar\eta_\rho})$; i.e. $\bar b\ \varphi_{n+1}\ 
\bar\sigma (\bar a_{\bar\eta_\rho})$. Finally,\footnote{
    As $\neg \bigl(\bar\sigma (\bar a_{\bar\eta_\nu})\ \varphi_{n+1}\ \bar\sigma (\bar a_{\bar\eta_\rho}) \bigr)$.
} 
using transitivity of the equivalence relation we
have $\neg \varphi_{n+1}\bigl(\bar b,\bar\sigma(\bar a_{\bar\eta_\rho})\bigr)$.

By the definition of $\bar\eta$, for every $\ell < n_0$ we have
\[
\bar\eta[\ell] = \bar\eta_\nu[\ell]\ \text{ or }\ \lh(\bar\eta[\ell]
\cap \eta_\nu[\ell]) < n.
\]
But since $n > k_0$, clearly 
$$
\big| \big\{ \bar\eta[\ell]\rest k_0 : \ell < n_0 \big\} \big| = \big|\big\{ \bar\eta_\nu[\ell] \rest k_0 : \ell < n_0 \big\}\big| = n_0.
$$

So by Fact \ref{3.1A}(2), as 
$\neg \bigl(\bar b\ \varphi_{n+1}\ \bar\sigma (\bar a_{\bar\eta_\nu})\bigr)$ 
(see above), we have 
$\neg \bigl(\bar b\ \varphi_n\ \bar\sigma (\bar a_{\bar\eta_\nu}) \bigr)$. 
But $\bar b\ \varphi_n\ \bar\sigma (\bar a_{\bar\eta_\rho})$ (see above) and
$\bar\sigma (\bar a_{\bar\eta_\rho})\ \varphi_n\ \bar\sigma (\bar a_{\bar\eta_\rho})$, a contradiction. 
\end{PROOF}

\noindent
So for proving theorem \ref{3.1} we can assume

\begin{hypothesis}\label{3.4bis}
$M_1,\LL \varphi_n : n < \omega\RR$, and $\bar a_\eta$ (for $\eta\in {}^\omega 2$)
are as in $(\beta)+(\gamma)$ of \ref{3.1A}(1) and $(\otimes)$ of \ref{3.2}.
\end{hypothesis}

\begin{lemma}\label{3.3}
Assume $\mu < \lambda \leq 2^{\aleph_0}$. We can find $S_\xi \subseteq {}^\omega2$ 
for $\xi < 2^{\aleph_0}$, pairwise disjoint, each of cardinality $\lambda$, such that
\mn
\begin{enumerate}
    \item[$\otimes$]  If $\xi < 2^{\aleph_0}$, $f : S_\xi \to {}^{{\omega >}}
    \bigl(\cM_{\mu,\omega}(\bigcup\limits_{\zeta\neq\xi} S_\zeta)\bigr)$ 
    and $\bfn$ is a function, 
    $$\bfn : \big\{ {\bar\sigma} : (\exists\bar x) \big[{\bar\sigma} = \LL
    \sigma_\ell(\bar x) : \ell < \ell^*\RR\big],\ \sigma_\ell\text{ a term of }
    \bbL_{\mu,\aleph_0}(\tau)\big\} \to \omega$$
    and $\tau$ is the vocabulary of 
    $\cM_{\mu,\omega}(\bigcup\limits_{\zeta\neq\xi} S_\xi)$, \underline{then} 
    we can find $m^*$ (see below) $S^* \subseteq S_\zeta$, $k_0 < \omega$,
    $n_0 = m_a + m_b < \omega$, a sequence $\bar\sigma (\bar x) = 
    \LL\sigma_\ell(\bar x) : \ell < \lh(\bar\sigma)\RR$, with $\lh(\bar x) = n_0$, 
    $\LL\bar\eta_\nu : \nu \in S^*\RR$ and $\bar\eta_0 \in {}^{n_0}({}^\omega 2)$
    with the following properties. Letting $\eta_{\nu,\ell}=\bar{\eta}_\nu[\ell]$:
    \begin{enumerate}
        \item[$(A)$]  $\eta \neq \nu \in S^{*}\ \Rightarrow\ \lh(\eta \cap \nu) > k_0$
\sn
        \item[$(B)$]  For $\nu \in S^*$ we have $\lh(\bar{\eta}_\nu) = n_0$,
        $(\forall\ell < n_0)\left[\bar\eta_{\nu,\ell} \rest k_0 = \bar\eta_{0,\ell} \rest k_0\right]$, and $\{\bar\eta_{\nu,\ell}\rest k_0 : \ell < n_0\} \cup 
        \{\nu \rest k_0\}$ are pairwise distinct.
\sn
        \item[$(C)$]  $k_0 > \bfn(\bar\sigma)$
\sn
        \item[$(D)$]  For each $\ell < n_0$, either 
        $\{\bar\eta_{\nu,\ell} : \nu \in S^*\} = \{\bar\eta_{0,\ell}\}$ or 
        $\{\bar\eta_{\nu,\ell} : \nu \in S^*\}$ are pairwise distinct.
\sn
        \item[$(E)$]  The sets $\{\lh(\nu_1 \cap \nu_2) : \nu_1 \neq \nu_2 \text{ from } S^*\}$ and
        \[
            \{\lh(\eta_{\nu_1,\ell_1} \cap \eta_{\nu_2,\ell_2}) : \nu_1,\nu_2\in S^*
            \text{ and }\ell_1,\ell_2 < n_0\}
        \]
        are disjoint.
\sn
        \item[$(F)$] For every $\nu\in S^*$, $f(\nu) = \bar\sigma(\bar\eta_\nu)$
        (i.e. equal to \\
        $\big\LL\sigma_\ell(\LL\eta_{\nu,n} : n < n_0\RR) : l < m^* \big\RR$).
\sn
        \item[$(G)$] For $\nu_1\neq\nu_2\in S^*$, we have
        \[
            \eta_{\nu_1,\ell} = \eta_{\nu_2,\ell} \Leftrightarrow \ell < m_a 
            \Leftrightarrow \eta_{\nu_1,\ell} = \eta_{0,\ell}.
        \]

        \item[$(H)$]  $S^*$ is $\mu^+$-large. (We say that $S\subseteq {}^\omega2$ 
        is $\chi$-\emph{large} \underline{iff} for every $n < \omega$ and $\nu \in S$ 
        we have $\big|\{\rho\in S : \rho \rest n = \nu \rest n\} \big| \geq \chi$.)
        We can replace $\mu^+$-large by $\lambda$-large if $\cf(\lambda) > \aleph_0$.
\sn
        \item[$(I)$]  $\nu_1,\nu_2 \in S^* \wedge \eta_{\nu_1,\ell_1} = \eta_{\nu_2,\ell_2}$ implies $\ell_1 = \ell_2$.
\sn
        \item[$(J)$]  For $\eta \in \bigcup\limits_\xi S_\xi$, let $\xi(\eta)$ 
        be the unique $\xi$ such that $\eta\in S_\xi$. Now, if
        $\xi(\eta_{\nu_1,\ell_1}) = \xi(\eta_{\nu_2,\ell_2})$ with 
        $\ell_1,\ell_2 < n_0$ and $\nu_1 \ne \nu_2 \in S^*$, then 
        $$\nu\in S^* \Rightarrow \xi(\eta_{\nu_1,\ell_1}) = 
        \xi(\eta_{\nu,\ell_1}) = \xi(\eta_{\nu,\ell_2}).$$
    \end{enumerate}
\end{enumerate}
\end{lemma}

\begin{remark}\label{5.5Anew}
1)  This claim is a version of the ``unembeddability'' 
results;\footnote{
    See Definitions in \cite[\S2]{Sh:331}, results (for example) in VI,  and here in \S1 for the tree $^{\omega\geq} 2$.
} 
well, they are necessarily somewhat weaker {than} in \S1 here.

\sn
2) Of course, we can replace $\bigcup\limits_{\zeta\neq \xi} S_\zeta$ by
$\sum\limits_{\zeta\neq\xi} S_\zeta$.
\end{remark}

\noindent
For proving \ref{3.3} we will use the following combinatorial fact, which is
slightly stronger than Sierpi\'nski's lemma on almost disjoint sets of
integers:

\begin{fact}\label{3.3B}
There are $W(*)$, $\{W_\eta : \eta \in {}^\omega2\}$, and 
$\{U_\eta : \eta \in {}^\omega 2\}$ such that for all $\eta \in {}^\omega 2$:
\mn
\begin{enumerate}
    \item  $W(*),W_\eta$ are infinite subsets of $\omega$.
\sn
    \item  $U_\eta$ is a perfect tree; i.e. $U_\eta\subseteq {}^{\omega>}2$ 
    is downward closed, $\LL\ \RR \in U_\eta$, and
    \[
        (\forall\rho \in U_\eta)(\exists \nu \in U_\eta)\big[\rho \unlhd \nu 
        \wedge \nu \caret \LL 0\RR \in U_\eta \wedge \nu \caret \LL 1 \RR 
        \in U_\eta\big].
    \]

    \item  If $\rho,\nu \in U_\eta$, $\rho \neq \nu$, and $\lh(\rho) = \ell(\nu)$ \underline{then} $\lh(\rho \cap \nu) \in W_\eta$, where $\rho \cap \nu$ is 
    the largest common initial segment of $\rho$ and $\nu$; i.e.
    \[
        \lh(\rho \cap \nu) \defeq \max\{n < \omega : \rho \rest n = \nu \rest n\}.
    \]

    \item  For all $\eta_1 \neq \eta_2\in {}^\omega 2$ and every 
    $\rho \in U_{\eta_1}$, $\nu \in U_{\eta_2}$, there are three possibilities:
    \begin{enumerate}
        \item  $\lh(\rho\cap\nu)\in W_{\eta_1}\cap W_{\eta_2}$
\sn
        \item  $\lh(\rho\cap\nu) \in W(*)$ and $\big(\forall \ell < \lh(\rho \cap \nu)\big) 
        \big[\ell \in W_{\eta_1} \equiv \ell \in W_{\eta_2} \big]$.
\sn
        \item  $\rho \unlhd \nu$ or $\nu \lhd \rho$.
    \end{enumerate}
\sn
    \item  $W(*) \cap W_\eta = \varnothing$
\sn
    \item  For distinct $\eta,\nu$ from $^\omega 2$, we have: 
    \begin{enumerate}
        \item  $W_\eta \cap W_\nu$ is finite (in fact, an initial segment of both).
\sn
        \item  If $\ell\in W\!(*)$ is above $W_\eta \cap W_\nu$ \underline{then} 
        $U_{\eta} \cap U_\nu$ is finite, contained in ${}^{\ell'>}2$ if 
        $\ell < \ell' \in W_\eta \cup W_\nu$, and has no splitting of level $\geq \ell$; 
        i.e. $$\neg(\exists \rho\in{}^{\omega>}2) \big[\lh(\rho) \geq \ell \wedge  
        \{\rho \caret \LL 0\RR,\rho \caret \LL 1 \RR\} \subseteq U_\eta\cap U_\nu \big].$$

        \item  If $\ell \in W(*)$ and $\ell < \sup (W_\eta \cap W_\nu)$ 
        \underline{then} $U_\eta \cap {}^{\ell\geq}2 = U_\nu \cap {}^{\ell\geq}2$.
    \end{enumerate}
\end{enumerate}
\end{fact}

\begin{PROOF}{\ref{3.3B}}
By induction on $n$, define $k(n) = k_n < \omega$, the set $W_n(*) \subseteq k(n)$ 
and the sets $U_\eta\subseteq {}^{k(n)\geq} 2$, $W_\eta\subseteq k(n)$, such that 
in the end (this imposes natural restrictions on them):
\[
\eta \in {}^\omega 2\ \Rightarrow\ W_\eta\cap k_n=W_{\eta\rest n},\quad
U_\eta \cap {}^{k(n)\geq}2 = U_{\eta\rest n},\quad W(*) \cap k(n) = W_n(*).
\]

For $n = 0$, let $k_0 = 0$, $W_n (*) = \varnothing$ and $W_\eta = \varnothing$, 
$U_\eta = \varnothing$ for $\eta\in {}^{n}2$.
For the induction step, choose $k^1(n) = k(n)+n+1$ and for $\eta \in {}^{n} 2$ let 
$$U^1_\eta = U_\eta \cup \{\nu \caret (\eta\rest\ell) : \nu \in U_\eta \cap {}^{k(n)}2,\ \ell \leq n\}.$$

Thus
$$\big(\forall \nu \in {}^{k(n)}2 \cap U_{\eta}\big) \big(\exists !\rho \in {}^{k^1(n)}2 \cap U^1_\eta\big) \big[\nu \unlhd \rho \big].$$

\mn
Define $W_{n+1}(*) = W_n(*) \cup \big[k(n),k^1(n) \big)$. Fix an enumeration
$\{\eta_k : k < 2^{n+1}\}$ of $^{n+1} 2$. Let $k(n+1) \defeq k^1(n)+2^{n+1}$. For
$\eta\in {}^{n+1}2$, there is {a} unique $k < 2^n$ such that $\eta = \eta_k$.
Let
\begin{align*} 
U_{\eta_k} \defeq U^1_{\eta_k\rest n} \cup \big\{ \nu \in {}^{k(n+1)\geq }2 : &\ 
\nu \rest k^1(n) \in U^1_{\eta_k\rest n}, \text{ and for } \ell < 2^n \text{ we have }\\
&\ k^1(n) + \ell < \lh(\nu) \wedge (\ell \neq 2k+1) \Rightarrow \nu (k^1(n) + \ell) = 0 \big\}
\end{align*}
and $W_{\eta_k} = W_{\eta_k} \cup \{k^1(n)+2k+1\}.$
\mn
It is easy to verify that the construction provides a family of sets as
required. 
\end{PROOF}

\mn
\textbf{Proof of Lemma \ref{3.3}}:\quad Let $W(*)$, $U_\eta$, $W_\eta$
be as in \ref{3.3B}. Fix an enumeration $\{\eta_\xi:\xi<2^{\aleph_0}\}=
{}^\omega 2$ and let $W^\xi \defeq W_{\eta_\xi}$. Let
\[
S_\xi \subseteq \lim(U_{\eta_\xi})\ \big(= \big\{\rho \in {}^\omega 2 : 
(\forall n < \omega)[\rho\rest n \in U_{\eta_\xi}] \big\} \big)
\]

\mn
be of cardinality $\lambda$. Fix $\{\rho^\xi_i : i < \lambda\} = S_\xi$, 
and without loss of generality $S_\xi$ is $\chi$-large.\footnote{
    Recall that we say $S \subseteq {}^\omega 2$ is $\chi$-\emph{large} if for every $n < \omega$ and $\nu \in S$,\\ $\big|\{\rho\in S : \rho \rest n = \nu \rest n\} \big| \geq \chi$. If $\chi \geq (|\tau_1|+\aleph_0)^+$ we may omit it.
}

Note that for every $S\subseteq {}^\omega 2$ of cardinality $>\mu$, for some 
$S_1 \subseteq S$, $|S_1| \leq \mu$ and $S \setminus S_1$ is $\mu^+$-large. 
Let $U^\zeta = U_{\eta_\zeta}$; note that by \ref{3.3B}(B)+(D), the sets
$S_\xi \setminus S$ are pairwise disjoint.

So let $\xi,f,\bfn$ be as in the assumption of \ref{3.3}$\otimes$.

For $\nu \in S_\xi$ let $f(\nu) = {\bar\sigma}_\nu(\bar\eta_\nu)$, where
$\bar\sigma_\nu$ is a finite sequence of terms and $\bar\eta_\nu$ is a finite
sequence of members of $\bigcup\limits_{\zeta\neq\xi} S_\zeta$ with no
repetitions. So there are $S^* \subseteq S_\xi$ which is $\mu^+$-large,
and $\bar\sigma$, and an integer $n_0$ such that
\[
\nu\in S^*\ \Rightarrow\ \bar\sigma_\nu = \bar\sigma \wedge \lh(\bar \eta_\nu) = n_0,
\]

\mn
and without loss of generality, for some $m_a \leq m_b < \omega$, we have  
$\bar\sigma(\bar\eta_\nu)\rest m_a = \bar\eta^*$ and
\[
\{\eta^*_\ell : \ell < m_a\} \cup \{\eta_{\nu,\ell} : \nu \in S^*
\text{ and } \ell \in [m_a,m_b)\}
\]

\mn
is without repetition (this is possible by the $\Delta$-system argument).

As $S_\xi \cap \bigcup\limits_{\zeta\neq\xi} S_\zeta = \varnothing$, clearly 
the sequence $\LL\nu\RR \caret \bar\eta_{\bar\nu}$ is without repetitions for any
$\nu \in S^*$. So for some $k = k_\nu < \omega$ large enough, we have:
\mn
\begin{enumerate}
    \item[$(i)$] $\LL\nu\rest k\RR \caret \LL\eta_{\nu,\ell} \rest k : \ell < \ell(*)\RR$ is without repetitions.
\sn
    \item[$(ii)$] Letting $\eta_{\nu,\ell}\in S_{\zeta(\nu,\ell)}$, we have 
    $W^\xi \cap W^{\zeta(\eta_{\nu,\ell})} \subseteq \{0,\ldots,k_\nu-1\}$.
    Moreover, $k_\nu > \min(W^\xi \setminus W^{\zeta(\eta_\nu,\ell)})$ 
    (remember clause (F) of \ref{3.3}).
\end{enumerate}
\mn
As we can shrink $S^*$ as long as it is $\mu^+$-large, without loss of generality for
some $k$:
\begin{enumerate}
    \item[$(iii)$]  $\nu_1 \neq \nu_2 \in S^*\ \Rightarrow\ \lh(\nu_1 \cap \nu_2) > k$
\sn
    \item[$(iv)$]  $\nu \in S^*\ \Rightarrow\ k_\nu < k < \omega$.
\end{enumerate}
\sn
So for $\nu_1 \neq \nu_2 \in S^*$, on the one hand 
$\lh(\nu_1\cap\nu_2)\in W^\xi\setminus k$
(as $\nu_1,\nu_2\in S_\xi\subseteq\lim (U_{\eta_\xi})$; see clause (iii) above and
\ref{3.3B}(C)) and on the other hand
\[
\lh(\eta_{\nu_1\cap\ell},\eta_{\nu_2,\ell})\in W(*)\cup U^{\zeta
(\nu_1,\ell)}\cup U^{\zeta(\nu_2,\ell)}
\]

\mn
which is disjoint to $W^\xi\setminus k$. So we have proved clause (E) of
\ref{3.3}; the other clauses can be checked. 

\begin{claim}\label{3.7}
If clauses $(\beta)$, $(\gamma)$, $(\delta)$ of \ref{3.1A}(1) hold, and 
\ref{3.2}$(\otimes)$ does as well, \underline{then} for $\lambda \leq 2^{\aleph_0}$:
\mn
\begin{enumerate}
    \item[$(*)_\lambda$] There is a family $\cP$ of subsets of ${}^\omega2$ each 
    of cardinality $\lambda$ (even their union has cardinality $\lambda$) with 
    $|\cP| = 2^\lambda$, such that (letting $N^1_S$ be the Skolem Hull of 
    $\{\bar a_\eta : \eta\in S\}$ for $S \in \cP$) we have:
\sn
    \begin{itemize}
        \item  For $Y_1 \ne Y_2$ from $\cP$, $N^1_{Y_1}$ has no 
        $\Delta$-embedding into $N^1_{Y_2}$.
        
        \item $\|N^1_Y\| = \lambda$ for $Y \in \cP$. 
    \end{itemize}
\end{enumerate}
\end{claim}

\begin{PROOF}{\ref{3.7}}
For $X \subseteq \lambda$, let $M^1_X$ be the Skolem Hull of 
$\{\bar a_{\eta} : \eta \in \bigcup\limits_{\xi\in X} S_\xi\}$ and\\ 
$M_X \defeq M^1_X \rest \tau_T$. 

In order to prove the theorem it is enough to
assume $X,Y \subseteq \lambda$ and $X \not \subseteq Y$, and show there does not
exist an elementary embedding $f$ from $M_X$ into $M_Y$. Let $\xi\in
X\setminus Y$. For the sake of contradiction suppose $f:M_X 
\to M_Y$ is an elementary embedding, or just one 
preserving the satisfaction of $\varphi_n,\neg \varphi_n$.

We can represent $M_Y$ in 
$\cM_{\mu,\omega}\bigl(\bigcup\limits_{\zeta \neq \xi} S_\zeta\bigr)$, and let 
us define $f' : S_\xi \to \cM_{\mu,\omega}\bigl(\bigcup\limits_{\zeta\neq\xi} 
S_\zeta\bigr)$ by $f'(\nu) = f(\bar a_\nu)$, let $\bfn$ be essentially as in 
\ref{3.1A}, but translated. Apply lemma \ref{3.3} to $f'$ and $\bfn$, and get 
$S^*$, $k_0$, $n_0$, $m_a$, $  m_b$, $\bar\sigma$, 
$\LL\bar\eta_\nu : \nu \in S^* \RR$ as there. Of course $n_0$, $m_a$, $m_b$ 
are predetermined as in \ref{3.1A}.

So we are done proving \ref{3.7}.
\end{PROOF}

\begin{PROOF}{\ref{3.1}}
\textbf{Proof of Theorem \ref{3.1}}:

When $\lambda\leq 2^{\aleph_0}$, the result follows from \ref{3.7} by \ref{3.2}.

So the proof of Theorem \ref{3.1} for the case $\lambda\leq 2^{\aleph_0}$ is
completed. How to deal with the case $\lambda>2^{\aleph_0}$? We just need to
use $(\delta)^+$; i.e. use \ref{3.9} (and Definition \ref{3.8}) below.
\end{PROOF}

\begin{definition}\label{3.8}
For any cardinal $\kappa$ and $M_1$ as in \ref{3.1A}(1)$(\beta)$-$(\delta)^+$, 
we define a model $M_{1,\kappa}$ as follows: it is a $\tau_1$-model generated by 
$\{\bar a_\eta : \eta \in {}^\omega 2\} \cup \{d_i : i < \kappa\}$ such that 
for every $n < \omega$, $i_1 < \ldots i_n < \kappa$, and 
$\eta_1,\ldots,\eta_m \in {}^\omega 2$, the quantifier-free type of 
${\bar a}_{\eta_1} \caret \ldots \caret {\bar a}_{\eta_m} \caret \LL d_{i_1},\ldots,d_{i_n}\RR$ in $M_{1,\kappa}$ is equal to the
quantifier-free type of ${\bar a}_{\eta_1} \caret \ldots \caret{\bar a}_{\eta_m}
\caret \LL d_1,\ldots,d_n\RR$ in $M_1$. (So if $M_1$ has Skolem
functions then $M_1 = M_{1,\mu}$ and they realize the same types.)
\end{definition}

\begin{claim}\label{3.9}
If clauses $(\beta),(\gamma),(\delta),(\delta)^+$ of \ref{3.1A}(1) hold, and
\ref{3.2}$(\otimes)$ does as well, \underline{then} for $\lambda\geq 2^{\aleph_0}$:
\mn
\begin{enumerate}
    \item[$(*)_\lambda$]  There is a family $\cP$ of subsets of ${}^\omega2$ each of cardinality $2^{\aleph_0}$ with $|\cP| = \beth_2$ such that, letting 
    $N^\lambda_S$ be the Skolem Hull of $\{\bar a_\eta : \eta \in S\} \cup 
    \{d_i : i < \kappa\}$ in $M_{1,\lambda}$ with $S \in \cP$ 
    (so $\|N^\lambda_S\| = \lambda$), we have:
\sn
    \begin{enumerate}
        \item[$(*)$]  For $Y_1\neq Y_2$ from $\cP$, $N^1_{Y_1}$ has no $\Delta$-embedding into $N^1_{Y_2}$ (i.e. no function from $N^1_{Y_1}$ into $N^1_{Y_2}$ preserves all the relations $\pm \varphi_n)$.
    \end{enumerate}
\end{enumerate}
\end{claim}

\noindent
We may consider using relations {$\varphi_n$} which are not equivalence
relations, and we may like to give another proof when $\mu > \aleph_0$ but
still $\MA_\mu$ holds.

\begin{claim}\label{3.10}
$[$Assume $\MA_\mu.]$ 

Suppose $M_1$, $\tau_1$, $\LL\bar a_\eta : \eta \in {}^\omega 2\RR$, $\varphi_n$ 
(for $n < \omega$), $\LL d_n : n < \omega\RR$ satisfy clauses
$(a)$, $(b)$, $(\beta)$, $(\gamma)$, $(\delta)$ of \ref{3.1A}, 
and $M_1$ is a $\tau_1$-model of the complete first order theory $T_1$. 
Also suppose $\bar a_\eta\in {}^{k}(M_1)$ for $\eta\in {}^{\omega>}2$ are
such that if $n < m < \omega$ and $\eta,\nu\in{}^m2$ then 
$\eta\rest n = \nu \rest n \Leftrightarrow M_1\models \bar a_\eta\ \varphi_n\ \bar a_\nu$.
(So $\varphi_n$ is not necessarily an equivalence relation and $|\tau_1| = \mu$ 
is not necessary countable).

\sn
$1)$  If we replaced $^{\omega\geq} 2$ by a perfect subtree (splitting
determined by level only) and replaced $\LL \varphi_n : n < \omega\RR$
by a subsequence, then we could add the statement of
\ref{3.2}$(\otimes)$ to the assumptions.

\sn
$2)$  So the conclusion of \ref{3.7} holds, and if we further assume
$(\delta)^+$ of \ref{3.1A}, the conclusion of \ref{3.9} also holds.
\end{claim}

\begin{PROOF}{\ref{3.10}}
We use Carlson and Simpson \cite{CS84}.

Let $W^*$ be the set of $\omega$-sequences $\eta$ from 
$\{0,1\} \cup \{x_i : i < \omega\}$ such that each $x_i$ appears infinitely often. 
For $\eta \in W^*$, let 
$$W_\eta = \{\nu\in W^*: \eta(\ell)\in\{0,1\} \Rightarrow \nu(\ell)=\eta(\ell),\ \eta(\ell_1) = \eta(\ell_2) \Rightarrow \nu(\ell_1) = \nu(\ell_2)\}.$$ 
As a set, $W \subseteq W^*$ is \emph{large} if it contains some $W_\eta$.
Let 
\begin{align*}
    I_W = \big\{\nu\in {}^{\omega>}2 : & \text{ for some $\eta\in W$, for every }
    \ell,\ell_1,\ell_2 < \lh(\nu),\\
    &\ \eta(\ell)\in\{0,1\} \Rightarrow \nu(\ell)=\eta(\ell) \wedge 
    \eta(\ell_1) = \eta(\ell_2) \Rightarrow \nu(\ell_1) = \nu(\ell_2) \big\}.
\end{align*}
Let 
$$\lev(W) = \big\{ \ell : \text{for some }\eta\in W,\ \eta(\ell)\notin \{0,1\}
\text{ but } \eta(0),\ldots,\eta(\ell-1)\in \{0,1\} \big\}.$$ 

{We say} $W_1 \subseteq^* W_2$ if for some $n$, 
$\{\nu\rest[n,\omega) : \nu\in W_1\} \supseteq \{\nu \rest (n,\omega) : \nu \in W_2\}$. 
By $\MA_\mu$, if $\LL W_i : i < \delta \leq \mu\RR$ is
$\subseteq^*$-decreasing sequence then there is $W$ such that
$\bigwedge\limits_i W_i \subseteq^* W$.

By the partition theorem {there}, if $n<\omega$, 
$\eta_1,\ldots,\eta_k \in {}^n 2$ are pairwise distinct and ${\bar\sigma}^1$, 
${\bar\sigma}^2$ are $\tau^+_1$-terms \underline{then} we can 
find large $W_1 \subseteq W$ such that $W_1 \rest n = W \rest n$ and:
\mn
\begin{enumerate}
    \item[$\circledast^n_{W_1,\bar \sigma}$]  If $n < m \in \lev(W_1)$, 
    ${}\rho^\nu_\ell \in \clT_{W_1} \cap {}^m 2$ for $\ell=1,\ldots,k$, and 
    $\nu_\ell=\eta_\ell \caret \rho_2\rest [n,\omega)$, \underline{then} 
    the truth value of 
    ${\bar\sigma}^1(\bar a_{\nu_1},\ldots,\bar a_{\nu_k})\ \varphi_n\  {\bar\sigma}^2(\bar a_{\nu_1},\ldots,\bar a_{\nu_k})$ is constant.
\end{enumerate}
\mn
Repeating it, we can get $W_1$ such that $\circledast^n_{W_1,\bar \sigma}$ 
for every $n$.
\mn
\begin{enumerate}
    \item[$(i)$]  Either $g$ is constant $< \min(\lev(W_1) \setminus n)$ or
    \[
        n \in \lev (W_1)\ \Rightarrow\ \big[ g(n),n \big) \cap \lev(W_1) = \varnothing.
    \]

    \item[$(ii)$] If $n < m \in \lev(W_1)$ and 
    $\eta_\ell \lhd \nu_\ell \in \clT_{W_1}\cap {}^m 2$ then
    \[
        \min\!\big\{i : \neg\big[\sigma^1(\bar a_{\nu_1},\ldots,\bar a_{\nu_k})\ \varphi_i\ 
        \sigma^2(\bar a_{\nu_1},\ldots,\bar a_{\nu_k})\big]\big\} = g(m).
    \]
\end{enumerate}
\mn
We apply such reasoning to the following statement: ``Given 
$\eta_1,\ldots,\eta_k \in \clT_{W_1} \cap {}^n 2$ pairwise distinct {and} 
$n < m \in \lev(W_1)$, and assuming $\eta_\ell \lhd \nu_\ell^i \in \clT_{W_1}\cap {}^m2$ for $\ell \in \{0,1,\ldots,k\}$ and $i \in \{0,1\}$, we have
$$\bar\sigma(\bar a_{\nu^0_1},\ldots,\bar a_{\nu^0_k})\ \varphi_\ell\ \bar\sigma (\bar a_{\nu^1_1},\ldots,\bar a_{\nu^1_k})."$$ 
We get that this depends only on $\lh(\nu^0_\ell \cap \nu^1_\ell)$ and $\nu^i_{\ell}\big(\lh(\nu^0_\ell\cap\nu^1_\ell) \big)$.
\end{PROOF}

\begin{discussion}\label{5.14}
The parallel (for a module $\dot \bbM$) concerning ``a surgery at" is extending
the ring $\dot{\bfR}$ to $\dot{\bfR}^+$; e.g. by $\{x_t : t\in I\}$ freely except 
some equation involving $x$ and the $x_i$-s and ``below $x$" is replaced by the
ideal generated by those equations. 
\end{discussion}


\bibliographystyle{amsalpha}
\bibliography{shlhetal}

\end{document}